\documentclass{elsarticle}
\usepackage[utf8]{inputenc}

\usepackage{graphicx}
\usepackage{url}
\usepackage[margin=1in]{geometry}
\usepackage{subcaption}
\usepackage{amsmath}
\usepackage{amssymb}
\newtheorem{remark}{Remark}
\usepackage[version=4]{mhchem}
\usepackage{siunitx}
\usepackage{longtable,tabularx}
\usepackage{booktabs}
\usepackage{xcolor}
\usepackage{float}
\usepackage{pgfplots}
\pgfplotsset{compat=1.9}
\usepgfplotslibrary{groupplots}
\usepgfplotslibrary{fillbetween}
\usetikzlibrary{arrows.meta}
\usetikzlibrary{matrix}

\usetikzlibrary{external}
\tikzifexternalizing{\PassOptionsToPackage{bookmarks=false}{hyperref}}{}

\pgfplotsset{
    colormap={fast}{
        rgb=(0.056,0.056,0.470)
        rgb=(0.243,0.460,0.810)
        rgb=(0.357,0.745,0.954)
        rgb=(0.688,0.930,0.918)
        rgb=(0.899,0.945,0.769)
        rgb=(0.957,0.834,0.509)
        rgb=(0.928,0.621,0.315)
        rgb=(0.800,0.352,0.160)
        rgb=(0.590,0.077,0.119)
    }
}

\newcommand{\logLogSlopeTriangle}[6][north]{%
  \pgfplotsextra{%
    \pgfkeysgetvalue{/pgfplots/xmin}{\xmin}%
    \pgfkeysgetvalue{/pgfplots/xmax}{\xmax}%
    \pgfkeysgetvalue{/pgfplots/ymin}{\ymin}%
    \pgfkeysgetvalue{/pgfplots/ymax}{\ymax}%
    \pgfmathsetmacro{\xArel}{#2}%
    \pgfmathsetmacro{\yArel}{#4}%
    \pgfmathsetmacro{\xBrel}{#2-#3}%
    \pgfmathsetmacro{\yBrel}{\yArel}%
    \pgfmathsetmacro{\xCrel}{\xArel}%
    \pgfmathsetmacro{\lnxB}{\xmin*(1-(#2-#3))+\xmax*(#2-#3)}%
    \pgfmathsetmacro{\lnxA}{\xmin*(1-#2)+\xmax*#2}%
    \pgfmathsetmacro{\lnyA}{\ymin*(1-#4)+\ymax*#4}%
    \pgfmathsetmacro{\lnyC}{\lnyA+#5*(\lnxA-\lnxB)}%
    \pgfmathsetmacro{\yCrel}{\lnyC-\ymin)/(\ymax-\ymin)}%
    \coordinate (A) at (rel axis cs:\xArel,\yArel);%
    \coordinate (B) at (rel axis cs:\xBrel,\yBrel);%
    \coordinate (C) at (rel axis cs:\xCrel,\yCrel);%
    \draw[#6] (A) -- node[pos=0.5,anchor=#1,scale=0.8] {1}
              (B) --
              (C) -- node[pos=0.5,anchor=west,scale=0.8] {$#5$}
              cycle;%
  }%
}
\newcommand{\func}[3]{\ensuremath{#1 \colon #2 \to #3}}
\newcommand{\Rbb}{\ensuremath{\mathbb{R}}}
\newcommand{\case}[1]{\texttt{#1}}

\pgfplotsset{
    every axis title/.append style={font=\small},
    tick label style={font=\normalsize},
    label style={font=\normalsize},
    legend style={font=\normalsize},
}
\tikzset{plabel/.style={anchor=north west, font=\small, inner sep=1pt, fill=white, fill opacity=0.75, text opacity=1}}
\newcommand{\figurelegend}[1]{\par\vspace{0.4em}\noindent\makebox[\linewidth][c]{\ref*{#1}}\par}

\usepackage[colorlinks, citecolor=blue, linkcolor=blue, urlcolor=blue]{hyperref}

\begin{document}

\begin{frontmatter}

\title{A Principled Framework for PDE-Based Shock Capturing in High-Order Discontinuous Galerkin Methods}

\author[nd]{Aaron Marques}
\ead{rmarque2@nd.edu}
\author[ut]{Alireza H.\ Razavi}
\ead{a.razavi@mail.utoronto.ca}
\author[ut]{Masayuki Yano}
\ead{masa.yano@utoronto.ca}
\author[nd]{Matthew J.\ Zahr\corref{cor1}}
\ead{mzahr@nd.edu}
\affiliation[nd]{organization={Department of Aerospace and Mechanical Engineering, University of Notre Dame},
            addressline={Fitzpatrick Hall of Engineering},
            city={Notre Dame},
            state={Indiana},
            postcode={46556},
            country={United States}}
\affiliation[ut]{organization={Institute for Aerospace Studies, University of Toronto},
            addressline={4925 Dufferin Street},
            city={Toronto},
            state={Ontario},
            postcode={M3H 5T6},
            country={Canada}}
\cortext[cor1]{Corresponding author}

\begin{abstract}
PDE-based artificial viscosity (AV) captures shocks in high-order discontinuous Galerkin (DG) methods by solving an auxiliary PDE for a smooth viscosity field. The behavior of that field is set by how the terms of the AV PDE scale with the element size $h$ and polynomial degree $p$, and by how its coefficients are parametrized. This work derives both choices from analysis and studies the impact these choices have on the implicit steady-state solver. First, an $L^2$-projection analysis shows that the thinnest shock a polynomial space can represent monotonically scales with the resolution length $\bar h = h/p$, and a scaling analysis of the AV PDE gives concrete conditions on its diffusion and source terms for the solution to be self-similar in $\bar h$. Prior formulations satisfy them only under $h$-refinement at fixed $p$; the AV PDE is scaled so that they hold under $hp$-refinement, and the captured shock thins accordingly. Second, a Green's-function analysis of the steady AV PDE gives closed-form expressions for the width and magnitude of the AV. In the standard parametrization the two are coupled: raising the diffusion coefficient widens the AV but lowers its peak. An alternative parametrization of the decay coefficient removes the coupling in the Green's-function response and inverts it in the coupled state--AV problem, so that widening the AV raises its peak instead of lowering it. Third, a parameter study on a Mach~5 viscous flow over a half-cylinder characterizes how the choice of solution variables and the AV PDE's time, width, and magnitude coefficients affect the robustness and cost of the implicit solve. The choice of variables dominates robustness, with entropy variables converging over the largest region of the width--magnitude parameter space. The time constant trades robustness against cost, measured in pseudo-transient and GMRES iterations, and conservative variables converge only within a narrow window of the time constant. The scaling and parametrization are validated on the 1D Burgers' equation, Mach~5 inviscid and viscous flow over a half-cylinder, Mach~11.68 viscous flow over a double wedge, and Mach~5 viscous flow over a sphere. Together, the three contributions make PDE-based AV more predictable to tune and its implicit solve more robust.
\end{abstract}


\end{frontmatter}

\section{Introduction}
\label{sec:intro}

Shock-dominated flows arise across aerospace applications, from transonic aerodynamics to hypersonic re-entry, and remain difficult to simulate accurately. High-order discontinuous Galerkin (DG) methods~\cite{hesthaven2007nodal} are attractive for such simulations because they combine high accuracy per degree of freedom with low numerical dissipation, geometric flexibility on unstructured meshes, and parallel scalability. Near shocks and contact discontinuities, however, the polynomial approximation admits spurious oscillations that erode solution quality and can produce negative pressure or density that crash the solver. A robust strategy for handling shocks is essential for these methods to be competitive in this flow regime.

Two broad strategies address this challenge. Shock \emph{capturing} regularizes discontinuities with numerical dissipation so that the discrete solution can represent them on a fixed mesh. Classical finite-volume shock capturing builds the dissipation into the reconstruction, from monotone second-order schemes~\cite{van1979towards} to essentially non-oscillatory and weighted essentially non-oscillatory reconstructions~\cite{harten1997uniformly, liu1994weighted, jiang1996efficient}. In the DG setting, shock capturing takes the form of slope limiters~\cite{cockburn2000development}, a posteriori sub-cell limiting~\cite{dumbser2014posteriori, vilar2019posteriori}, or artificial viscosity, which adds an explicit dissipative term scaled to the sub-cell resolution of the high-order approximation~\cite{persson2006sub, barter2010shock, persson2013shock, moro2016dilation, fernandez2018physics, ching2019shock}. Captured shocks are spread over a few elements and the formal order of accuracy is lost in their vicinity, but the approach is robust and applies to shock systems of arbitrary complexity. Shock \emph{fitting} instead aligns mesh faces with discontinuities so that the approximation never represents a jump. Classical shock fitting treats each shock as an explicit boundary advanced by the Rankine--Hugoniot relations~\cite{moretti2002thirty, salas2009shock}, and unstructured variants extend the idea to interacting shocks and shock--boundary-layer interaction~\cite{bonfiglioli2016unsteady, zou2021embedded, assonitis2021numerical}. Implicit shock tracking recasts fitting as an optimization problem in which the mesh and solution are determined together, so that faces align with discontinuities without prior knowledge of their location~\cite{corrigan2019moving, kercher2021least, zahr2018shktrk, zahr2020implicit, huang2022robust}, and it recovers high-order accuracy for inviscid shock-dominated flows, including unsteady ones~\cite{shi2022implicit, naudet2024space}.

This work pursues artificial-viscosity (AV) shock capturing~\cite{persson2006sub, barter2010shock}, and the central question is how the AV field should be defined over the domain. A natural starting point is to derive it directly from the finite-dimensional solution, either as a pointwise function of the flow state or as a single value per element. In the DG setting, the solution itself is discontinuous at element interfaces, so an AV field defined directly from it inherits those discontinuities. Persson~\cite{persson2013shock} and Barter~\cite{barter2010shock} argued that the AV should nevertheless be continuous, because a discontinuous AV induces oscillations in the state gradients that corrupt the downstream flow. Recovering a continuous AV through this route would require discretizing it on a continuous ($H^1$-conforming) space, which sacrifices the element-local structure of DG methods.

PDE-based AV avoids this trade-off. Barter and Darmofal~\cite{barter2010shock} append an auxiliary diffusive PDE to the system of governing equations: the AV field becomes the solution of its own equation, with a shock sensor as the source and a diffusion term that smooths the field. The auxiliary equation is discretized on the same DG mesh as the conservation law, so smoothness is obtained within the DG framework. A later formulation~\cite{yano2012optimization} replaces the bounding-box diffusion tensor with a mesh-implied metric tensor, which propagates the AV consistently under arbitrary element orientation and sharpens shock capture on anisotropic meshes. The present work adopts this PDE-based AV formulation as the starting point for the design choices developed in Sections~\ref{sec:avpde_design} and~\ref{sec:reparam}.

PDE-based AV formulations introduce a small set of free parameters that together control the AV magnitude, width, propagation, decay, and boundary behavior. The formulations of~\cite{barter2010shock, yano2012optimization} adopt a single fixed set of constants across all the configurations they consider. More recently, Nguyen et al.~\cite{nguyen2024optimal} automate the choice of an amplitude and a width parameter by constrained optimization, minimizing the artificial viscosity subject to positivity of density and pressure and to a smoothness constraint on the numerical solution, in combination with mesh adaptation. Razavi et al.~\cite{razavi2026output} combine PDE-based AV with output-based error estimation and anisotropic mesh adaptation, and study how the diffusion coefficient of the AV PDE affects the adapted mesh.

Three aspects of this framework have received less attention, and this work addresses each. First, prior formulations~\cite{barter2010shock, yano2012optimization, nguyen2024optimal} tie the length scale of the AV PDE to the element size $h$ alone, so the AV's width does not shrink under $p$-refinement even though the enriched polynomial space can then resolve thinner profiles. We use $L^2$ projection of a 1D model shock onto a polynomial element to identify the smallest thickness the basis can monotonically resolve as a function of $h$, $p$, and shock position, and find that it scales with $\bar h = h/p$. A scaling analysis of the AV PDE then gives concrete conditions on the powers of $\bar h$ carried by its diffusion and source terms for the solution to be self-similar in $\bar h$. We scale the AV PDE so that these conditions hold under both $h$- and $p$-refinement, and the AV and the captured shock then thin under $hp$-refinement.

Second, in the standard parametrization shared by these formulations, increasing the diffusion parameter widens the AV but also reduces its magnitude, which can leave oscillations in the captured solution. A Green's-function analysis of the steady AV PDE gives closed-form expressions for the width and magnitude of the AV and motivates a reparametrization of the decay coefficient under which widening no longer suppresses the magnitude in the Green's-function response and raises it in the coupled state--AV problem.

Third, the effect of the solver state variables, such as conservative, entropy, and logarithmic variables, on convergence has been studied for DG methods~\cite{alberti2024comparative}, but not in the presence of PDE-based AV, whose own solver parameters add to the question. A systematic study on a Mach~5 viscous flow over a half-cylinder characterizes how the variable formulation and the AV PDE parameters $c_1$, $c_2$, and $c_3$ affect the robustness and cost of the implicit solve. Variable choice has a substantial effect on robustness, with entropy variables converging on the largest fraction of the $(c_2, c_3)$ space tested. The time constant $c_1$ trades robustness against solver cost, measured in pseudo-transient and GMRES iterations, with the conservative formulation converging only in a narrow $c_1$ window.

The remainder of the paper is organized as follows. Section~\ref{sec:govern} presents the governing equations, the AV PDE, the resulting augmented system, and the DG discretization. Section~\ref{sec:avpde_design} develops the resolution-aware scaling of the AV PDE and verifies it on the 1D Burgers' equation and a Mach~5 inviscid flow over a half-cylinder. Section~\ref{sec:reparam} develops the reparametrization of the decay coefficient and verifies it on the same problems and on a Mach~5 viscous flow over a half-cylinder. Section~\ref{sec:applications} applies the complete method to a Mach~11.68 viscous flow over a double wedge and a Mach~5 viscous flow over a sphere. Section~\ref{sec:var_robustness} turns to the implicit solver, examining how the variable formulation and the time constant $c_1$ govern its robustness and cost. Section~\ref{sec:conclude} summarizes the contributions and outlines future work.

\section{Governing equations and high-order discretization}
\label{sec:govern}

This section defines the discrete problem solved in the rest of the paper. Starting from a system of conservation laws (Section~\ref{subsec:claws}), we append the artificial-viscosity PDE, whose solution supplies a smooth viscosity field to the conservation law it regularizes (Section~\ref{subsec:scpde}), and combine the two into a single augmented system (Section~\ref{subsec:augsys}) with boundary conditions for both (Section~\ref{subsec:av_bnd}). Section~\ref{subsec:claws_used} specializes the augmented system to the two physical models used in this work, Burgers' equation and the compressible Navier--Stokes equations, and introduces the variable formulations in which the latter is solved. Section~\ref{subsec:dg} discretizes the augmented system with an interior-penalty DG method and solves it by pseudo-transient continuation. The scaling of the AV PDE and the choice of its decay coefficient are stated here and justified in Sections~\ref{sec:avpde_design} and~\ref{sec:reparam}.

\subsection{System of conservation laws}
\label{subsec:claws}

Consider a general system of $m$ conservation laws on the spatial domain $\Omega \subset \Rbb^d$ and temporal domain $\mathcal{T} \subset \Rbb_{\geq 0}$,
\begin{equation}
\frac{\partial u}{\partial t} + \nabla \cdot \left[ F^{\mathrm{I}}(u) - F^{\mathrm{V}}(u, \nabla u) \right] = S(u, \nabla u) \quad \text{in} \quad \Omega \times \mathcal{T},
\label{eq:claws}
\end{equation}
where $\func{u}{\Omega \times \mathcal{T}}{\Rbb^m}$ is the state, $\func{F^{\mathrm{I}}}{\Rbb^m}{\Rbb^{m \times d}}$ the inviscid flux, $\func{F^{\mathrm{V}}}{\Rbb^m \times \Rbb^{m \times d}}{\Rbb^{m \times d}}$ the viscous flux, and $\func{S}{\Rbb^m \times \Rbb^{m \times d}}{\Rbb^m}$ the source. Throughout, we seek steady solutions of~\eqref{eq:claws} by pseudo-transient continuation, so the time term is a pseudo-time march toward the steady state and is not used to capture transient physics.

\subsection{Shock-capturing PDE and modified system of conservation laws}
\label{subsec:scpde}

To stabilize the solution near shocks, we introduce an AV working variable $\func{\hat{\nu}}{\Omega \times \mathcal{T}}{\Rbb}$ governed by the auxiliary PDE of Barter and Darmofal~\cite{barter2010shock} with anisotropic length scaling~\cite{yano2012optimization}, recast as
\begin{equation}
    c_1 \, \frac{\partial \hat{\nu}}{\partial t} - c_2 \, \nabla \cdot \left[\eta \, \nabla \hat{\nu}\right] = c_3 \, \bar{s}(u, \nabla u) - c_4 \, \hat{\nu} \quad \text{in} \quad \Omega \times \mathcal{T},
    \label{eq:avpde}
\end{equation}
where $\func{\bar{s}}{\Rbb^m \times \Rbb^{m \times d}}{\Rbb}$ is a shock sensor, specified for each model in Section~\ref{subsec:claws_used}, $\func{\eta}{\Omega}{\Rbb^{d \times d}}$ is a diffusion tensor, and $c_1$, $c_2$, $c_3$, $c_4$ are constants. The diffusion tensor is built from the element geometry. Let $\func{G}{\Omega}{\Rbb^{d \times d}}$ be the Jacobian of the map from the physical element $\Omega_e$ containing the point $x$ to an ideal element $\Omega_\star$ of unit size, and let $X$ be the coordinate in the ideal element, so that the infinitesimal line element satisfies $dX = G\, dx$, and define pointwise $H = (G^T G)^{-1/2}$. For an element that is affinely mapped to the ideal element, $G$ is constant and $\partial/\partial X = G^{-T}\, \partial/\partial x$, so unit isotropic diffusion in the ideal element pulls back to physical space as
\begin{equation}
    \nabla_X \cdot \nabla_X \hat\nu = \nabla \cdot \big( G^{-1} G^{-T} \nabla \hat\nu \big) = \nabla \cdot \big( H^2 \, \nabla \hat\nu \big).
    \label{eq:H_pullback}
\end{equation}
When the map is not affine, as on a curved element, $G$ varies over the element and the pullback is not exact, but we use the same $H$ pointwise. The eigenvectors of $H$ are the principal directions of the element and its eigenvalues $h_i$ the corresponding element lengths, so $H^2 = h^2 I$ on an isotropic element of size $h$, and on an elongated element $H^2$ diffuses in proportion to the squared length along each principal direction. The diffusion tensor of the AV PDE is this geometric tensor divided by the square of the polynomial degree,
\begin{equation}
    \eta = \frac{H^2}{p^2},
    \label{eq:eta_tensor}
\end{equation}
so that the diffusion length is the resolution length $h/p$ rather than the element size $h$; the division by $p^2$ is what makes the AV self-similar under $hp$-refinement, as derived in Section~\ref{subsubsec:hbar2}. The constant $c_1$ is a time constant that affects only the pseudo-transient iteration (Section~\ref{subsec:c1_study}); $c_2$ and $c_3$ set the width and magnitude of the AV; and the decay coefficient is tied to the diffusion coefficient as $c_4 = 1/c_2$ (Section~\ref{subsec:c4_reparam}), except where the standard choice $c_4 = 1$ of prior formulations is compared.

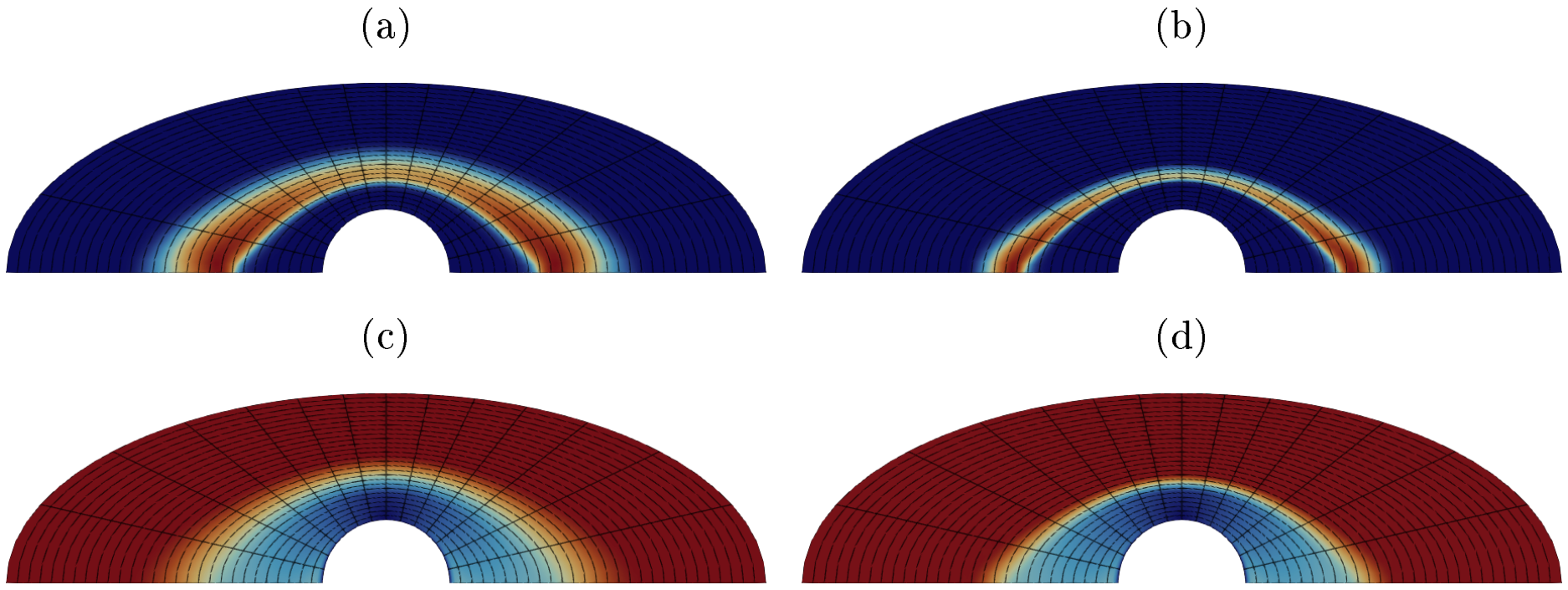
\begin{figure}[hbt!]
    \centering
    \begin{tikzpicture}
        \tikzset{cell/.style={inner sep=2pt, outer sep=0pt, anchor=north west}}
        \node[cell] (A) at (0,0) {%
            \IfFileExists{_tikz/av_iso.tex}{\begin{tikzpicture}[baseline=(current axis.north)]
\begin{axis}[title={(a)},
    tick label style={font=\normalsize},
    width=6cm,
    axis equal image,
    scale only axis,
    hide axis,
    xmin=0, xmax=902, ymin=0, ymax=225,
    enlargelimits=false,
    colormap name=fast,
    point meta min=0,
    point meta max=0.681209,
]
\addplot graphics[xmin=0, xmax=902, ymin=0, ymax=225] {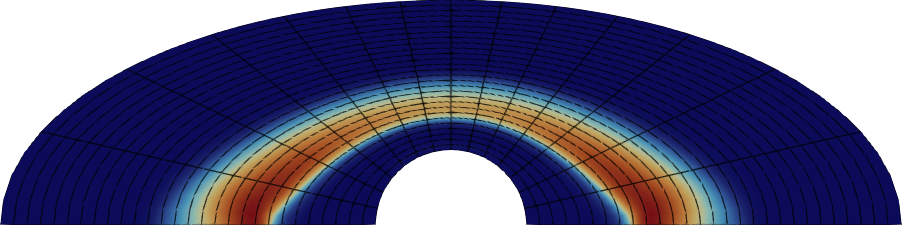};
\end{axis}
\path[use as bounding box] (current axis.south west) rectangle (current axis.north east);
\end{tikzpicture}}{\fbox{\parbox[c][3cm][c]{0.45\linewidth}{\centering\scriptsize AV, scalar diffusion}}}};
        \node[cell] (B) at ([xshift=0.15cm]A.north east) {%
            \IfFileExists{_tikz/av_aniso.tex}{\begin{tikzpicture}[baseline=(current axis.north)]
\begin{axis}[title={(b)},
    tick label style={font=\normalsize},
    width=6cm,
    axis equal image,
    scale only axis,
    hide axis,
    xmin=0, xmax=902, ymin=0, ymax=225,
    enlargelimits=false,
    colormap name=fast,
    point meta min=0,
    point meta max=0.365219,
]
\addplot graphics[xmin=0, xmax=902, ymin=0, ymax=225] {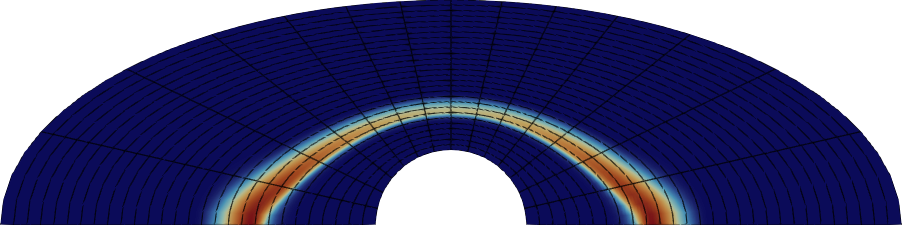};
\end{axis}
\path[use as bounding box] (current axis.south west) rectangle (current axis.north east);
\end{tikzpicture}}{\fbox{\parbox[c][3cm][c]{0.45\linewidth}{\centering\scriptsize AV, diffusion tensor}}}};
        \node[cell] (C) at ([yshift=-0.15cm]A.south west) {%
            \IfFileExists{_tikz/mach_iso.tex}{\begin{tikzpicture}[baseline=(current axis.north)]
\begin{axis}[title={(c)},
    tick label style={font=\normalsize},
    width=6cm,
    axis equal image,
    scale only axis,
    hide axis,
    xmin=0, xmax=902, ymin=0, ymax=225,
    enlargelimits=false,
    colormap name=fast,
    point meta min=0.0028880,
    point meta max=5.00755,
]
\addplot graphics[xmin=0, xmax=902, ymin=0, ymax=225] {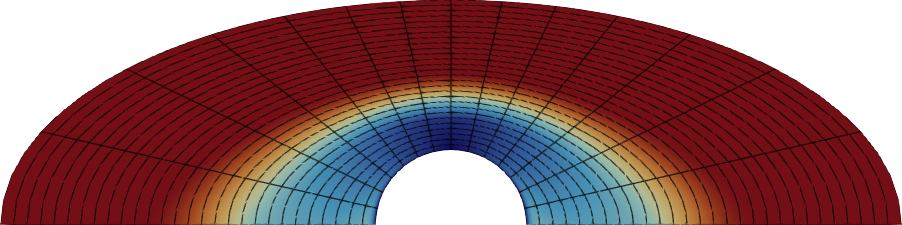};
\end{axis}
\path[use as bounding box] (current axis.south west) rectangle (current axis.north east);
\end{tikzpicture}}{\fbox{\parbox[c][3cm][c]{0.45\linewidth}{\centering\scriptsize Mach, scalar diffusion}}}};
        \node[cell] (D) at ([xshift=0.15cm]C.north east) {%
            \IfFileExists{_tikz/mach_aniso.tex}{\begin{tikzpicture}[baseline=(current axis.north)]
\begin{axis}[title={(d)},
    tick label style={font=\normalsize},
    width=6cm,
    axis equal image,
    scale only axis,
    hide axis,
    xmin=0, xmax=902, ymin=0, ymax=225,
    enlargelimits=false,
    colormap name=fast,
    point meta min=0.003279,
    point meta max=5.03928,
]
\addplot graphics[xmin=0, xmax=902, ymin=0, ymax=225] {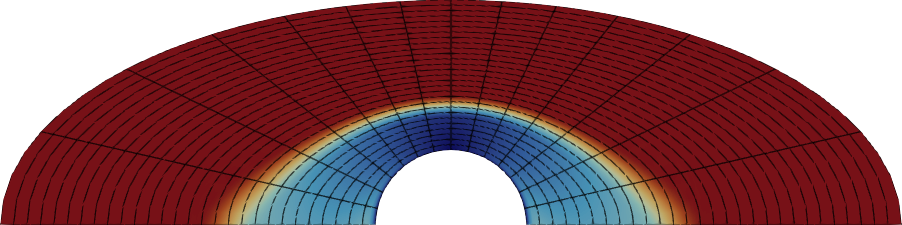};
\end{axis}
\path[use as bounding box] (current axis.south west) rectangle (current axis.north east);
\end{tikzpicture}}{\fbox{\parbox[c][3cm][c]{0.45\linewidth}{\centering\scriptsize Mach, diffusion tensor}}}};
    \end{tikzpicture}
    \caption{Mach~5 viscous flow over a half-cylinder (Section~\ref{subsec:c4_practice}): (a,b)~artificial viscosity $\nu$; (c,d)~Mach number, with (a,c)~the scalar diffusion $\eta = (h_V/p)^2 I$ and (b,d)~the diffusion tensor~\eqref{eq:eta_tensor}.}
    \label{fig:tensor_comparison}
\end{figure}

\begin{remark}
On elements with large aspect ratio, in particular on shock-aligned meshes that are thin normal to the shock and long along it, a scalar diffusion $\eta = (h_V/p)^2 I$, with $h_V$ an isotropic size from the element volume, cannot follow the element shape. On a Mach~5 viscous flow over a half-cylinder, the scalar diffusion spreads the AV isotropically around the bow shock and smears it, whereas the tensor~\eqref{eq:eta_tensor} confines the AV to a thin band aligned with the shock (Fig.~\ref{fig:tensor_comparison}). The bounding-box tensor $\mathrm{diag}(h_1^2, \ldots, h_d^2)$ of Barter and Darmofal~\cite{barter2010shock} is coordinate dependent and over-diffuses shocks that cross elements diagonally, an artifact documented on tilted meshes in~\cite{yano2012optimization} and removed by the metric tensor.
\end{remark}

The AV that enters the conservation law is the filtered working variable $\nu = \Theta(\hat{\nu})$, where $\func{\Theta}{\Rbb}{\Rbb_{\geq 0}}$ is the $C^1$ approximation of $\max(\hat{\nu} - \nu_0, 0)$
\begin{equation}
    \Theta(\hat{\nu}) = \begin{cases}
        0 & \hat{\nu} \le \nu_0, \\
        z^2 \, (2a - z) / a^2 & \nu_0 < \hat{\nu} < \nu_0 + a, \\
        z & \hat{\nu} \ge \nu_0 + a,
    \end{cases}
    \label{eq:filter}
\end{equation}
where $z = \hat{\nu} - \nu_0$, with transition width $a > 0$, which is small relative to the AV in shocks, and shift $\nu_0 \geq 0$. The values of $a$ and $\nu_0$ are given with each test case. The AV enters the conservation law~\eqref{eq:claws} through an artificial viscous flux $\func{F^{\mathrm{AV}}}{\Rbb^m \times \Rbb^{m \times d} \times \Rbb}{\Rbb^{m \times d}}$, specified for each model in Section~\ref{subsec:claws_used}, that adds to the physical viscous flux. The modified conservation law reads
\begin{equation}
    \frac{\partial u}{\partial t} + \nabla \cdot \left[ F^{\mathrm{I}}(u) - F^{\mathrm{V}}(u, \nabla u) - F^{\mathrm{AV}}(u, \nabla u; \nu) \right] = S(u, \nabla u) \quad \text{in} \quad \Omega \times \mathcal{T}.
    \label{eq:modclaws}
\end{equation}

\begin{remark}
The filter~\eqref{eq:filter} guarantees $\nu \geq 0$: a negative AV would act as an antidiffusive viscous flux. The $C^1$ cubic, rather than $\max(\hat\nu - \nu_0, 0)$ itself, keeps the augmented system differentiable for the implicit solver of Section~\ref{subsec:dg}.
\end{remark}

\subsection{Augmented system}
\label{subsec:augsys}

The modified conservation law~\eqref{eq:modclaws} and the AV PDE~\eqref{eq:avpde}, divided by $c_1$, form a coupled system for the augmented state $\func{U = (u, \hat{\nu})}{\Omega \times \mathcal{T}}{\Rbb^{m+1}}$,
\begin{equation}
    \frac{\partial U}{\partial t} + \nabla \cdot \left[ \mathcal{F}^{\mathrm{I}}(U) - \mathcal{F}^{\mathrm{V}}(U, \nabla U) \right] = \mathcal{S}(U, \nabla U) \quad \text{in} \quad \Omega \times \mathcal{T},
    \label{eq:augmented}
\end{equation}
with augmented fluxes and source
\begin{align}
    \mathcal{F}^{\mathrm{I}}(U) &= \begin{pmatrix} F^{\mathrm{I}}(u) \\ 0 \end{pmatrix}, \qquad
    \mathcal{F}^{\mathrm{V}}(U, \nabla U) = \begin{pmatrix} F^{\mathrm{V}}(u, \nabla u) + F^{\mathrm{AV}}(u, \nabla u; \Theta(\hat{\nu})) \\ \tfrac{c_2}{c_1} \, (\eta \, \nabla \hat{\nu})^T \end{pmatrix}, \notag \\
    \mathcal{S}(U, \nabla U) &= \begin{pmatrix} S(u, \nabla u) \\ \tfrac{c_3}{c_1} \, \bar{s}(u, \nabla u) - \tfrac{c_4}{c_1} \, \hat{\nu} \end{pmatrix}.
    \label{eq:augmented_fluxes}
\end{align}

\subsection{Boundary conditions}
\label{subsec:av_bnd}

The augmented system~\eqref{eq:augmented} requires boundary conditions for the flow state $u$ and for the AV working variable $\hat\nu$. For the flow state, we retain the boundary conditions of the original conservation law~\eqref{eq:claws}, imposed weakly through a boundary state $\func{u_\Gamma}{\Rbb^m \times \Rbb^d}{\Rbb^m}$, $(u, N) \mapsto u_\Gamma(u, N)$, with $N$ the outward unit normal, and a boundary flux, as in the DG discretization of Section~\ref{subsec:dg}. The boundary state encodes the physical condition (freestream, outflow, wall, or symmetry), and the physical conditions are stated for each case in Sections~\ref{sec:avpde_design}--\ref{sec:applications}. Wherever the boundary condition involves the normal viscous flux $F^{\mathrm{V}}(u, \nabla u) \cdot N$, the modified conservation law~\eqref{eq:modclaws} uses $\big(F^{\mathrm{V}}(u, \nabla u) + F^{\mathrm{AV}}(u, \nabla u;\, \nu)\big) \cdot N$ with $\nu = \Theta(\hat\nu)$ evaluated from the interior working variable.

For the AV working variable, we split $\partial\Omega$ into no-slip walls, on which we impose the essential condition $\hat\nu = 0$, and all other boundaries, on which we impose the homogeneous Neumann condition $\big(\eta \nabla \hat\nu\big) \cdot N = 0$.

\begin{remark}
At no-slip walls the essential condition gives $\nu = \Theta(0) = 0$, so the AV contribution to the wall flux vanishes and the AV does not compete with the physical viscosity in the boundary layer. The Neumann condition on the remaining boundaries has been adequate on every test problem considered here and introduces no additional shock-capturing parameter.
\end{remark}

\begin{remark}
Barter and Darmofal~\cite{barter2010shock} impose a Robin condition on all boundaries. Because the diffusion in the AV PDE is anisotropic through $\eta$, a Robin condition written in physical coordinates would entangle its strength with the local mesh metric; it is instead imposed in the rescaled coordinates $X = \eta^{-1/2} x$, in which the diffusion is the identity and the rate is a single dimensionless constant $c_5 \ge 0$. Transformed back to physical space, the condition reads
\begin{equation}
    \frac{c_2}{c_1}\,\big(\eta \nabla \hat\nu\big) \cdot N = \frac{c_2}{c_1}\,\bigl\|\eta^{1/2} N\bigr\|\,c_5\,\bigl(\hat\nu_\infty - \hat\nu\bigr),
    \label{eq:avpde_robin}
\end{equation}
with $\hat\nu_\infty$ the ambient AV value. The homogeneous Neumann condition used in this work is the special case $c_5 = 0$; a positive $c_5$ recovers the anisotropic Robin condition of Barter and Darmofal with the metric-tensor length scale~\cite{yano2012optimization}.
\end{remark}

\subsection{Conservation laws considered in this work}
\label{subsec:claws_used}

This subsection specializes the conservation law~\eqref{eq:claws}, the artificial viscous flux $F^{\mathrm{AV}}$, and the shock sensor $\bar{s}$ to the two physical models used in this work. Throughout, $h$ denotes the smallest eigenvalue of $H$, $p$ the polynomial degree of the DG approximation, and $\bar h = h/p$ the resolution length.

\subsubsection{Viscous Burgers' equation}
\label{subsubsec:burgers}

The 1D viscous Burgers' equation is the scalar conservation law with quadratic inviscid flux and Newtonian viscous flux,
\begin{equation}
    d = 1, \quad m = 1, \qquad
    F^{\mathrm{I}}(u) = \tfrac{1}{2} u^2, \qquad
    F^{\mathrm{V}}(u, \nabla u) = \mu \, \frac{\partial u}{\partial x}, \qquad
    S(u, \nabla u) = 0,
    \label{eq:burgers_pde}
\end{equation}
with physical viscosity $\mu \ge 0$. The artificial viscous flux mirrors $F^{\mathrm{V}}$ with the AV $\nu$ in place of $\mu$,
\begin{equation}
    F^{\mathrm{AV}}(u, \nabla u; \nu) = \nu \, \frac{\partial u}{\partial x},
    \label{eq:burgers_avflux}
\end{equation}
and the shock sensor is the magnitude of the state gradient scaled by the squared resolution length,
\begin{equation}
    \bar{s}(u, \nabla u) = \bar h^2 \, \left|\frac{\partial u}{\partial x}\right|.
    \label{eq:burgers_sensor}
\end{equation}
With $\eta = \bar h^2$ in 1D and $\nu = \Theta(\hat\nu)$, the modified conservation law~\eqref{eq:modclaws} and the AV PDE~\eqref{eq:avpde} reduce to the coupled pair
\begin{gather}
    \frac{\partial u}{\partial t} + \frac{\partial}{\partial x}\!\left(\frac{u^2}{2}\right) - \frac{\partial}{\partial x}\!\left[\big(\mu + \Theta(\hat\nu)\big) \frac{\partial u}{\partial x}\right] = 0, \label{eq:burgers_aug_state} \\
    c_1 \, \frac{\partial \hat\nu}{\partial t} - c_2 \, \frac{\partial}{\partial x}\!\left(\bar h^2 \, \frac{\partial \hat\nu}{\partial x}\right) = c_3 \, \bar h^2 \, \left|\frac{\partial u}{\partial x}\right| - c_4 \, \hat\nu, \label{eq:burgers_aug_avpde}
\end{gather}
solved for the augmented state $U = (u, \hat\nu)$.

\subsubsection{Compressible Navier--Stokes and Euler equations}
\label{subsubsec:nseuler}

The compressible Navier--Stokes equations describe the flow of a viscous, heat-conducting, calorically perfect ideal gas. The conservative state and inviscid flux are
\begin{equation}
    u = \begin{pmatrix} \rho \\ \rho \mathbf{v} \\ \rho E \end{pmatrix}, \qquad
    F^{\mathrm{I}}(u) = \begin{pmatrix} \rho \mathbf{v}^T \\ \rho \, \mathbf{v} \otimes \mathbf{v} + P \, I \\ (\rho E + P) \mathbf{v}^T \end{pmatrix},
    \label{eq:euler_flux}
\end{equation}
with density $\rho$, velocity $\mathbf{v}$, total energy per unit mass $E$, and pressure $P = (\gamma - 1) \big(\rho E - \tfrac{1}{2} \rho |\mathbf{v}|^2\big)$ with ratio of specific heats $\gamma$. The viscous flux carries the Newtonian stress tensor $\tau$ and the Fourier heat flux $q$,
\begin{equation}
    F^{\mathrm{V}}(u, \nabla u) = \begin{pmatrix} 0 \\ \tau \\ (\tau \mathbf{v} - q)^T \end{pmatrix}, \qquad
    \tau = \mu \big[\nabla \mathbf{v} + (\nabla \mathbf{v})^T - \tfrac{2}{3} (\nabla \cdot \mathbf{v})\, I\big], \qquad
    q = -\kappa \, \nabla T,
    \label{eq:ns_visc_flux}
\end{equation}
with dynamic viscosity $\mu$, thermal conductivity $\kappa = c_p \mu / \mathrm{Pr}$, temperature $T = P / (\rho R)$, ideal-gas constant $R$, and Prandtl number $\mathrm{Pr}$. Throughout this work, $\gamma = 1.4$ and $\mathrm{Pr} = 0.72$; the source is $S(u, \nabla u) = 0$, and the Euler equations are the inviscid limit $\mu = \kappa = 0$. The artificial viscous flux $F^{\mathrm{AV}}(u, \nabla u;\, \Theta(\hat\nu))$ takes the Newtonian--Fourier structure of $F^{\mathrm{V}}$ with the AV $\nu$ in place of the dynamic viscosity, zero bulk viscosity, and an artificial Prandtl number $\mathrm{Pr}_{\mathrm{av}} = 3/4$~\cite{persson2006sub}. The shock sensor is the dilation-based sensor of Moro~et~al.~\cite{moro2016dilation}, bounded by the local wave-speed envelope and activated in compressive regions through the smooth, saturating cutoff of Fern\'andez et al.~\cite{fernandez2018physics}:
\begin{equation}
    \bar{s}(u, \nabla u) = k_h \, \bar h \, \sqrt{|\mathbf{v}|^2 + c^2} \, \Phi(S_{\mathrm{d}}), \qquad
    S_{\mathrm{d}} = -k_h \, \bar h \, \frac{\nabla \cdot \mathbf{v}}{c_{\mathrm{ref}}},
    \label{eq:nseuler_sensor}
\end{equation}
where $c = \sqrt{\gamma P / \rho}$ is the local sound speed, $c_{\mathrm{ref}}$ is the critical sound speed computed pointwise from the local stagnation temperature, and $k_h$ is a dimensionless prefactor. The nondimensional dilation $S_{\mathrm{d}}$ is positive in compressive regions, and $\func{\Phi}{\Rbb}{\Rbb_{\geq 0}}$ is a smoothed form of the clipped ramp $\min(\max(S_{\mathrm{d}} - S_0, 0), S_{\max})$,
\begin{equation}
\begin{gathered}
    \Phi(S_{\mathrm{d}}) = S_{\max} + m^-\big(m^+(S_{\mathrm{d}} - S_0) - S_{\max}\big), \\
    m^+(y) = \frac{y}{2} + \frac{y}{\pi}\arctan(k_\Phi\, y) + \frac{1}{2} - \frac{\arctan k_\Phi}{\pi}, \qquad
    m^-(y) = y - m^+(y),
\end{gathered}
    \label{eq:sensor_cutoff}
\end{equation}
where $m^+$ and $m^-$ are smooth approximations of the positive and negative parts $\max(y, 0)$ and $\min(y, 0)$ that converge to them as $k_\Phi \to \infty$, $S_0$ is an activation threshold, $S_{\max} = 2/\sqrt{\gamma^2 - 1}$ is the limiting compressive dilation of an ideal-gas normal shock, and $k_\Phi$ sets the sharpness of the smoothing. Throughout this work, $k_h = 3$, $S_0 = 0.01$, and $k_\Phi = 10^4$. With $\nu = \Theta(\hat\nu)$, the modified conservation law~\eqref{eq:modclaws} and the AV PDE~\eqref{eq:avpde} reduce to
\begin{gather}
    \frac{\partial u}{\partial t} + \nabla \cdot F^{\mathrm{I}}(u) - \nabla \cdot \!\big[ F^{\mathrm{V}}(u, \nabla u) + F^{\mathrm{AV}}(u, \nabla u;\, \Theta(\hat\nu)) \big] = 0, \label{eq:nseuler_aug_state} \\
    c_1 \, \frac{\partial \hat\nu}{\partial t} - c_2 \, \nabla \cdot (\eta \, \nabla \hat\nu) = c_3 \, k_h \, \bar h \, \sqrt{|\mathbf{v}|^2 + c^2} \, \Phi(S_{\mathrm{d}}) - c_4 \, \hat\nu, \label{eq:nseuler_aug_avpde}
\end{gather}
solved for the augmented state $U = (u, \hat\nu)$, with the Euler equations recovered by setting $F^{\mathrm{V}} = 0$.

The Euler and Navier--Stokes equations are solved in one of three variable formulations. In each, the solver state $Z$ is related to the conservative state by an invertible map $u = g(Z)$ and the DG discretization is applied to the transformed system
\begin{equation}
    \frac{\partial g(Z)}{\partial t} + \nabla \cdot \left[ F^{\mathrm{I}}(g(Z)) - F^{\mathrm{V}}(g(Z), \nabla u) - F^{\mathrm{AV}}(g(Z), \nabla u;\, \nu) \right] = S(g(Z), \nabla u), \qquad
    \nabla u = \tfrac{\partial g}{\partial Z}(Z) \, \nabla Z.
    \label{eq:cov}
\end{equation}
The change of variables acts on the flow state only; the AV working variable $\hat\nu$ is left untransformed, so the augmented solver state is $(Z, \hat\nu)$. The three formulations are listed in Table~\ref{tab:formulations}, where $s = \ln(P / \rho^\gamma)$ is the specific thermodynamic entropy.

\begin{table}[hbt!]
\centering
\caption{Variable formulations for the Euler and Navier--Stokes equations considered in this work.}
\label{tab:formulations}
\begin{tabular}{l l l}
    \toprule
    Formulation & State variables $Z$ & Motivation \\
    \midrule
    Conservative & $(\rho,\, \rho \mathbf{v},\, \rho E)$ & Standard DG formulation \\
    Entropy & $\left(\dfrac{\gamma - s}{\gamma - 1} - \dfrac{\rho |\mathbf{v}|^2}{2P},\, \dfrac{\rho \mathbf{v}}{P},\, -\dfrac{\rho}{P}\right)$ & Entropy stability \\
    Logarithmic & $(\ln \rho,\, \mathbf{v},\, \ln P)$ & Positivity of $\rho$ and $P$ \\
    \bottomrule
\end{tabular}
\end{table}

\begin{remark}
The conservative state can develop negative density or pressure during the nonlinear iteration, and comparative studies on entropy-stable DG schemes report that alternative variable formulations reduce the incidence of these failures~\cite{alberti2024comparative}; Section~\ref{sec:var_robustness} carries out the analogous comparison for PDE-based shock capturing. The entropy variables symmetrize the governing equations, which is the basis of entropy-stable DG discretizations~\cite{hughes1986new}. The logarithmic variables guarantee $\rho = \exp(Z_1) > 0$ and $P = \exp(Z_{d+2}) > 0$ for any finite $Z$, which removes one common source of solver failure.
\end{remark}

\subsection{Discontinuous Galerkin discretization}
\label{subsec:dg}

The augmented system~\eqref{eq:augmented} is discretized with the interior-penalty discontinuous Galerkin (IPDG) method of Hartmann and Houston~\cite{hartmann2008optimal}. Let $\mathcal{K}_h$ be a mesh of $\Omega$ into non-overlapping, potentially curved elements $K$, and let
\begin{equation}
    \mathcal{V}_h^p = \big\{ v \in [L^2(\Omega)]^{m+1} \;:\; v|_K \in [\mathcal{P}^p(K)]^{m+1}, \; \forall \, K \in \mathcal{K}_h \big\},
    \label{eq:dg_space}
\end{equation}
where $\mathcal{P}^p(K)$ is the space of polynomials of total degree at most $p$ on $K$, be the trial and test space for the augmented state $U_h = (u_h, \hat\nu_h)$. The discrete steady-state problem reads: find $U_h \in \mathcal{V}_h^p$ such that
\begin{equation}
    r_h(U_h;\, w_h) = 0 \qquad \forall \, w_h \in \mathcal{V}_h^p,
    \label{eq:dg_problem}
\end{equation}
where the IPDG residual splits into volume, interior-face, and boundary-face contributions,
\begin{equation}
    r_h(U_h;\, w_h) = r_{h,\Omega} + r_{h,\Sigma} + r_{h,\Gamma},
    \label{eq:dg_residual_split}
\end{equation}
each a semilinear form of the trial $U_h$ and the test $w_h$. The face terms depend on the outward unit normal $N$ to $\partial K$, the two traces $(\cdot)^+$ from inside $K$ and $(\cdot)^-$ from the neighboring element, and the boundary state $U_\Gamma$ on $\partial \Omega$. The viscous fluxes considered in this work are linear in the state gradient,
\begin{equation}
    \mathcal{F}^{\mathrm{V}}(U, \nabla U) = \mathcal{G}(U) : \nabla U, \qquad
    \big(\mathcal{G}(U) : A\big)_{ik} = \sum_{j=1}^{m+1} \sum_{l=1}^{d} \mathcal{G}_{ikjl}(U)\, A_{jl}, \quad A \in \Rbb^{(m+1) \times d},
    \label{eq:visc_flux_tensor}
\end{equation}
where $\func{\mathcal{G}}{\Rbb^{m+1}}{\Rbb^{(m+1) \times d \times (m+1) \times d}}$ is a fourth-order tensor and $:$ denotes the double contraction over its last two indices. With the average operator $\{\!\{v\}\!\} = \tfrac{1}{2}(v^+ + v^-)$ on interior faces, the three contributions read
\begin{align}
    r_{h,\Omega} &= \sum_{K \in \mathcal{K}_h} \int_K \Big[ -\nabla w_h : \big( \mathcal{F}^{\mathrm{I}}(U_h) - \mathcal{F}^{\mathrm{V}}(U_h, \nabla U_h) \big) - w_h \cdot \mathcal{S}(U_h, \nabla U_h) \Big]\, dV, \label{eq:dg_vol} \\
    r_{h,\Sigma} &= \sum_{K \in \mathcal{K}_h} \int_{\partial K \setminus \partial \Omega} \Big[ w_h^+ \cdot \Big( \mathcal{H}(U_h^+, U_h^-;\, N) - \{\!\{ \mathcal{F}^{\mathrm{V}}(U_h, \nabla U_h) \}\!\}\, N + \sigma_{\Sigma} \, \big( \{\!\{ \mathcal{G}(U_h) \}\!\} : (U_h^+ - U_h^-) \otimes N \big)\, N \Big) \nonumber \\
    &\quad - \tfrac{1}{2}\, \nabla w_h^+ : \big( \mathcal{G}(U_h^+) : (U_h^+ - U_h^-) \otimes N \big) \Big]\, dS, \label{eq:dg_iface} \\
    r_{h,\Gamma} &= \sum_{K \in \mathcal{K}_h} \int_{\partial K \cap \partial \Omega} \Big[ w_h \cdot \Big( \mathcal{F}^{\mathrm{I}}_\Gamma(U_h, U_\Gamma;\, N) - \mathcal{F}^{\mathrm{V}}_\Gamma(U_\Gamma, \nabla U_h)\, N + \sigma_{\Gamma}\, \big( \mathcal{G}(U_\Gamma) : (U_h - U_\Gamma) \otimes N \big)\, N \Big) \nonumber \\
    &\quad - \nabla w_h : \big( \mathcal{G}_\Gamma(U_\Gamma) : (U_h - U_\Gamma) \otimes N \big) \Big]\, dS, \label{eq:dg_bface}
\end{align}
with penalty coefficients $\sigma_\Sigma = C^\Sigma_{\mathrm{ip}}\, p^2 / h_f$ and $\sigma_\Gamma = C^\Gamma_{\mathrm{ip}}\, p^2 / h_f$, where $h_f$, a measure of the element size normal to the face~\cite{hartmann2008optimal}, is $\min(|K^+|, |K^-|)/|f|$ on an interior face $f$ shared by the elements $K^+$ and $K^-$, and $|K|/|f|$ on a boundary face, with $C^\Sigma_{\mathrm{ip}} = C^\Gamma_{\mathrm{ip}} = 10$ in all cases. Here $\func{\mathcal{H}}{\Rbb^{m+1} \times \Rbb^{m+1} \times \Rbb^d}{\Rbb^{m+1}}$ is the inviscid numerical flux, taken as Roe's approximate Riemann solver~\cite{roe1981approximate} for every model, with the Harten--Hyman entropy fix~\cite{harten1983self} for the Euler and Navier--Stokes systems. The inviscid boundary flux $\mathcal{F}^{\mathrm{I}}_\Gamma(U, U_\Gamma;\, N) = \mathcal{F}^{\mathrm{I}}(U_\Gamma)\, N$ imposes the inviscid boundary conditions through the boundary state $U_\Gamma = U_\Gamma(U, N)$ of Section~\ref{subsec:av_bnd}. The boundary viscous flux $\mathcal{F}^{\mathrm{V}}_\Gamma(U, \nabla U) = \mathcal{G}_\Gamma(U) : \nabla U$ imposes the viscous boundary conditions through $\mathcal{G}_\Gamma$, the corresponding modification of $\mathcal{G}$. The symmetrizing term uses $\mathcal{G}_\Gamma$ and the penalty term uses the unmodified $\mathcal{G}$. When an alternative variable formulation of Section~\ref{subsubsec:nseuler} is used, the residual~\eqref{eq:dg_vol}--\eqref{eq:dg_bface} is assembled in the solver state $(Z, \hat\nu)$: the jumps and penalty terms act on $Z$, and $\mathcal{G}$ is the gradient tensor of the transformed viscous flux in~\eqref{eq:cov}, which includes the map Jacobian $\partial g / \partial Z$.

\begin{remark}
The boundary symmetrizing term uses $\mathcal{G}_\Gamma$ so that it is the adjoint of the consistency term and the viscous discretization remains adjoint consistent~\cite{hartmann2008optimal}. The boundary penalty uses the unmodified $\mathcal{G}$ because lifting the jump through a closure that zeroes rows of $\mathcal{G}$, as at slip walls and symmetry planes, under-penalizes the corresponding components; the standard penalty acts on every component of the jump regardless of the closure.
\end{remark}

The nonlinear system~\eqref{eq:dg_problem} is solved by pseudo-transient continuation (PTC)~\cite{kelley1998convergence}, initialized from a $p = 0$ solution computed on the same mesh, with the update limiting and CFL control of Modisette~\cite{modisette2011automated}: each pseudo-time step takes one Newton step, limited by a physicality check and a line search on the unsteady residual, and the CFL number doubles after a fully accepted step and is reduced tenfold after a rejected one. The linear system of each Newton step is solved by restarted GMRES, preconditioned with a block incomplete-LU factorization with no fill-in and minimum-discarded-fill reordering~\cite{persson_newton-gmres_2008}. PTC is declared converged when $\|r_h(U_h;\,\cdot)\|_\infty < 10^{-8}$, stalled (the CFL$\to 0$ failure mode of Section~\ref{sec:var_robustness}) if the CFL number drops below $10^{-8}$, and terminated as \emph{maxit} at a cap of 1000 iterations. The time constant $c_1$ does not affect the converged solution, only the PTC iteration count; its effect on the robustness and cost of the solve is studied in Section~\ref{subsec:c1_study}.

\section{Resolution-aware scaling of the AV PDE}
\label{sec:avpde_design}
\label{subsec:scaling}

A captured shock should be thin enough to localize the discontinuity but thick enough that the polynomial space can represent it monotonically: too thin, and the oscillations the regularization was introduced to suppress reappear; too thick, and the shock is smeared over more elements than necessary. Because the AV dominates the physical viscosity inside a captured shock, the captured-shock thickness scales with the AV~\cite{whitham1974linear, persson2006sub}, so a target shock thickness is a target AV thickness. Prior PDE-based formulations~\cite{barter2010shock, yano2012optimization, nguyen2024optimal} tie the diffusion length of the AV PDE to the element size $h$ alone, even though $p$-refinement shrinks the smallest thickness the polynomial space can represent monotonically, which scales as $\bar h = h/p$ (Section~\ref{subsubsec:dmin}).

Section~\ref{subsubsec:dmin} identifies the target thickness, Section~\ref{subsubsec:hbar2} derives the scaling of the AV PDE that makes its solution self-similar in that target, and Section~\ref{subsec:scaling_verify} verifies the prediction. The self-similarity analysis of Section~\ref{subsubsec:hbar2} is the central contribution of this section: it yields concrete conditions, on the powers of $\bar h$ carried by the diffusion and source terms, that any AV PDE must satisfy for its solution to be self-similar under $hp$-refinement, and Remark~\ref{rem:conditions} measures the prior formulations against them.

\subsection{Minimum resolvable shock thickness}
\label{subsubsec:dmin}

\begin{figure}[hbt!]
    \centering
    \IfFileExists{_tikz/l2proj_motivation.tex}%
        {\input{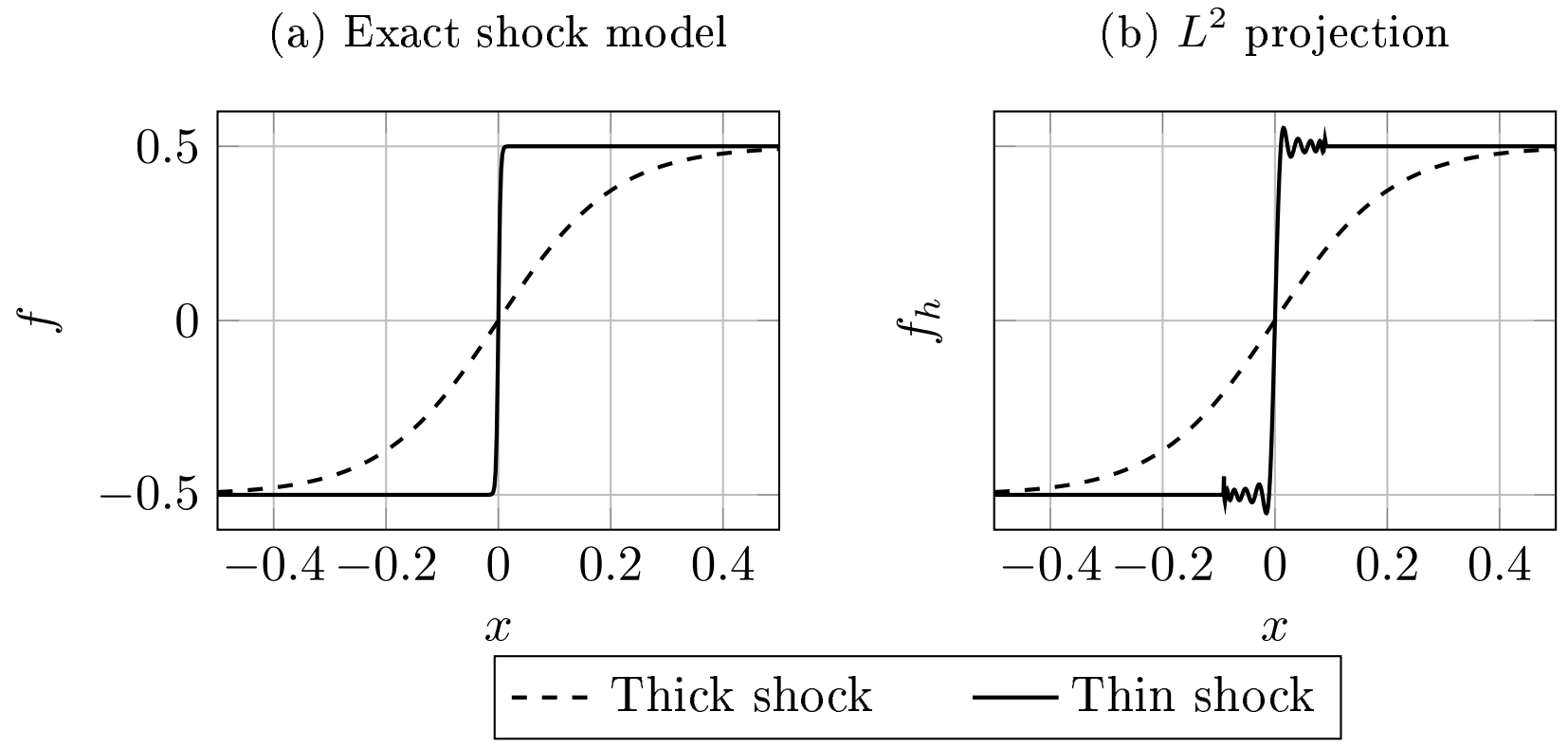}}%
        {\fbox{\parbox[c][4cm][c]{0.7\linewidth}{\centering\scriptsize $L^2$ projection motivation\\[0.3em]\textbf{[2-panel: exact tanh shock (left) vs element-local $L^2$ projection (right); two thicknesses $\delta$]}\\[0.3em]\textit{thin oscillates, thick smears}}}}
    \caption{(a)~A tanh shock model and (b)~its $L^2$ projection onto a piecewise polynomial space.}
    \label{fig:l2proj_motivation}
\end{figure}

Let $\delta_{\min}(h, p)$ denote the smallest shock thickness that the piecewise polynomial space of degree $p$ on elements of size $h$ can represent monotonically. We characterize it by the element-local $L^2$ projection of the model shock
\begin{equation}
    f(x) = \frac{\alpha}{2}\,\tanh\!\left(\frac{2\,\operatorname{atanh}(1-\beta)}{\delta}\,(x - x_{\mathrm{s}})\right),
    \label{eq:tanh_profile}
\end{equation}
where $\alpha$ is the jump amplitude, $x_{\mathrm{s}}$ the shock center, and $\delta$ the thickness, defined so that $f(x_{\mathrm{s}} \pm \delta/2) = \pm \alpha (1 - \beta)/2$ with $\beta = 0.01$. Two regimes arise: a thin profile whose projection overshoots and oscillates, and a thick profile whose projection is monotonic but spans many elements (Fig.~\ref{fig:l2proj_motivation}). A projection is classified as monotonic if the slope of the projected profile does not change sign, beyond a small numerical tolerance, within $|x - x_{\mathrm{s}}| \le \delta$. For a given shock position $x_{\mathrm{s}}$, the monotonicity threshold $\vartheta(h, p, x_{\mathrm{s}})$ is the smallest $\delta$ for which the projection is monotonic. The minimum resolvable thickness is its worst case over positions, $\delta_{\min}(h, p) = \max_{x_{\mathrm{s}}} \vartheta(h, p, x_{\mathrm{s}})$, and we write $\vartheta_{\min}(h, p) = \min_{x_{\mathrm{s}}} \vartheta(h, p, x_{\mathrm{s}})$ for the best case.

\begin{figure}[hbt!]
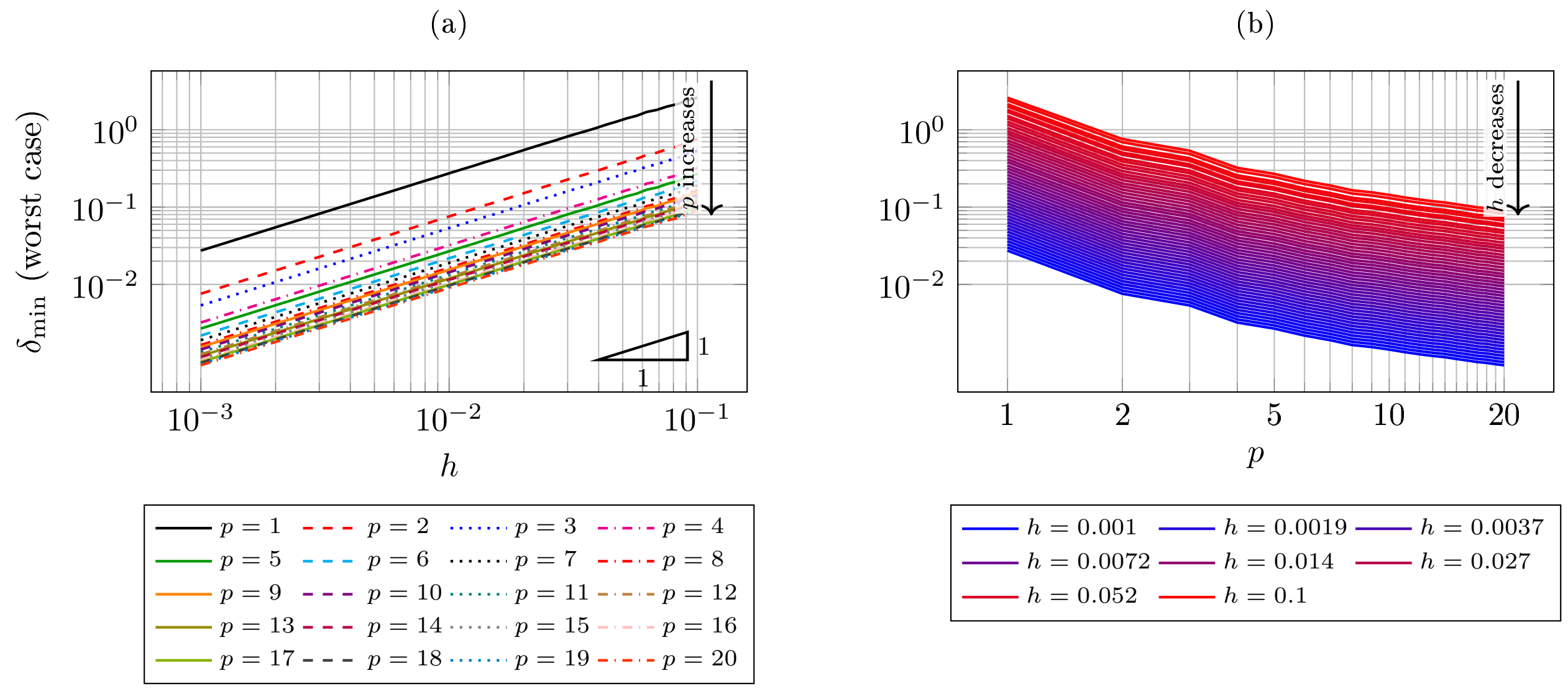

    \centering
    \begin{minipage}[t]{0.48\textwidth}
        \centering
        \IfFileExists{_tikz/delta_vs_h.tex}%
            {\input{_tikz/delta_vs_h.tex}}%
            {\fbox{\parbox[c][5cm][c]{0.95\linewidth}{\centering\scriptsize $\delta_{\min}$ vs.\ $h$\\[0.3em]\textbf{[one curve per $p$, log-log]}}}}    \end{minipage}\hfill
    \begin{minipage}[t]{0.48\textwidth}
        \centering
        \IfFileExists{_tikz/delta_vs_p.tex}%
            {\input{_tikz/delta_vs_p.tex}}%
            {\fbox{\parbox[c][5cm][c]{0.95\linewidth}{\centering\scriptsize $\delta_{\min}$ vs.\ $p$\\[0.3em]\textbf{[one curve per $h$, log-log]}}}}    \end{minipage}
    \caption{Minimum resolvable shock thickness $\delta_{\min}$ from the $L^2$ projection study. (a)~$\delta_{\min}$ vs.\ $h$, one curve per $p$. (b)~$\delta_{\min}$ vs.\ $p$, one curve per $h$.}
    \label{fig:delta_min}
\end{figure}

The study sweeps $p \in \{1, \ldots, 20\}$, 50 log-spaced values of $h \in [10^{-3}, 10^{-1}]$, and 20 shock positions within an element (Fig.~\ref{fig:delta_min}). At fixed $p$, $\delta_{\min}$ increases linearly with $h$; at fixed $h$, it decreases with $p$ slightly faster than $1/p$; and at higher $p$ it falls below the element size, so the polynomial space resolves sub-cell shock profiles. A joint least-squares fit yields
\begin{equation}
    \delta_{\min} \approx 16.56 \cdot h^{1.0} \cdot p^{-1.04},
    \label{eq:deltamin_fit}
\end{equation}
so $\delta_{\min}$ is proportional to the resolution length $\bar h = h/p$ to within a few percent in each exponent. The design target that follows is an AV whose magnitude scales with $\bar h$, which produces a captured-shock thickness $\delta \sim \bar h$.

\begin{figure}[hbt!]
    \centering
    \begin{minipage}[t]{0.48\textwidth}
        \centering
        \IfFileExists{_tikz/collapse_test.tex}%
            {\input{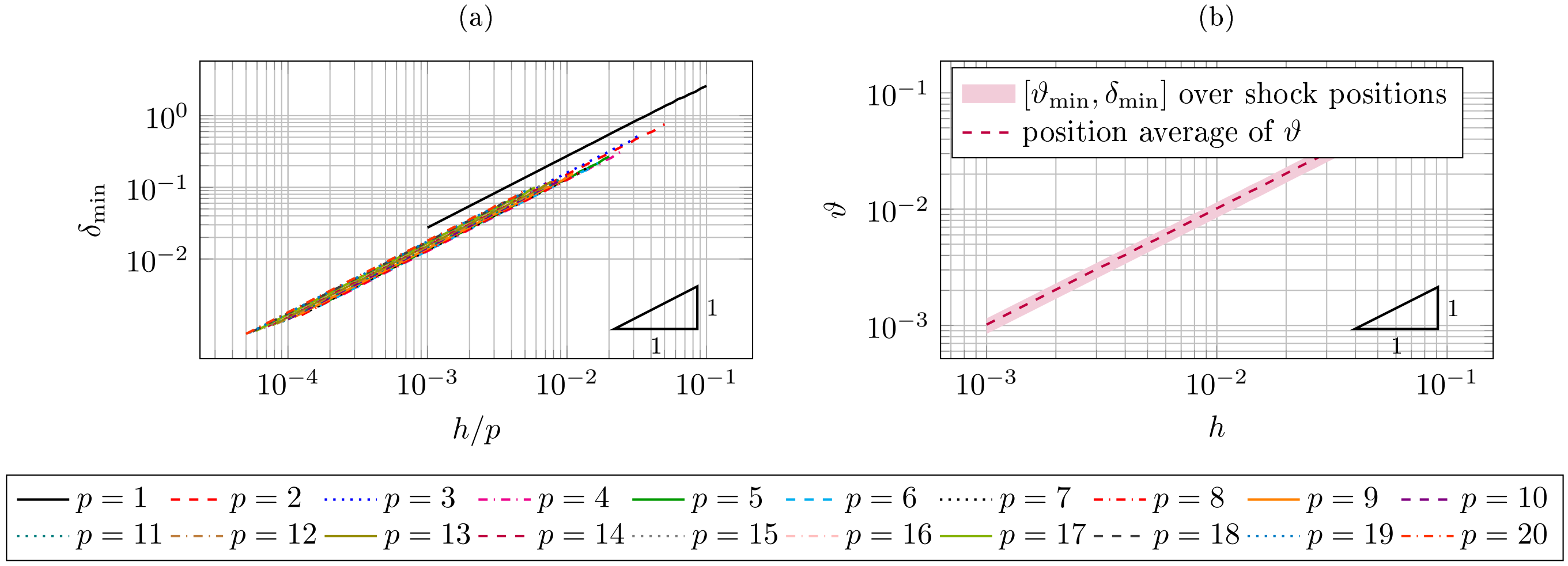}}%
            {\fbox{\parbox[c][5cm][c]{0.95\linewidth}{\centering\scriptsize $\delta_{\min}$ vs.\ $h/p$\\[0.3em]\textbf{[collapse test, one curve per $p$, log-log]}}}}    \end{minipage}\hfill
    \begin{minipage}[t]{0.48\textwidth}
        \centering
        \IfFileExists{_tikz/position_envelope.tex}%
            {\begin{tikzpicture}[baseline=(current bounding box.north)]
\begin{loglogaxis}[title={(b)},
    tick label style={font=\normalsize},
  xlabel={$h$},
  ylabel={$\vartheta$},
  grid=both,
  width=\linewidth,
  height=5cm,
  legend style={at={(0.02, 0.98)}, anchor=north west},
  legend cell align={left},
]
\addplot[name path=upper, draw=none, forget plot] coordinates {
  (0.001, 0.0011481717228889465)
  (0.0010985411419875584, 0.0012610480189323425)
  (0.0012067926406393288, 0.0013849511742591858)
  (0.0013257113655901094, 0.001521781086921692)
  (0.0014563484775012444, 0.0016712397336959839)
  (0.0015998587196060573, 0.0018355995416641235)
  (0.0017575106248547913, 0.0020178407430648804)
  (0.0019306977288832496, 0.0022164732217788696)
  (0.0021209508879201904, 0.0024350732564926147)
  (0.002329951810515372, 0.0026731938123703003)
  (0.002559547922699536, 0.0029364973306655884)
  (0.002811768697974231, 0.0032250583171844482)
  (0.0030888435964774815, 0.003543645143508911)
  (0.00339322177189533, 0.0038920342922210693)
  (0.0037275937203149418, 0.004276186227798462)
  (0.004094915062380423, 0.004695802927017212)
  (0.004498432668969444, 0.005160123109817505)
  (0.004941713361323833, 0.005670040845870972)
  (0.005428675439323859, 0.006223022937774658)
  (0.005963623316594642, 0.006833970546722412)
  (0.006551285568595509, 0.007504522800445557)
  (0.0071968567300115215, 0.00826030969619751)
  (0.007906043210907703, 0.009076297283172607)
  (0.008685113737513529, 0.009940564632415771)
  (0.009540954763499945, 0.010934770107269287)
  (0.010481131341546858, 0.012022852897644043)
  (0.011513953993264475, 0.013197064399719238)
  (0.012648552168552958, 0.014441609382629395)
  (0.013894954943731374, 0.015946030616760254)
  (0.015264179671752334, 0.01739680767059326)
  (0.016768329368110083, 0.01913607120513916)
  (0.018420699693267165, 0.021067261695861816)
  (0.020235896477251564, 0.02319514751434326)
  (0.022229964825261943, 0.02551436424255371)
  (0.02442053094548651, 0.028003454208374023)
  (0.02682695795279726, 0.03061652183532715)
  (0.029470517025518114, 0.033768415451049805)
  (0.032374575428176434, 0.03703713417053223)
  (0.03556480306223129, 0.04028677940368652)
  (0.03906939937054617, 0.0441586971282959)
  (0.04291934260128778, 0.0488591194152832)
  (0.04714866363457394, 0.053403377532958984)
  (0.0517947467923121, 0.05887746810913086)
  (0.05689866029018296, 0.06378412246704102)
  (0.06250551925273973, 0.0717616081237793)
  (0.06866488450043001, 0.07654428482055664)
  (0.07543120063354619, 0.08504629135131836)
  (0.08286427728546843, 0.09185075759887695)
  (0.09102981779915217, 0.10437488555908203)
  (0.1, 0.11481761932373047)
};
\addplot[name path=lower, draw=none, forget plot] coordinates {
  (0.001, 0.0008450821042060852)
  (0.0010985411419875584, 0.0009281560778617859)
  (0.0012067926406393288, 0.0010194256901741028)
  (0.0013257113655901094, 0.0011200830340385437)
  (0.0014563484775012444, 0.0012301281094551086)
  (0.0015998587196060573, 0.0013510510325431824)
  (0.0017575106248547913, 0.0014851614832878113)
  (0.0019306977288832496, 0.0016314536333084106)
  (0.0021209508879201904, 0.0017923861742019653)
  (0.002329951810515372, 0.0019676238298416138)
  (0.002559547922699536, 0.0021613389253616333)
  (0.002811768697974231, 0.002373829483985901)
  (0.0030888435964774815, 0.0026082247495651245)
  (0.00339322177189533, 0.0028646737337112427)
  (0.0037275937203149418, 0.0031472742557525635)
  (0.004094915062380423, 0.003456324338912964)
  (0.004498432668969444, 0.003798156976699829)
  (0.004941713361323833, 0.004173070192337036)
  (0.005428675439323859, 0.004580169916152954)
  (0.005963623316594642, 0.005030184984207153)
  (0.006551285568595509, 0.005523413419723511)
  (0.0071968567300115215, 0.006079375743865967)
  (0.007906043210907703, 0.006680190563201904)
  (0.008685113737513529, 0.007316768169403076)
  (0.009540954763499945, 0.00804811716079712)
  (0.010481131341546858, 0.008848607540130615)
  (0.011513953993264475, 0.009713470935821533)
  (0.012648552168552958, 0.010629594326019287)
  (0.013894954943731374, 0.011736452579498291)
  (0.015264179671752334, 0.01280367374420166)
  (0.016768329368110083, 0.01408398151397705)
  (0.018420699693267165, 0.015504956245422363)
  (0.020235896477251564, 0.01707136631011963)
  (0.022229964825261943, 0.018778443336486816)
  (0.02442053094548651, 0.020610690116882324)
  (0.02682695795279726, 0.022533535957336426)
  (0.029470517025518114, 0.024853944778442383)
  (0.032374575428176434, 0.027257204055786133)
  (0.03556480306223129, 0.029648542404174805)
  (0.03906939937054617, 0.03250002861022949)
  (0.04291934260128778, 0.03595709800720215)
  (0.04714866363457394, 0.039299726486206055)
  (0.0517947467923121, 0.04333138465881348)
  (0.05689866029018296, 0.04694104194641113)
  (0.06250551925273973, 0.05280733108520508)
  (0.06866488450043001, 0.05632638931274414)
  (0.07543120063354619, 0.06258249282836914)
  (0.08286427728546843, 0.06758928298950195)
  (0.09102981779915217, 0.07680177688598633)
  (0.1, 0.08447885513305664)
};
\addplot[purple!20] fill between[of=upper and lower];
\addlegendentry{$[\vartheta_{\min}, \delta_{\min}]$ over shock positions}

\addplot[purple, dashed, thick] coordinates {
  (0.001, 0.0010150782763957977)
  (0.0010985411419875584, 0.0011148490011692047)
  (0.0012067926406393288, 0.0012244470417499542)
  (0.0013257113655901094, 0.0013453606516122818)
  (0.0014563484775012444, 0.001477530226111412)
  (0.0015998587196060573, 0.0016228128224611282)
  (0.0017575106248547913, 0.0017839614301919937)
  (0.0019306977288832496, 0.001959577202796936)
  (0.0021209508879201904, 0.002152852714061737)
  (0.002329951810515372, 0.002363361418247223)
  (0.002559547922699536, 0.0025960803031921387)
  (0.002811768697974231, 0.002851303666830063)
  (0.0030888435964774815, 0.003132931888103485)
  (0.00339322177189533, 0.0034095607697963715)
  (0.0037275937203149418, 0.0037804991006851196)
  (0.004094915062380423, 0.004113808274269104)
  (0.004498432668969444, 0.004562109708786011)
  (0.004941713361323833, 0.005012691020965576)
  (0.005428675439323859, 0.005451634526252747)
  (0.005963623316594642, 0.0060421377420425415)
  (0.006551285568595509, 0.006634503602981567)
  (0.0071968567300115215, 0.007302612066268921)
  (0.007906043210907703, 0.008024245500564575)
  (0.008685113737513529, 0.008788436651229858)
  (0.009540954763499945, 0.009667396545410156)
  (0.010481131341546858, 0.010629013180732727)
  (0.011513953993264475, 0.01166754961013794)
  (0.012648552168552958, 0.01276804506778717)
  (0.013894954943731374, 0.0140981525182724)
  (0.015264179671752334, 0.015379846096038818)
  (0.016768329368110083, 0.016917705535888672)
  (0.018420699693267165, 0.018625080585479736)
  (0.020235896477251564, 0.020506441593170166)
  (0.022229964825261943, 0.022556990385055542)
  (0.02442053094548651, 0.02475762367248535)
  (0.02682695795279726, 0.027068346738815308)
  (0.029470517025518114, 0.029854893684387207)
  (0.032374575428176434, 0.03244590759277344)
  (0.03556480306223129, 0.035616397857666016)
  (0.03906939937054617, 0.03904104232788086)
  (0.04291934260128778, 0.043194472789764404)
  (0.04714866363457394, 0.04721254110336304)
  (0.0517947467923121, 0.05205410718917847)
  (0.05689866029018296, 0.0563923716545105)
  (0.06250551925273973, 0.06344127655029297)
  (0.06866488450043001, 0.06767034530639648)
  (0.07543120063354619, 0.07518982887268066)
  (0.08286427728546843, 0.08120489120483398)
  (0.09102981779915217, 0.0922774076461792)
  (0.1, 0.10150539875030518)
};
\addlegendentry{position average of $\vartheta$}

\logLogSlopeTriangle{0.9}{0.15}{0.1}{1}{thick, black}
\end{loglogaxis}
\end{tikzpicture}}%
            {\fbox{\parbox[c][5cm][c]{0.95\linewidth}{\centering\scriptsize $\delta_{\min}$ envelope vs.\ $h$, $p=14$\\[0.3em]\textbf{[position envelope]}}}}    \end{minipage}
    \figurelegend{collapse_position_legend}
    \caption{(a)~Minimum resolvable thickness $\delta_{\min}$ against the resolution length $h/p$, one curve per $p \in \{1, \ldots, 20\}$, with a slope-one reference triangle. (b)~Dependence of the monotonicity threshold on shock position at $p = 14$: for each $h$, the threshold $\vartheta(h, p, x_{\mathrm{s}})$ is computed at 20 shock positions $x_{\mathrm{s}}$ within the element, and the shaded band spans its range $[\vartheta_{\min}, \delta_{\min}]$ from the best-case position (lower edge) to the worst-case position (upper edge, which is the $\delta_{\min}$ plotted in panel~(a)); the dashed line is the average of $\vartheta$ over the 20 positions.}
    \label{fig:collapse_test}
\end{figure}

\begin{remark}
The element-local $L^2$ projection is used as a proxy for the DG approximation: on elements with weak inter-element coupling, which is the case in smooth regions and across well-resolved shocks, the DG solution is close to the element-local projection, so the projection's resolution limit transfers to the DG approximation. Two diagnostics in Fig.~\ref{fig:collapse_test} support the fit~\eqref{eq:deltamin_fit}. In panel~(a) the $\delta_{\min}$ curves for $p \geq 2$ collapse onto a tight cluster of slope one when plotted against $h/p$, with the $p = 1$ curve slightly above it (linear elements being the most restrictive space for thin transitions). Panel~(b) shows that the dependence on shock position is modest and uniform: at $p = 14$, the degree with the largest spread between the best and worst shock positions, the ratio $\delta_{\min}/\vartheta_{\min}$ is $1.36$ at every $h$, and it is smaller at every other $p$. The worst-case threshold $\delta_{\min}$ therefore exceeds the best-case one by at most $36\%$, and because the factor is independent of $h$ the scaling of $\delta_{\min}$ with $h/p$ does not depend on where the shock sits within the element.
\end{remark}

\subsection{Self-similar scaling of diffusion and source terms}
\label{subsubsec:hbar2}

The target thickness of Section~\ref{subsubsec:dmin} sets the length scale of the AV PDE. Consider the steady 1D linear AV PDE with the power $q$ of the resolution length in the diffusion term left open,
\begin{equation}
    -\bar h^{q} \frac{\partial^2 \hat\nu}{\partial x^2} = \bar s - \hat\nu.
    \label{eq:pde_1d_unscaled}
\end{equation}
The sensor $\bar s$ is the complete source of the AV PDE~\eqref{eq:avpde}, including whatever dependence on $\bar h$ its definition carries, so no separate power of $\bar h$ multiplies it. The decay term carries no power of $\bar h$ either, because a factor common to the whole equation is absorbed into the coefficients $c_2$, $c_3$, $c_4$, which are omitted here for the same reason. A viscous shock with viscosity $\nu$ and jump $\Delta u$ has thickness $\delta = \mathcal{O}(\nu / \Delta u)$~\cite{whitham1974linear}, and inside a captured shock the AV dominates the physical viscosity, so $\delta = \mathcal{O}(\hat\nu / \Delta u)$; the target~\eqref{eq:deltamin_fit} requires $\delta = \mathcal{O}(\bar h)$. Since $\Delta u$ is independent of $\bar h$, the two together require $\hat\nu = \mathcal{O}(\bar h)$. We therefore seek a solution of the form
\begin{equation}
    \hat\nu(x) = \bar h \, \tilde\nu(\xi), \qquad \xi = x/\bar h,
    \label{eq:ansatz}
\end{equation}
where $\tilde\nu$ is a universal profile independent of the discretization, so that the AV has width and amplitude proportional to $\bar h$ in physical space. We call the coupled solution $(u, \hat\nu)$ self-similar if the AV has the form~\eqref{eq:ansatz} and the state is a fixed profile of the rescaled coordinate, $u(x) = \tilde u(\xi)$, with $\tilde\nu$ and $\tilde u$ independent of $\bar h$: the captured shock has a fixed jump and a thickness proportional to $\bar h$. On a self-similar solution the sensor is a fixed function of $\xi$ up to a power of $\bar h$,
\begin{equation}
    \bar s(x) = \bar h^{r} \, \tilde s(\xi),
    \label{eq:sensor_class}
\end{equation}
where $\tilde s$ is independent of $\bar h$ and the exponent $r$ depends on the sensor. For example, for the Burgers' sensor in~\eqref{eq:burgers_sensor}, the chain rule on a self-similar solution gives
\begin{equation}
    \bar s(x) = \bar h^2 \left|\frac{\partial u}{\partial x}(x)\right| = \bar h^2 \, \frac{1}{\bar h} \, |\tilde u'(\xi)| = \bar h \, |\tilde u'(\xi)|,
    \label{eq:sensor_scaling}
\end{equation}
so $r = 1$ with $\tilde s = |\tilde u'|$. Substituting~\eqref{eq:ansatz} and~\eqref{eq:sensor_class} into~\eqref{eq:pde_1d_unscaled} gives
\begin{equation}
    -\bar h^{q - 1} \, \frac{\partial^2 \tilde\nu}{\partial \xi^2} = \bar h^{r} \, \tilde s - \bar h \, \tilde\nu.
    \label{eq:pde_1d_unscaled_rescaled}
\end{equation}
The equation for $\tilde\nu$ is free of $\bar h$ if and only if the three powers coincide,
\begin{equation}
    q = 2, \qquad r = 1.
    \label{eq:exponents}
\end{equation}
These are the conditions an AV PDE must satisfy for its solution to be self-similar in $\bar h$: the diffusion must carry $\bar h^2$, and the sensor must scale as $\bar h$ on a self-similar solution. For the sensors used in this work, which satisfy the second condition (Remark~\ref{rem:sensor_scaling}), the equation is
\begin{equation}
    -\bar h^2 \, \frac{\partial^2 \hat\nu}{\partial x^2} = \bar s - \hat\nu,
    \label{eq:pde_1d_scaled}
\end{equation}
which under~\eqref{eq:ansatz} and~\eqref{eq:sensor_class} reduces to the $\bar h$-free equation
\begin{equation}
    -\frac{\partial^2 \tilde\nu}{\partial \xi^2} = \tilde s - \tilde\nu.
    \label{eq:pde_1d_scaled_rescaled}
\end{equation}
Solving it once yields a universal profile $\tilde\nu(\xi)$ that, via $\hat\nu(x) = \bar h \, \tilde\nu(x/\bar h)$, describes the AV at every discretization.

\begin{remark}
\label{rem:conditions}
A gradient sensor, which contributes a factor $1/\bar h$ as in~\eqref{eq:sensor_scaling}, multiplied by $\bar h^2$, and an $\bar h$-independent switch multiplied by $\bar h$ both have $r = 1$ in~\eqref{eq:sensor_class} and satisfy~\eqref{eq:exponents}. The formulations of~\cite{barter2010shock, yano2012optimization} use the switch form. Their diffusion scales as $h^2$ and their source as $h/p$, so the terms have the balance~\eqref{eq:exponents} requires in $h$ but not in $p$: the AV is self-similar under $h$-refinement at fixed $p$, and under $p$-refinement the width stays at $h$ while the amplitude falls as $1/p$. Nguyen et al.~\cite{nguyen2024optimal} place $h_{\min}^2$ in front of the diffusion and use the dilation as the source without any compensating length-scale factor, so their formulation is not self-similar in $h$ (or $p$) by construction. The formulation of Section~\ref{sec:govern} is the only one of the four that satisfies both conditions under $h$-, $p$-, and $hp$-refinement.
\end{remark}

\begin{remark}
\label{rem:sensor_scaling}
The sensors of Section~\ref{subsec:claws_used} satisfy the second condition in~\eqref{eq:exponents}. For Burgers' equation, the steady form of~\eqref{eq:burgers_aug_avpde} is $-c_2 \bar h^2 \hat\nu'' = c_3 \bar h^2 |u'| - c_4 \hat\nu$, so $q = 2$, and $r = 1$ by~\eqref{eq:sensor_scaling}. For the Navier--Stokes equations, the nondimensional dilation $S_{\mathrm{d}}$ in~\eqref{eq:nseuler_sensor} is independent of $\bar h$ on a self-similar solution, because $\nabla \cdot \mathbf{v} = \tilde v'(\xi)/\bar h$ across the captured shock and $S_{\mathrm{d}}$ carries a compensating factor of $\bar h$, so $S_{\mathrm{d}} = -k_h \tilde v'(\xi)/c_{\mathrm{ref}}$. The cutoff $\Phi(S_{\mathrm{d}})$ is therefore a function of $\xi$ alone, in both its ramp and its saturated regime, and the single explicit factor $\bar h$ in the amplitude of~\eqref{eq:nseuler_sensor} gives $r = 1$. Together with the diffusion tensor $\eta = H^2/p^2$ of~\eqref{eq:eta_tensor}, the AV PDE of Section~\ref{subsec:scpde} is therefore self-similar in $\bar h$ for both models.
\end{remark}

\begin{remark}
The self-similarity is a statement about the continuum equation~\eqref{eq:pde_1d_scaled_rescaled}. The DG approximation inherits it only in the well-resolved regime of Section~\ref{subsubsec:dmin}, where the captured shock stays sufficiently thicker than $\delta_{\min}$.
\end{remark}

\subsection{Verification}
\label{subsec:scaling_verify}

This subsection tests the self-similarity prediction on the 1D Burgers' equation and on a Mach~5 inviscid flow over a half-cylinder. Every Euler and Navier--Stokes case in this work is solved in entropy variables unless otherwise stated, which Section~\ref{sec:var_robustness} shows afford the widest convergence window of the three formulations.

\subsubsection{Self-similarity under joint \texorpdfstring{$hp$}{hp}-refinement on Burgers' equation}
\label{subsec:val_burgers}
\label{subsubsec:selfsim_burgers}

Consider the inviscid Burgers' equation ($\mu = 0$) of Section~\ref{subsubsec:burgers} on $\Omega = [-40, 40]$, with the Dirichlet boundary states $u_{\mathrm{L}} = +1$ at the left boundary and $u_{\mathrm{R}} = -1$ at the right. Both boundaries are inflow boundaries for this data, and the solution is a stationary shock at $x = 0$ with jump $\Delta u = 2$. The AV working variable takes the homogeneous Neumann condition of Section~\ref{subsec:av_bnd} at both boundaries. We use $c_1 = 1$ throughout and no AV filter ($\nu = \hat\nu$), except in the under-resolved case of Section~\ref{subsubsec:collapse_limit}, whose oscillations drive $\hat\nu$ negative; there the filter uses $a = 10^{-3}$ and $\nu_0 = 0$.

The first experiment on this problem (test case: \case{burgers-hp}) tests the prediction of Section~\ref{subsubsec:hbar2} that the converged AV is self-similar under joint $hp$-refinement: the physical AV $\hat\nu(x) = \bar h\, \tilde\nu(x/\bar h)$ depends on the discretization only through $\bar h$, so that in the rescaled coordinate $\xi = x/\bar h$ the profiles from different $(h, p)$ pairs collapse onto the single shape $\tilde\nu(\xi)$, which depends only on the coefficients $(c_2, c_3, c_4)$ and on the problem through the rescaled sensor $\tilde s$, not on the discretization. Five $(h, p)$ pairs are chosen so that $\bar h = h/p$ matches the target sweep $\bar h \in \{0.4, 0.2, 0.1, 0.05, 0.025\}$ while $h$ and $p$ both vary (Table~\ref{tab:burgers_joint_sweep}), with $c_2 = c_3 = 5$ and $c_4 = 1/c_2$.

\begin{table}[hbt!]
    \centering
    \caption{The five $(h, p)$ pairs of the joint refinement sweep, targeting matched $\bar h = h/p$ values. Approximated entries reflect rounding of $h = 2L / n_\mathrm{el}$, with $L = 40$ the half-length of the domain, to integer $n_\mathrm{el}$.}
    \label{tab:burgers_joint_sweep}
    \begin{tabular}{c c c c}
        \toprule
        $p$ & $h$ & $n_\mathrm{el}$ & $\bar h$ \\
        \midrule
        2 & $0.800$         & 100 & $0.4$ \\
        3 & $\approx 0.602$ & 133 & $\approx 0.2$ \\
        5 & $0.500$         & 160 & $0.1$ \\
        7 & $\approx 0.349$ & 229 & $\approx 0.05$ \\
        9 & $\approx 0.225$ & 356 & $\approx 0.025$ \\
        \bottomrule
    \end{tabular}
\end{table}

\begin{figure}[hbt!]
    \centering
    \IfFileExists{_tikz/282s_4panel.tex}{\input{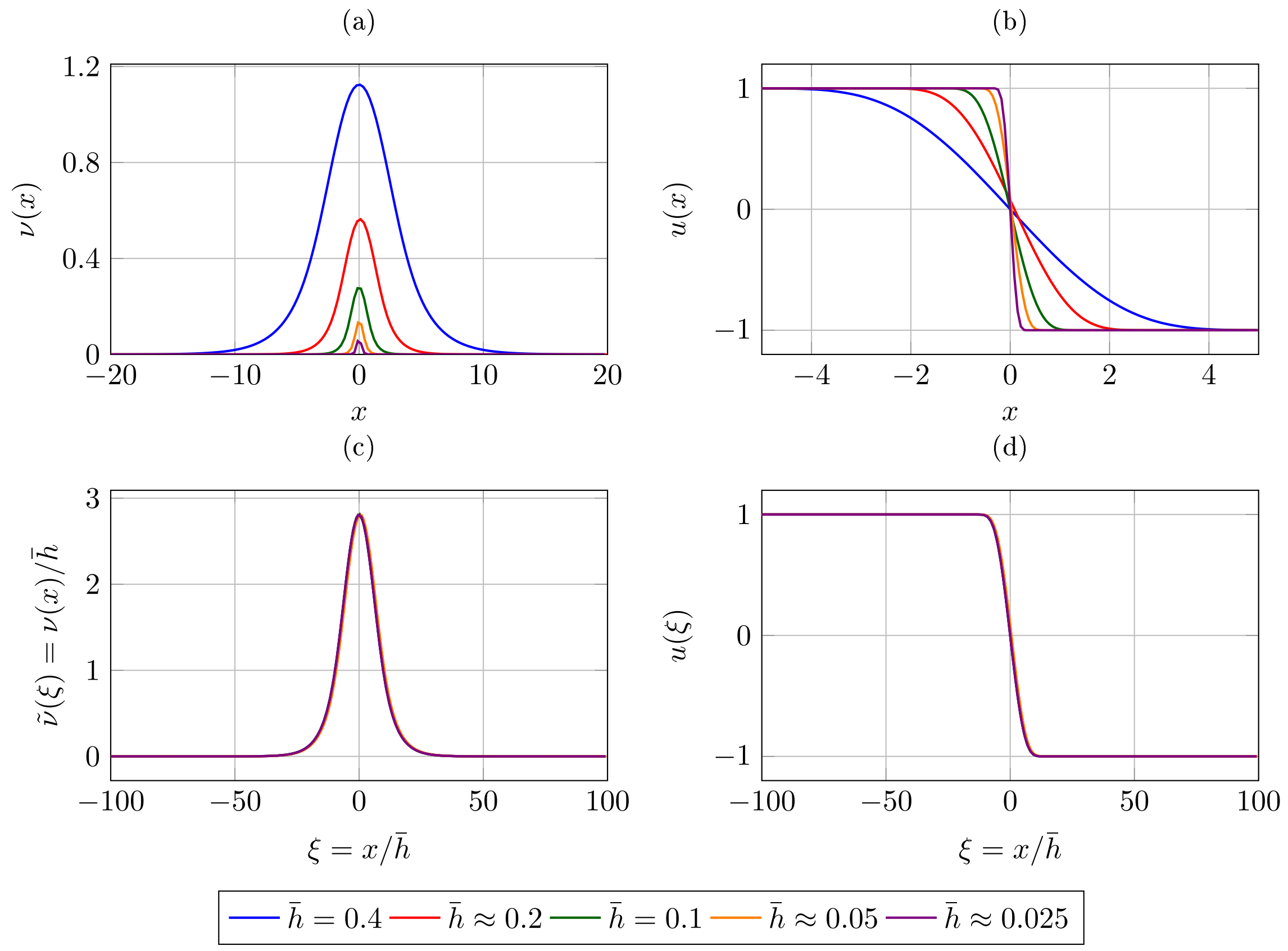}}{\fbox{\parbox[c][6cm][c]{0.85\linewidth}{\centering\scriptsize Burgers' equation 4-panel: top row raw $\nu(x)$, $u(x)$; bottom row rescaled, collapsed onto single shape; joint $(h, p)$ sweep}}}
    \figurelegend{hbarlegendSSjointfix}
    \caption{\case{burgers-hp}: converged AV and state across the joint $(h, p)$ sweep of Table~\ref{tab:burgers_joint_sweep}. (a)~$\nu(x)$; (b)~$u(x)$; (c)~$\nu/\bar h$ and (d)~$u$ against $\xi = x/\bar h$.}
    \label{fig:selfsim_burgers}
\end{figure}

The AV peak and width both shrink as $\bar h$ decreases and the shock in $u$ thins correspondingly, with the asymptotic values $\pm 1$ unchanged, and in rescaled coordinates the five curves, each from a distinct $(h, p)$ pair, collapse onto a single profile for both the AV and the state (Fig.~\ref{fig:selfsim_burgers}). For fixed $(c_2, c_3, c_4)$ the rescaled state--AV solution is therefore independent of the discretization, as predicted, and the physical solution depends on it only through $\bar h$.

\begin{remark}
With $\mu = 0$ the AV is the sole regularization and the self-similarity is exact. This exactness does not carry to viscous flows, but the AV still dominates the regularization inside the shock for viscosity levels of practical interest, so the shock-region scaling established here carries over approximately; Section~\ref{sec:applications} tests it on viscous flows.
\end{remark}

\subsubsection{Loss of the collapse below the well-resolved regime}
\label{subsubsec:collapse_limit}

\begin{figure}[hbt!]
    \centering
    \IfFileExists{_tikz/324s_4panel.tex}{\input{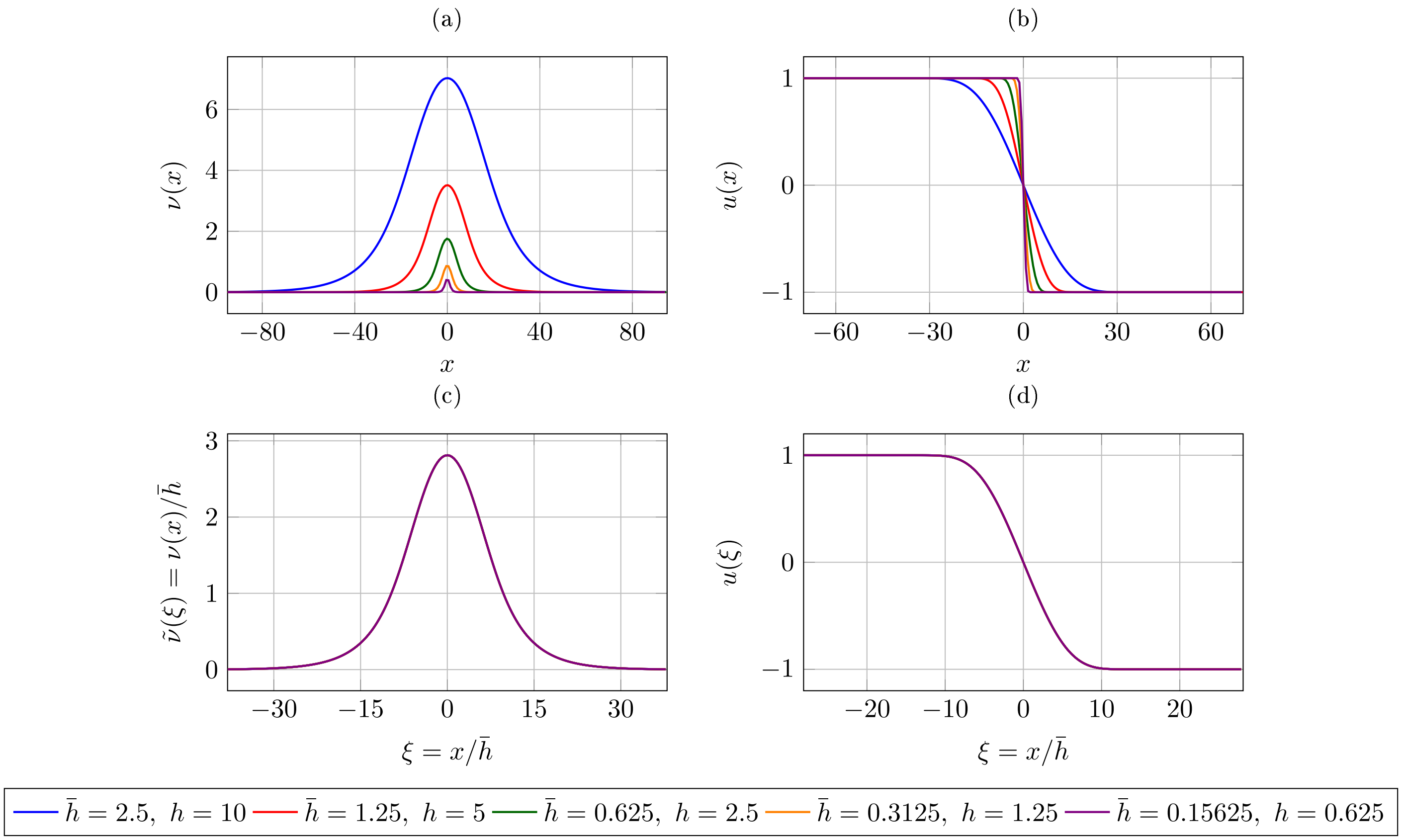}}{\fbox{\parbox[c][6cm][c]{0.85\linewidth}{\centering\scriptsize Resolved Burgers' 4-panel ($c_2=c_3=5$): raw and rescaled profiles collapse}}}
    \figurelegend{hbarlegendRESp4}
    \caption{\case{burgers-limit}, well-resolved case ($c_2 = c_3 = 5$, $p = 4$, $h$ varied): (a)~$\nu(x)$; (b)~$u(x)$; (c)~$\nu/\bar h$ and (d)~$u$ against $\xi = x/\bar h$.}
    \label{fig:burgers_overp_resolved}
\end{figure}

\begin{figure}[hbt!]
    \centering
    \IfFileExists{_tikz/314s_collapse_4panel.tex}{\input{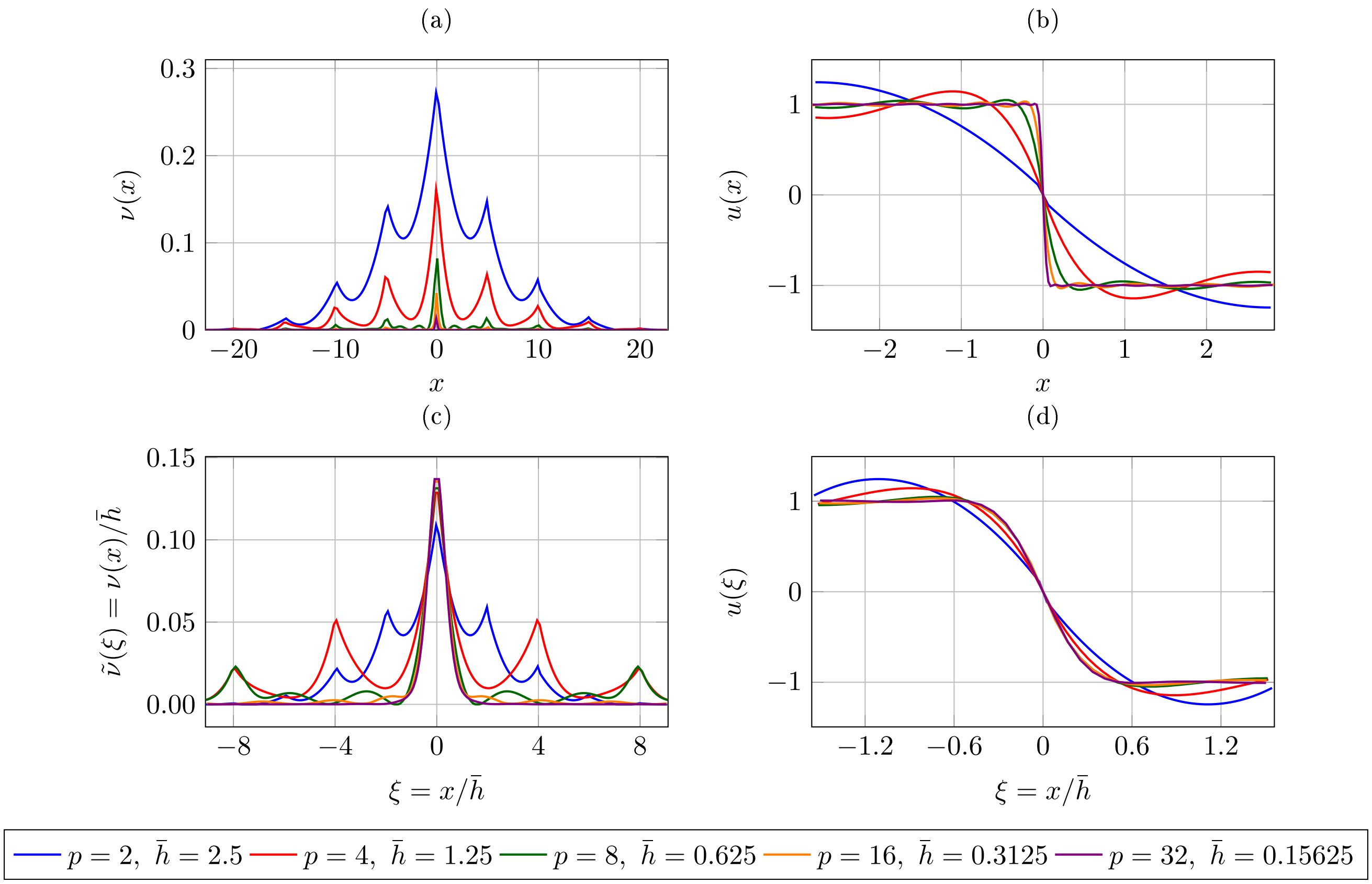}}{\fbox{\parbox[c][6cm][c]{0.85\linewidth}{\centering\scriptsize Under-resolved Burgers' 4-panel ($c_2=c_3=0.25$): rescaled profiles fail to collapse}}}
    \figurelegend{hbarlegendFILTh5c0p25coll}
    \caption{\case{burgers-limit}, under-resolved case ($c_2 = c_3 = 0.25$, $h = 5$, $p$ varied) at the same $\bar h$ as Fig.~\ref{fig:burgers_overp_resolved}: (a)~$\nu(x)$; (b)~$u(x)$; (c)~$\nu/\bar h$ and (d)~$u$ against $\xi = x/\bar h$.}
    \label{fig:burgers_overp_filter_collapse}
\end{figure}

The collapse of Section~\ref{subsubsec:selfsim_burgers} holds while the AV keeps the captured shock monotone, which requires a thickness sufficiently above $\delta_{\min}(h, p)$. To reach the boundary of that regime we vary the AV PDE coefficients rather than the discretization (test case: \case{burgers-limit}). Two settings are compared at the five matched resolution lengths $\bar h = 2.5$, $1.25$, $0.625$, $0.3125$, and $0.15625$. A well-resolved case with $c_2 = c_3 = 5$ sweeps $\bar h$ with $p = 4$ fixed on the larger domain $\Omega = [-200, 200]$, where its wider shock stays clear of the boundaries, and an under-resolved one with $c_2 = c_3 = 0.25$ fixes $h = 5$ on $\Omega = [-40, 40]$ and varies $p$. With the larger coefficients the rescaled profiles collapse onto a single curve (Fig.~\ref{fig:burgers_overp_resolved}), reproducing the self-similarity of Section~\ref{subsubsec:selfsim_burgers} at these resolution lengths. With the smaller coefficients the converged shock thins to about one resolution length and its representation degrades (Fig.~\ref{fig:burgers_overp_filter_collapse}): the state overshoots its asymptotic values $\pm 1$ near the shock and the rescaled profiles no longer collapse.

The contrast is quantified through three measures (Table~\ref{tab:collapse_quality}): the captured-shock thickness $\delta$, the distance between the points where $u$ crosses $+(1-\beta)$ and $-(1-\beta)$, which for the unit half-amplitude of this problem applies the thickness definition of the model shock~\eqref{eq:tanh_profile} to the computed solution; the AV peak $\nu_{\mathrm{peak}}$; and the AV integral width $\ell_{\mathrm{int}} = \int\nu\,\mathrm{d}x/(2\,\nu_{\mathrm{peak}})$. The thickness $\delta$ and the integral width $\ell_{\mathrm{int}}$ are used here, rather than the $e$-folding half-width $\ell$ of Section~\ref{subsubsec:c4_test_burgers}, because they remain well defined for the oscillatory profiles of the under-resolved case, whereas $\ell$ presumes a profile that decays monotonically from a single peak. In the well-resolved case the three $\bar h$-rescaled ratios are constant to all printed digits. With $c_2 = c_3 = 0.25$, the AV width varies across the sweep, with $\ell_{\mathrm{int}}/\bar h$ falling from 1.68 to 0.54, the AV peak varies less ($\nu_{\mathrm{peak}}/\bar h$ from 0.109 to 0.137, still near $c_3/2 = 0.125$), and the shock thickness stays at about one resolution length, far below $\delta_{\min}$. The worst monotonicity violation is $\approx 0.16$ of the jump against $\sim 10^{-5}$ in the well-resolved case. The self-similar collapse is therefore conditional on the AV holding the shock sufficiently above $\delta_{\min}$ that the state stays monotone, as it does at $c_2 = c_3 = 5$ ($\delta \approx 19.6\,\bar h$ against $\delta_{\min} \approx 16.6\,\bar h$).

\begin{table}[hbt!]
    \centering
    \caption{\case{burgers-limit}: rescaled shock thickness, AV peak, and AV integral width at matched $\bar h$ for the well-resolved case ($h$-refinement at $p = 4$), whose three ratios are the same at every $\bar h$, and the under-resolved one ($p$-refinement at $h = 5$), whose ratios do not collapse.}
    \label{tab:collapse_quality}
    \IfFileExists{_tikz/collapse_ratios_table.tex}{{\small
\begin{tabular}{l c c c rrr}
\toprule
case & $\bar h$ & $h$ & $p$ & $\delta/\bar h$ & $\nu_{\mathrm{peak}}/\bar h$ & $\ell_{\mathrm{int}}/\bar h$ \\
\midrule
well-resolved, $c_2 = c_3 = 5$     & 2.5      & 10        & 4  & 19.59 & 2.81  & 8.89 \\
                                   & 1.25     & 5         & 4  & 19.59 & 2.81  & 8.89 \\
                                   & 0.625    & 2.5       & 4  & 19.59 & 2.81  & 8.89 \\
                                   & 0.3125   & 1.25      & 4  & 19.59 & 2.81  & 8.89 \\
                                   & 0.15625  & 0.625     & 4  & 19.59 & 2.81  & 8.89 \\
\midrule
under-resolved, $c_2 = c_3 = 0.25$ & 2.5      & 5         & 2  & 1.19  & 0.109 & 1.68 \\
                                   & 1.25     & 5         & 4  & 1.01 & 0.129 & 1.68 \\
                                   & 0.625    & 5         & 8  & 0.97 & 0.131 & 1.09 \\
                                   & 0.3125   & 5         & 16 & 0.99 & 0.135 & 0.83 \\
                                   & 0.15625  & 5         & 32 & 0.98 & 0.137 & 0.54 \\
\bottomrule
\end{tabular}
}
}{\fbox{\parbox[c][3cm][c]{0.7\linewidth}{\centering\scriptsize collapse\_ratios\_table.tex}}}
\end{table}

\subsubsection{Mach~5 inviscid flow over a half-cylinder}
\label{subsec:val_euler_cyl}
\label{subsubsec:selfsim_cyl}

Consider Mach~5 inviscid flow over a half-cylinder of unit radius, governed by the steady Euler equations of Section~\ref{subsubsec:nseuler}. The domain lies between the cylinder and an outer half-ellipse centered on it, with semi-axis 3 along the freestream and 6 across it~\cite{dong2025rpadapt}. The freestream approaches along the axis of symmetry: the cylinder is a slip wall, the outer arc carries the freestream state, and the two straight edges downstream of the cylinder are supersonic outflow boundaries. The AV working variable takes the homogeneous Neumann condition on every boundary, there being no no-slip wall. The AV filter uses $a = \nu_0 = 10^{-6}$, and the AV PDE uses $c_1 = 10^{-4}$.

\begin{figure}[hbt!]
    \centering
    \includegraphics[width=0.5\linewidth]{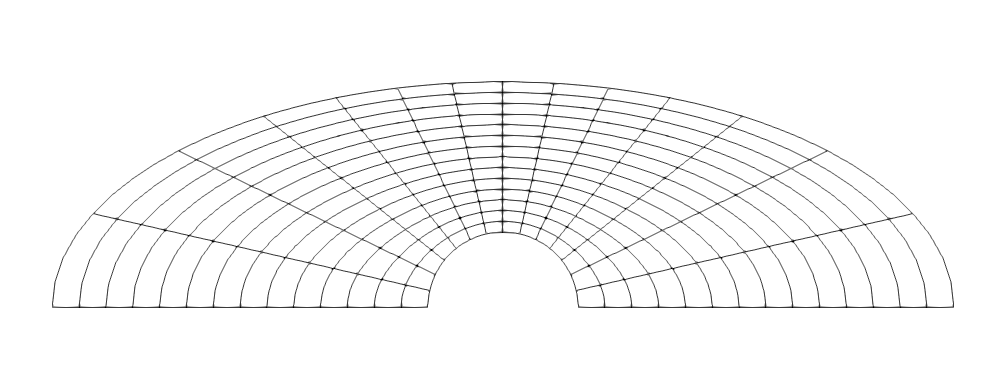}
    \caption{Structured $n_{\mathrm{r}} \times n_{\mathrm{a}}$ half-cylinder mesh at the coarsest level ($n_{\mathrm{r}} = n_{\mathrm{a}} = 14$, 196 elements), with $n_{\mathrm{r}}$ elements along the stagnation line and $n_{\mathrm{a}}$ along the surface. Boundaries: wall on the cylinder, freestream on the outer arc, outflow on the two straight edges.}
    \label{fig:halfcyl_mesh}
\end{figure}

The domain is discretized by the structured $n_{\mathrm{r}} \times n_{\mathrm{a}}$ mesh of Fig.~\ref{fig:halfcyl_mesh}, uniform in the radial and angular directions. The experiment on this problem (test case: \case{cyl-euler-hp}) tests whether the scaling of Section~\ref{subsubsec:hbar2}, which the Navier--Stokes sensor inherits as explained in Remark~\ref{rem:sensor_scaling}, describes the AV around a curved bow shock, and in particular whether $\bar h$ alone, not $h$ or $p$ individually, governs the converged profile. We fix $n_{\mathrm{a}} = 14$ and choose three $(\bar h, p)$ combinations with distinct $p$, $(p, n_{\mathrm{r}}) \in \{(3, 20), (2, 60), (4, 60)\}$, corresponding to $\bar h \in \{1/30, 1/60, 1/120\}$, with $h = 2/n_{\mathrm{r}}$ the radial element size along the stagnation line. The AV PDE parameters are $(c_2, c_3) = (1, 1)$ with $c_4 = 1/c_2$.

\begin{figure}[hbt!]
    \centering
    \begin{tikzpicture}
\pgfplotsset{
  compat=1.18,
  every axis plot/.append style={thick},
}
\begin{groupplot}[
  group style={
    group size=2 by 2,
    horizontal sep=0.07\textwidth,
    vertical sep=1.6cm,
  },
  width=0.45\textwidth, height=5cm,
  grid=both,
  tick label style={font=\normalsize},
]

\nextgroupplot[title={(a)}, xlabel={$y$}, ylabel={$\nu(y)$},
               xmin=1.0, xmax=3.0, ymin=0, ymax=0.31, ytick={0,0.1,0.2,0.3},
               legend to name=hbar3legend288,
               legend columns=3, legend cell align=left]
\addplot[color=red]   table[x=y, y=AV] {_tikz/288s_p3_n20.dat};
\addlegendentry{$p=3,\ n_{\mathrm{r}}=20$}
\addplot[color=black] table[x=y, y=AV] {_tikz/288s_p2_n60.dat};
\addlegendentry{$p=2,\ n_{\mathrm{r}}=60$}
\addplot[color=blue]  table[x=y, y=AV] {_tikz/288s_p4_n60.dat};
\addlegendentry{$p=4,\ n_{\mathrm{r}}=60$}

\nextgroupplot[title={(b)}, xlabel={$y$}, ylabel={$M(y)$},
               xmin=1.0, xmax=3.0, ymin=0, ymax=5.5, ytick={0,1.5,3,4.5}]
\addplot[color=red]   table[x=y, y=Mach] {_tikz/288s_p3_n20.dat};
\addplot[color=black] table[x=y, y=Mach] {_tikz/288s_p2_n60.dat};
\addplot[color=blue]  table[x=y, y=Mach] {_tikz/288s_p4_n60.dat};

\nextgroupplot[title={(c)}, xlabel={$(y-y_\mathrm{shock})/\bar h$}, ylabel={$\nu(y)/\bar h$},
               xmin=-10, xmax=10, ymin=0, ymax=10, ytick={0,3,6,9},
               ]
\addplot[color=red]   table[x expr={(\thisrow{y} - 1 - 0.5187408704) / 0.0333333333}, y expr={\thisrow{AV} / 0.0333333333}] {_tikz/288s_p3_n20.dat};
\addplot[color=black] table[x expr={(\thisrow{y} - 1 - 0.4967518759) / 0.0166666667}, y expr={\thisrow{AV} / 0.0166666667}] {_tikz/288s_p2_n60.dat};
\addplot[color=blue]  table[x expr={(\thisrow{y} - 1 - 0.4887558779) / 0.0083333333}, y expr={\thisrow{AV} / 0.0083333333}] {_tikz/288s_p4_n60.dat};

\nextgroupplot[title={(d)}, xlabel={$(y-y_\mathrm{shock})/\bar h$}, ylabel={$M(y)$},
               xmin=-10, xmax=10, ymin=0, ymax=5.5, ytick={0,1.5,3,4.5}]
\addplot[color=red]   table[x expr={(\thisrow{y} - 1 - 0.5187408704) / 0.0333333333}, y expr={\thisrow{Mach}}] {_tikz/288s_p3_n20.dat};
\addplot[color=black] table[x expr={(\thisrow{y} - 1 - 0.4967518759) / 0.0166666667}, y expr={\thisrow{Mach}}] {_tikz/288s_p2_n60.dat};
\addplot[color=blue]  table[x expr={(\thisrow{y} - 1 - 0.4887558779) / 0.0083333333}, y expr={\thisrow{Mach}}] {_tikz/288s_p4_n60.dat};

\end{groupplot}
\end{tikzpicture}
    \figurelegend{hbar3legend288}
    \caption{\case{cyl-euler-hp}: stagnation-line profiles at three $(\bar h, p)$ combinations with distinct $p$. (a)~$\nu$ and (b)~Mach number against the physical coordinate $y$; (c)~$\nu/\bar h$ and (d)~Mach number against the shifted, rescaled coordinate $(y - y_\mathrm{shock})/\bar h$, where $y_\mathrm{shock}$ is the location of the AV maximum.}
    \label{fig:cyl_stagprofiles_288s}
\end{figure}
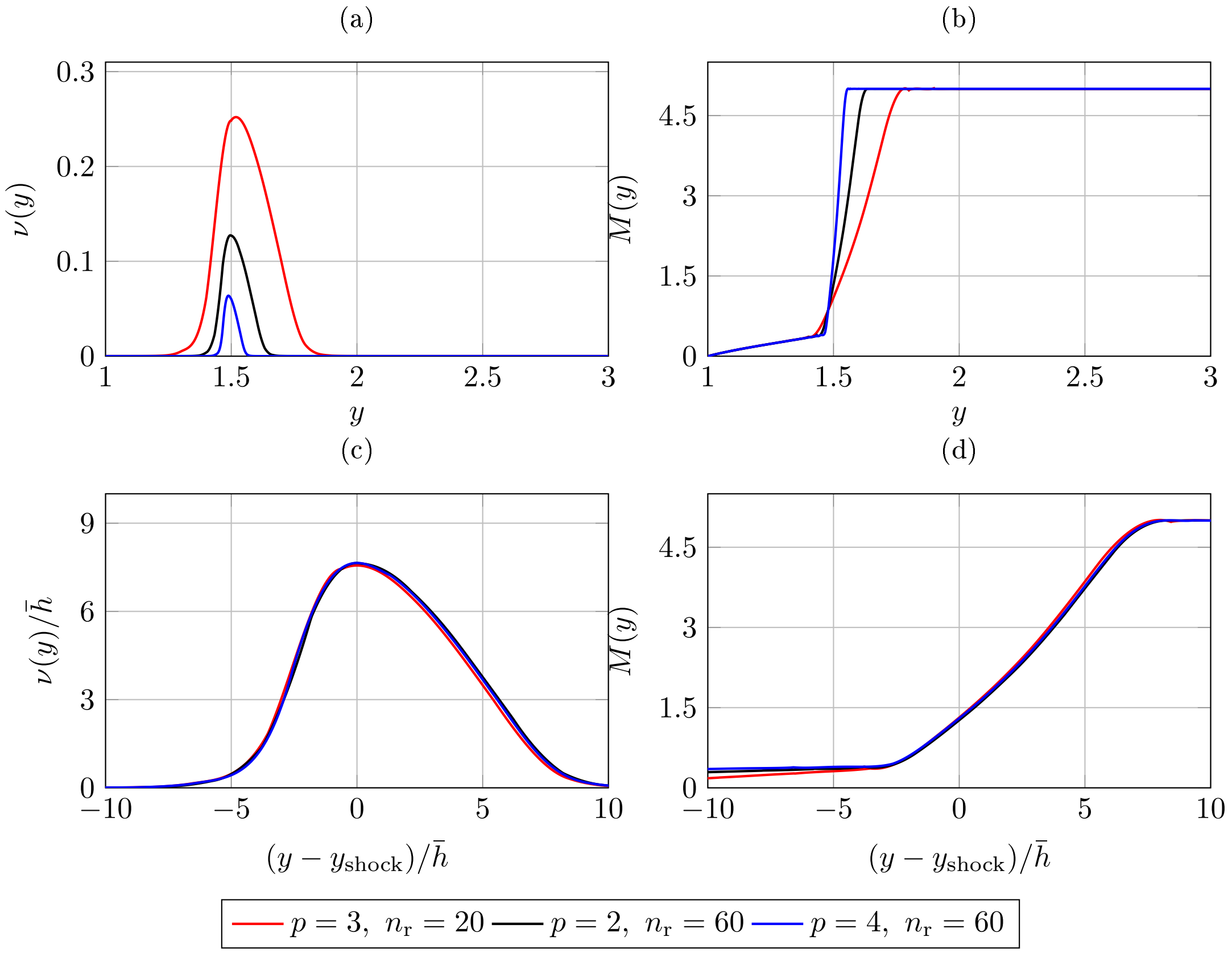

At each level $\nu$ concentrates in a band aligned with the bow shock, and as $\bar h$ shrinks the band thins and weakens while the shock thickness decreases (Fig.~\ref{fig:cyl_2d_refinement}). The captured shock position $y_\mathrm{shock}$, taken as the stagnation-line location of the AV maximum, depends on $\bar h$, because the shock relaxes toward its inviscid position as the AV weakens; the profiles are therefore shifted to $y - y_\mathrm{shock}$ before rescaling by the local resolution length at the shock. The AV profiles shrink in peak and width with $\bar h$ and collapse onto a single curve once rescaled, and the Mach profiles thin under refinement and collapse onto a single shock transition in the rescaled coordinate (Fig.~\ref{fig:cyl_stagprofiles_288s}). The collapse of three different $(h, p)$ decompositions onto one profile establishes that the converged AV depends on the discretization only through $\bar h = h/p$, within the well-resolved regime of Section~\ref{subsubsec:dmin}.

\begin{figure}[hbt!]
    \centering
    \begin{tikzpicture}
        \tikzset{cell/.style={inner sep=2pt, outer sep=0pt, anchor=north west}}
        \node[cell] (A) at (0,0) {\includegraphics[width=0.48\linewidth]{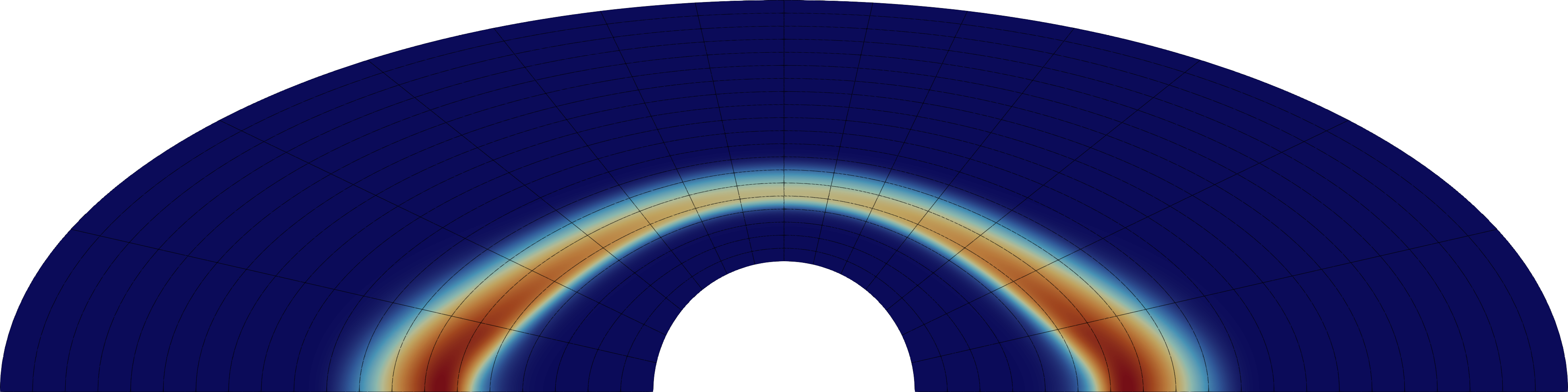}};
        \node[cell] (B) at ([xshift=0.15cm]A.north east) {\includegraphics[width=0.48\linewidth]{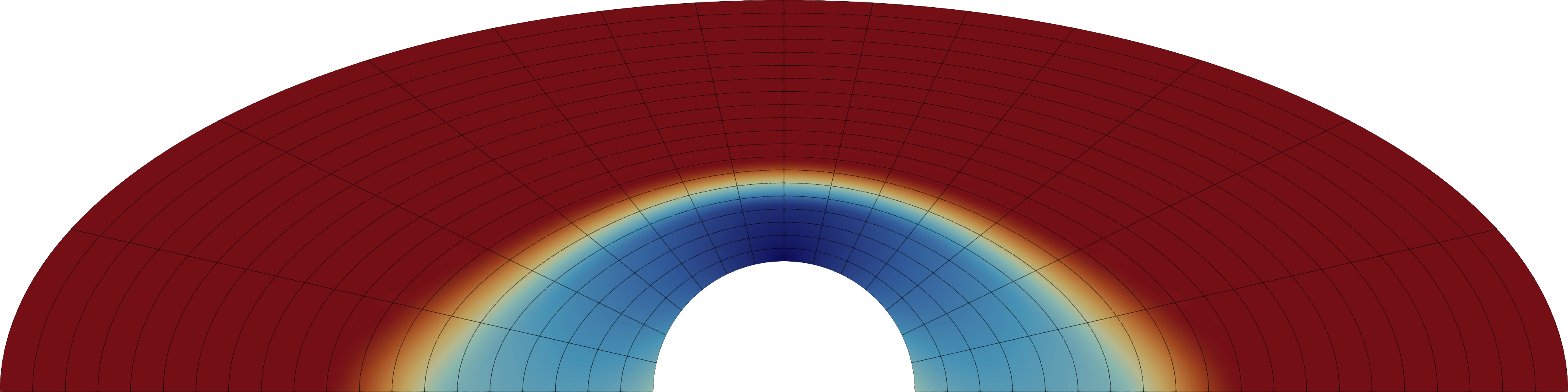}};
        \node[cell] (C) at ([yshift=-0.15cm]A.south west) {\includegraphics[width=0.48\linewidth]{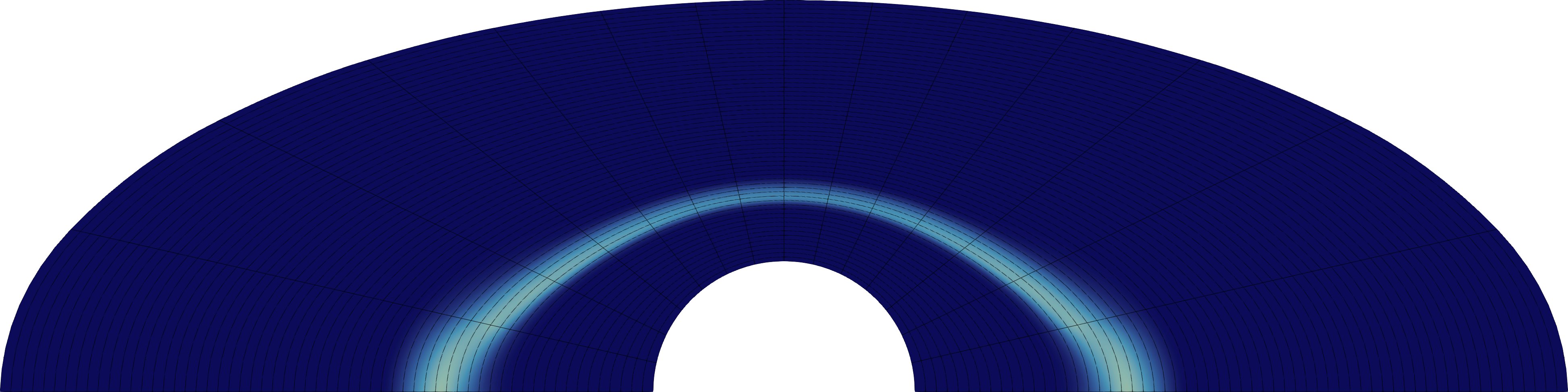}};
        \node[cell] (D) at ([xshift=0.15cm]C.north east) {\includegraphics[width=0.48\linewidth]{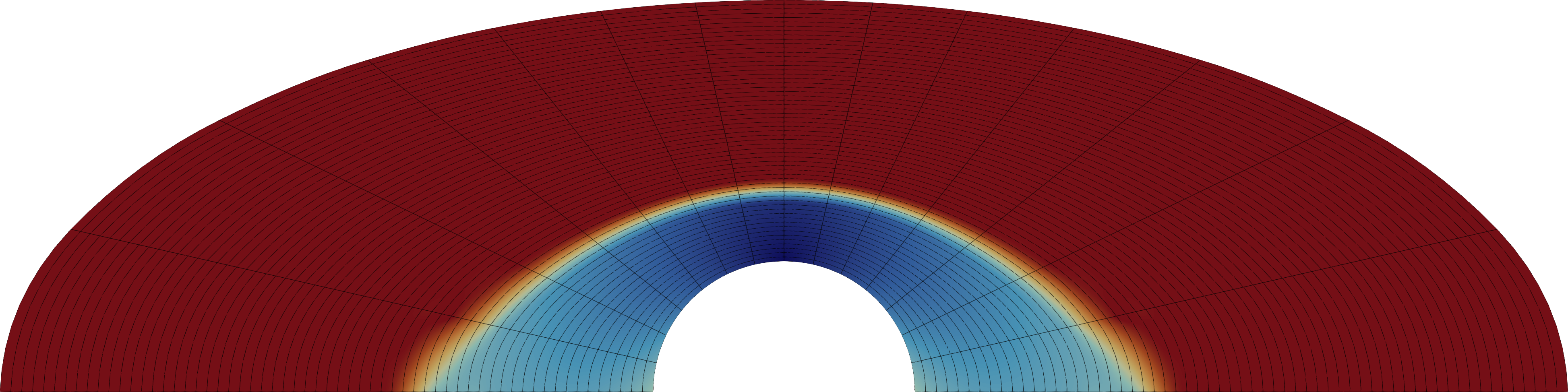}};
        \node[cell] (E) at ([yshift=-0.15cm]C.south west) {\includegraphics[width=0.48\linewidth]{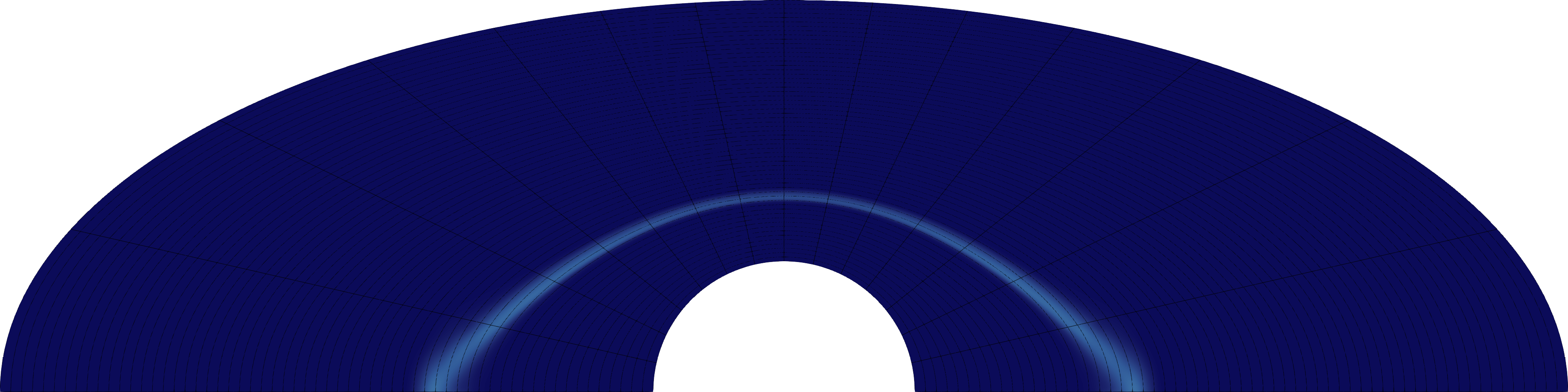}};
        \node[cell] (F) at ([xshift=0.15cm]E.north east) {\includegraphics[width=0.48\linewidth]{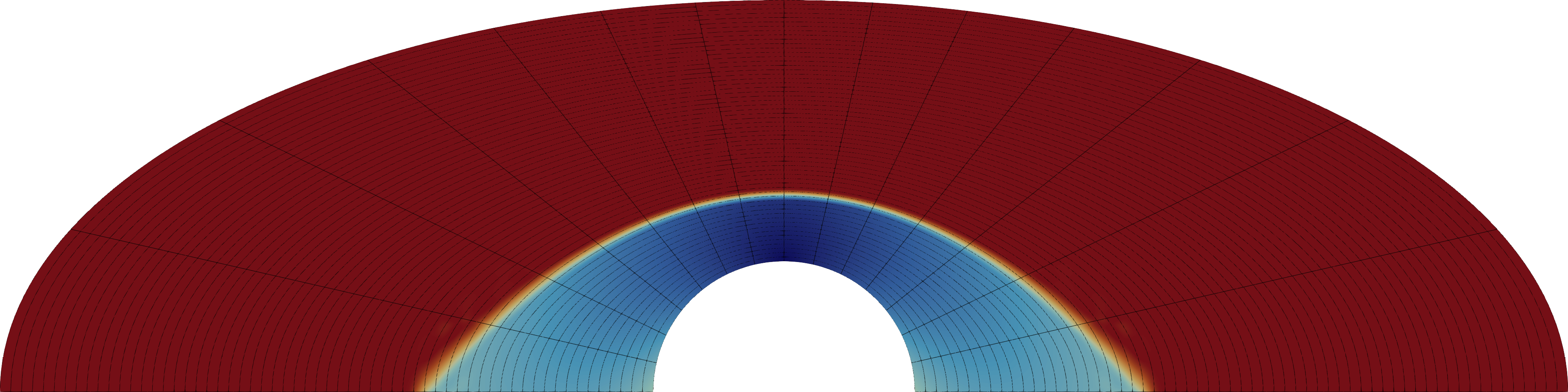}};
        \node[plabel] at (A.north west) {(a)};
        \node[plabel] at (B.north west) {(b)};
        \node[plabel] at (C.north west) {(c)};
        \node[plabel] at (D.north west) {(d)};
        \node[plabel] at (E.north west) {(e)};
        \node[plabel] at (F.north west) {(f)};
    \end{tikzpicture}
    \caption{\case{cyl-euler-hp}: AV and Mach fields at three $(\bar h, p)$ combinations. (a,c,e)~AV $\nu$, on a common color scale from $0$ to $0.477$; (b,d,f)~Mach number. Rows top to bottom: (a,b)~$(p, n_{\mathrm{r}}) = (3, 20)$ at $\bar h = 1/30$; (c,d)~$(p, n_{\mathrm{r}}) = (2, 60)$ at $\bar h = 1/60$; (e,f)~$(p, n_{\mathrm{r}}) = (4, 60)$ at $\bar h = 1/120$.}
    \label{fig:cyl_2d_refinement}
\end{figure}
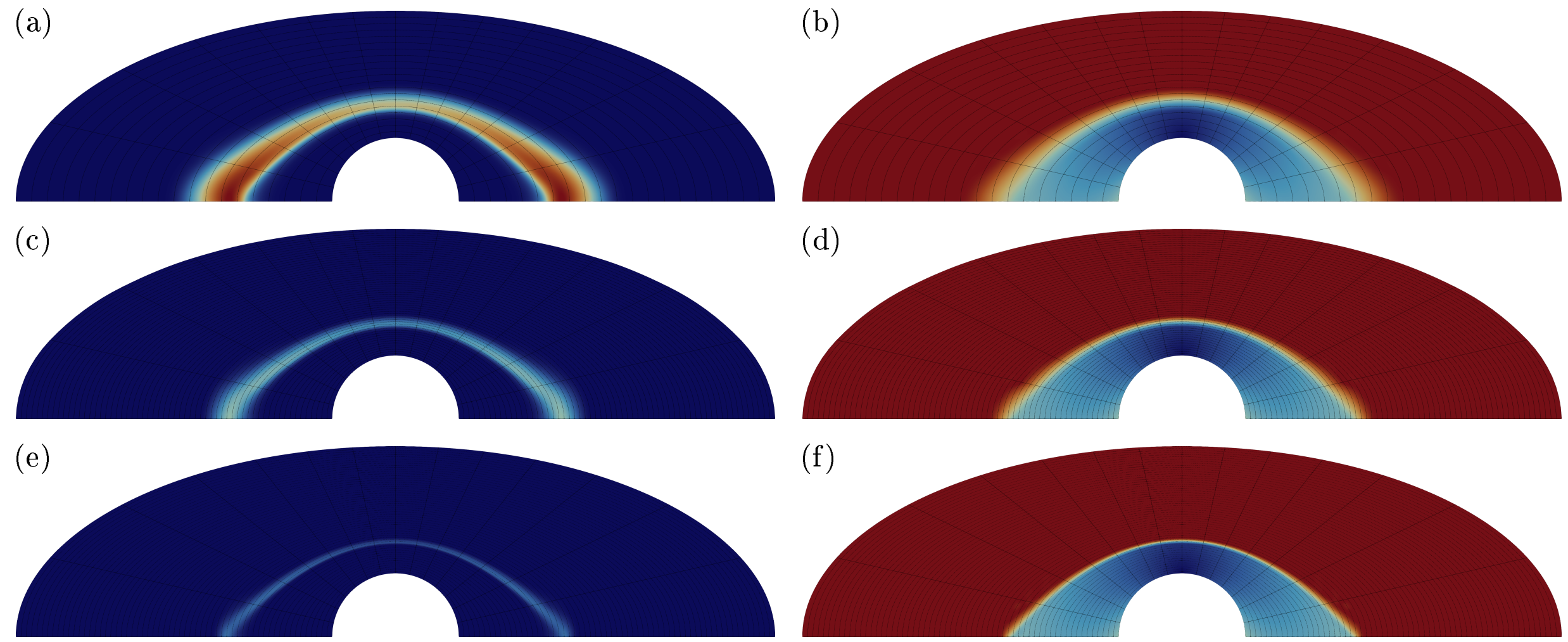

\section{Width--magnitude reparametrization of the decay coefficient}
\label{sec:reparam}
\label{subsec:c4_reparam}

This section justifies the choice $c_4 = 1/c_2$ of Section~\ref{subsec:scpde} through a Green's-function analysis of the steady AV PDE (Section~\ref{subsec:greens}) and verifies its consequences on the 1D Burgers' equation and on inviscid and viscous flow over a half-cylinder (Section~\ref{subsec:reparam_verify}).

\subsection{Green's-function analysis}
\label{subsec:greens}

\begin{figure}[hbt!]
    \centering
    \IfFileExists{_tikz/widening_puzzle.tex}{\begin{tikzpicture}[baseline=(current bounding box.north)]

\begin{axis}[
    name=left,
    width=6.7cm, height=4.8cm,
    xmin=-6, xmax=6,
    ymin=0, ymax=0.6,
    xlabel={$\xi$}, ylabel={$\tilde{\nu}$},
    tick label style={font=\normalsize},
    title={(a) Standard ($c_4 = 1$)},
    grid=both,
    grid style={gray!20},
    legend style={draw=none, fill=white, fill opacity=0.85, text opacity=1, at={(0.97,0.97)}, anchor=north east, inner sep=2pt, row sep=-1pt},
    legend cell align=left,
]
\addplot[gray, dashed, thin, forget plot, domain=-6:6, samples=2] {0.5};
\addplot[gray, dashed, thin, forget plot, domain=-6:6, samples=2] {0.25};
\addplot[red, thick, domain=-6:6, samples=200] {0.5*exp(-abs(x))};
\addlegendentry{$c_2 = 1$}
\addplot[blue, thick, domain=-6:6, samples=200] {0.25*exp(-abs(x)/2)};
\addlegendentry{$c_2 = 4$}
\end{axis}

\begin{axis}[
    at={(left.east)}, anchor=west, xshift=1.0cm,
    width=6.7cm, height=4.8cm,
    xmin=-6, xmax=6,
    ymin=0, ymax=0.6,
    xlabel={$\xi$}, yticklabels={},
    tick label style={font=\normalsize},
    title={(b) Derived ($c_4 = 1/c_2$)},
    grid=both,
    grid style={gray!20},
    legend style={draw=none, fill=white, fill opacity=0.85, text opacity=1, at={(0.97,0.97)}, anchor=north east, inner sep=2pt, row sep=-1pt},
    legend cell align=left,
]
\addplot[gray, dashed, thin, forget plot, domain=-6:6, samples=2] {0.5};
\addplot[red, thick, domain=-6:6, samples=200] {0.5*exp(-abs(x))};
\addlegendentry{$c_2 = 1$}
\addplot[blue, thick, domain=-6:6, samples=200] {0.5*exp(-abs(x)/4)};
\addlegendentry{$c_2 = 4$}
\end{axis}

\end{tikzpicture}}{\fbox{\parbox[c][4.5cm][c]{8cm}{\centering\scriptsize widening puzzle: two side-by-side panels, AV curves vs $\xi$ for $c_2 \in \{1, 4\}$; left standard, right derived}}}
    \caption{$\tilde\nu(\xi)$ at $c_3 = 1$ for $c_2 \in \{1, 4\}$. (a)~standard~\eqref{eq:greens_std}; (b)~derived~\eqref{eq:greens_reparam}. Dashed horizontals: peak heights.}
    \label{fig:widening_puzzle}
\end{figure}
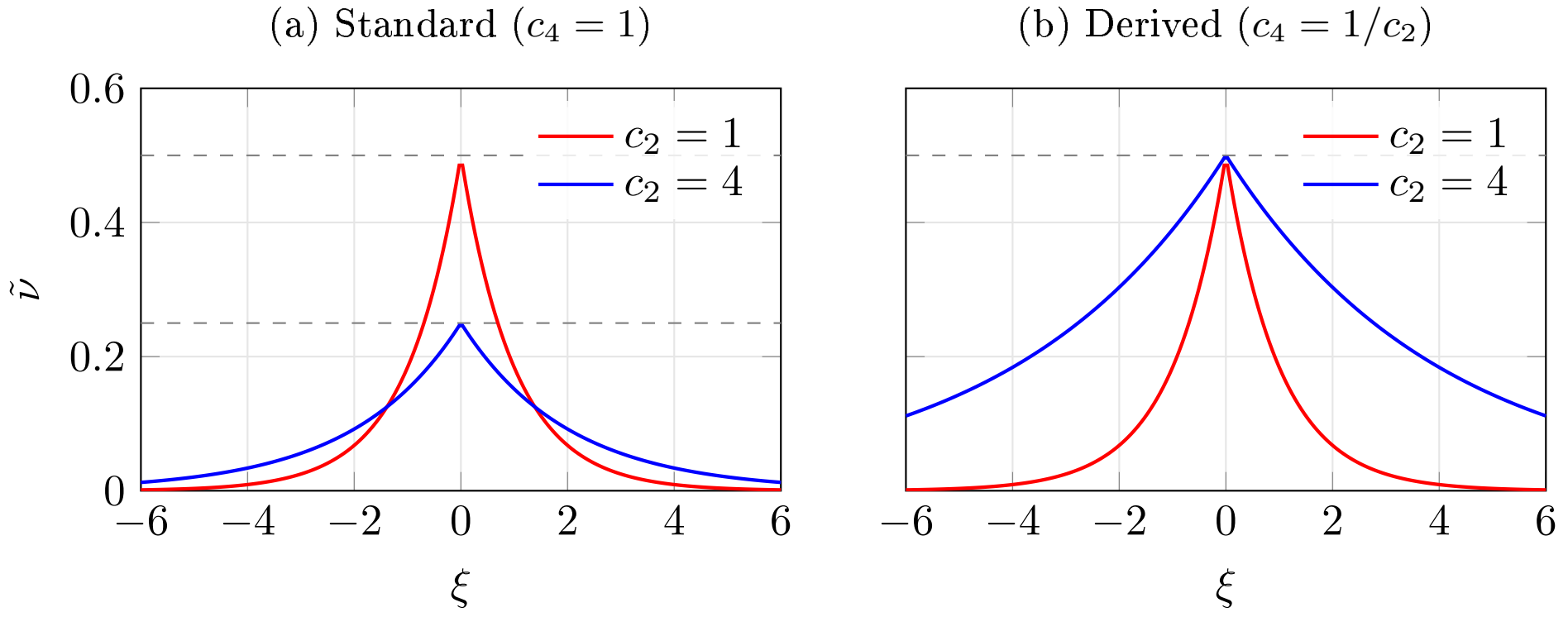

In the recast notation of~\eqref{eq:avpde}, the prior formulations~\cite{barter2010shock, yano2012optimization, nguyen2024optimal} take $c_4 = 1$, so the decay term scales independently of the diffusion coefficient $c_2$; we call this the standard parametrization. Take the steady 1D AV PDE in the rescaled coordinate $\xi = x/\bar h$ of Section~\ref{subsubsec:hbar2}, in which the $\bar h^2$-scaling is absorbed and only $c_2$, $c_3$, $c_4$ remain,
\begin{equation}
    -c_2 \, \frac{\partial^2 \tilde\nu}{\partial \xi^2} = c_3 \, \tilde s - c_4 \, \tilde\nu.
    \label{eq:avpde_1d_linearized}
\end{equation}
For a unit-impulse sensor $\tilde s(\xi) = \delta_{\mathrm{D}}(\xi)$, with $\delta_{\mathrm{D}}$ the Dirac delta, the Green's function is
\begin{equation}
    \tilde\nu(\xi) = \frac{c_3}{2 \sqrt{c_2 \, c_4}} \, \exp\!\left( -\frac{|\xi|}{\sqrt{c_2/c_4}} \right),
    \label{eq:greens}
\end{equation}
with peak magnitude $\tilde\nu_{\mathrm{peak}} = c_3 / (2 \sqrt{c_2 \, c_4})$ and $e$-folding half-width $\ell = \sqrt{c_2/c_4}$, the distance from the peak at which $\tilde\nu$ decays by a factor of $e$. In the standard parametrization, peak and width reduce to
\begin{equation}
    \tilde\nu_{\mathrm{peak}} = \frac{c_3}{2 \sqrt{c_2}}, \qquad \ell = \sqrt{c_2},
    \label{eq:greens_std}
\end{equation}
and are coupled through $c_2$: their product $\tilde\nu_{\mathrm{peak}} \, \ell = c_3/2$ is independent of $c_2$, so widening the AV by raising $c_2$ lowers the peak by the same factor. Setting $c_4 = 1/c_2$ in~\eqref{eq:greens} removes the coupling,
\begin{equation}
    \tilde\nu_{\mathrm{peak}} = \frac{c_3}{2}, \qquad \ell = c_2,
    \label{eq:greens_reparam}
\end{equation}
so that $c_3$ sets the magnitude and $c_2$ the width; we call this the derived parametrization. For $c_3 = 1$, quadrupling $c_2$ from $1$ to $4$ doubles the width and halves the peak under the standard parametrization, whereas it quadruples the width at constant peak under the derived one (Fig.~\ref{fig:widening_puzzle}).

\begin{remark}
The standard parametrization is benign in most cases. It becomes a concern when widening the AV by raising $c_2$ drops the peak below what the local sensor demands, so that oscillations re-emerge in the captured shock. Barter and Darmofal~\cite{barter2010shock} use a Green's-function analysis of the same model to inform the boundary condition of the AV PDE but retain $c_4 = 1$.
\end{remark}

\begin{remark}
The decoupling~\eqref{eq:greens_reparam} holds exactly only for the 1D Green's function. In the coupled state--AV problem the sensor responds to the state, which depends on the AV, so the peak still varies with $c_2$ under the derived parametrization. Section~\ref{subsec:reparam_verify} shows that the residual dependence has the opposite sign to the standard one: the peak grows with $c_2$ rather than shrinks.
\end{remark}

\subsection{Verification}
\label{subsec:reparam_verify}

\subsubsection{1D Burgers' equation}
\label{subsubsec:c4_test_burgers}

The purpose of this experiment (test case: \case{burgers-c4}) is to test two claims of Section~\ref{subsec:greens} on the coupled state--AV problem: that the derived parametrization removes the width--magnitude coupling, as the Green's function predicts, and that the residual $c_2$-dependence of the peak has the opposite sign to the standard one. Consider the Burgers' problem of Section~\ref{subsubsec:selfsim_burgers} with $p = 5$ on $n_\mathrm{el} = 400$ elements, $\bar h = 0.04$. The approximation space is fine enough to resolve the AV-regularized shock without spurious oscillations, so differences between the parametrizations reflect their behavior rather than under-resolution. For both parametrizations we sweep a logarithmically spaced grid in $c_2 \in [1, 50]$ and $c_3 \in [2, 150]$ and extract from each converged solution the peak AV magnitude $\nu_\mathrm{peak} = \max_x \nu(x)$ and the $e$-folding half-width $\ell$.

\begin{figure}[hbt!]
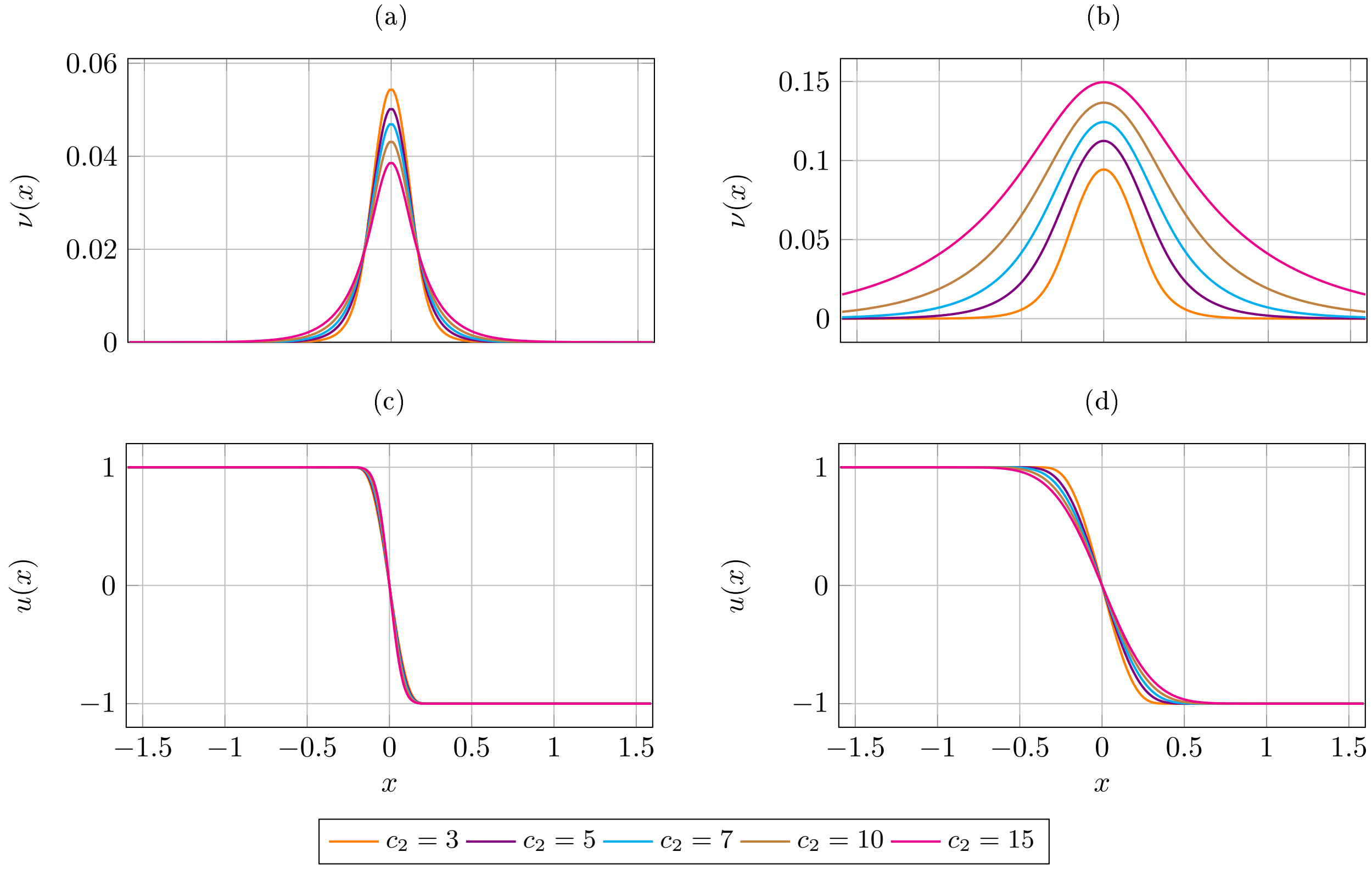

    \centering
    \begin{subfigure}[t]{0.48\textwidth}
        \centering
        \IfFileExists{_tikz/av_profile_std.tex}{\input{_tikz/av_profile_std.tex}}{\fbox{\parbox[c][3.8cm][c]{0.95\linewidth}{\centering\scriptsize AV, standard ($c_4=1$)}}}    \end{subfigure}\hfill
    \begin{subfigure}[t]{0.48\textwidth}
        \centering
        \IfFileExists{_tikz/av_profile_der.tex}{\input{_tikz/av_profile_der.tex}}{\fbox{\parbox[c][3.8cm][c]{0.95\linewidth}{\centering\scriptsize AV, derived ($c_4=1/c_2$)}}}    \end{subfigure}
    \par\vspace{0.4em}
    \begin{subfigure}[t]{0.48\textwidth}
        \centering
        \IfFileExists{_tikz/q_profile_std.tex}{\input{_tikz/q_profile_std.tex}}{\fbox{\parbox[c][3.8cm][c]{0.95\linewidth}{\centering\scriptsize $U$, standard ($c_4=1$)}}}    \end{subfigure}\hfill
    \begin{subfigure}[t]{0.48\textwidth}
        \centering
        \IfFileExists{_tikz/q_profile_der.tex}{\input{_tikz/q_profile_der.tex}}{\fbox{\parbox[c][3.8cm][c]{0.95\linewidth}{\centering\scriptsize $U$, derived ($c_4=1/c_2$)}}}    \end{subfigure}
    \figurelegend{xspaceprofileslegend}
    \caption{\case{burgers-c4}: converged (a,b)~AV $\nu(x)$ and (c,d)~state $u(x)$ across the $c_2$ sweep at fixed $c_3 = 5$, for the (a,c)~standard ($c_4 = 1$) and (b,d)~derived ($c_4 = 1/c_2$) parametrizations.}
    \label{fig:c4_profiles_x}
\end{figure}

Along the $c_2$ sweep at fixed $c_3 = 5$, raising $c_2$ under the standard parametrization widens the AV and lowers its peak, and the shock thickness in $u$ changes only mildly because the two effects partially compensate (Fig.~\ref{fig:c4_profiles_x}). Under the derived parametrization the Green's-function prediction of a peak independent of $c_2$ does not survive the state--AV coupling: raising $c_2$ widens the AV \emph{and} grows its peak, so the shock thickness in $u$ broadens visibly across the sweep.

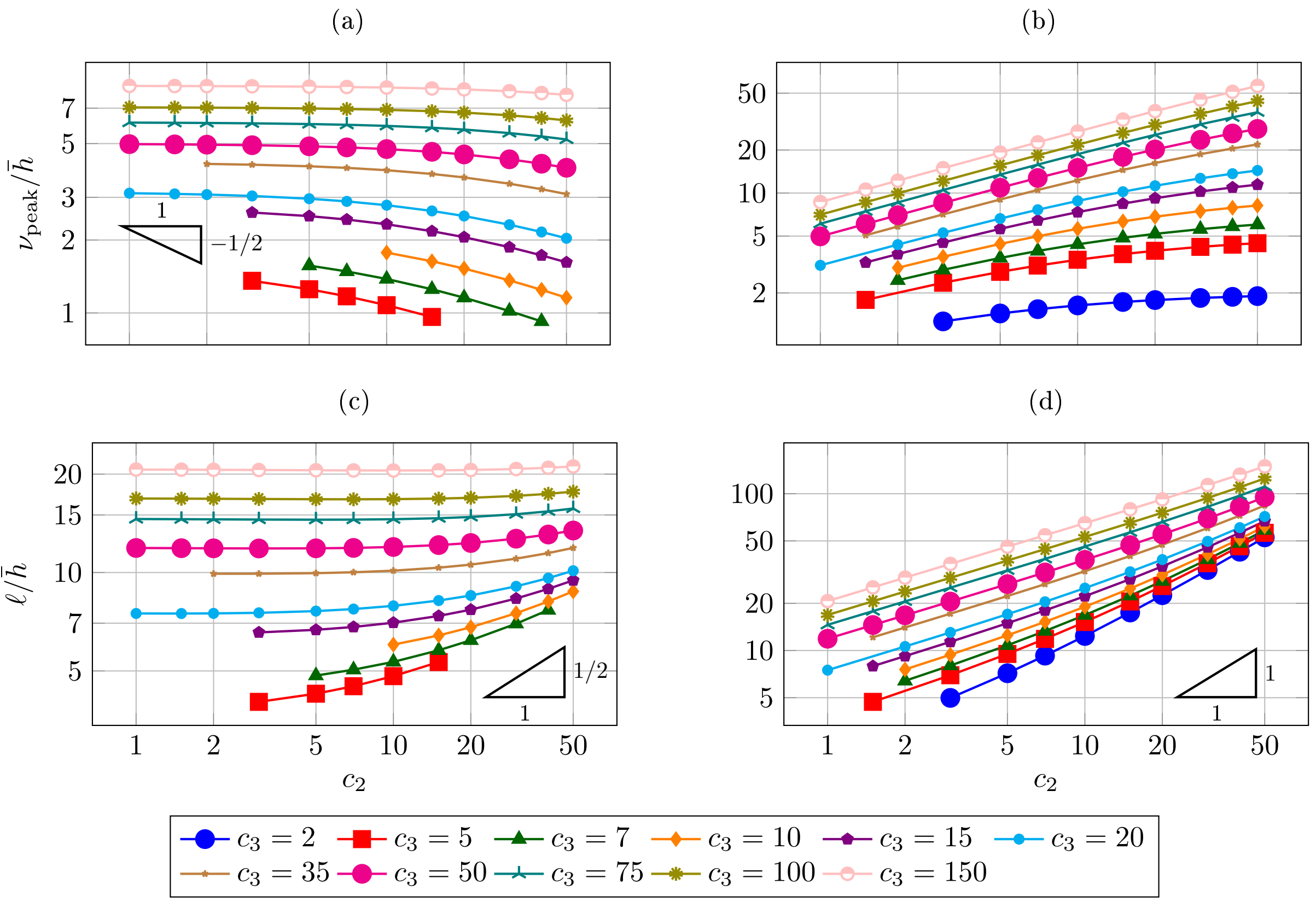
\begin{figure}[hbt!]
    \centering
    \begin{subfigure}[t]{0.48\textwidth}
        \centering
        \begin{tikzpicture}[baseline=(current axis.north)]
\begin{axis}[title={(a)},
    tick label style={font=\normalsize},
    ylabel={$\nu_{\mathrm{peak}} / \bar{h}$},
    grid=both, width=\linewidth, height=5cm,
    xmode=log, ymode=log,
    xtick={1,2,5,10,20,50}, xticklabel=\empty,
    ytick={1,2,3,5,7},
    log ticks with fixed point,
]
\addplot[red, mark=square*, mark size=2.7pt, thick] coordinates {
    (3.0, 1.3573742728685891)
    (5.0, 1.252993994699859)
    (7.0, 1.171559836810747)
    (10.0, 1.0768276688109726)
    (15.0, 0.9629271938957602)
};
\addplot[black!60!green, mark=triangle*, mark size=2.4pt, thick] coordinates {
    (5.0, 1.5688247022562312)
    (7.0, 1.483899032720856)
    (10.0, 1.3812523683890245)
    (15.0, 1.2525363158640905)
    (20.0, 1.1563409224122223)
    (30.0, 1.018519813365892)
    (40.0, 0.9231417057713385)
};
\addplot[orange, mark=diamond*, mark size=2.1pt, thick] coordinates {
    (10.0, 1.7748517059869944)
    (15.0, 1.6341396911917803)
    (20.0, 1.524972906373158)
    (30.0, 1.363100136813772)
    (40.0, 1.2460652333736042)
    (50.0, 1.1558645245818309)
};
\addplot[violet, mark=pentagon*, mark size=1.8pt, thick] coordinates {
    (3.0, 2.5973926521650084)
    (5.0, 2.5096239538931773)
    (7.0, 2.4293449993205276)
    (10.0, 2.32269649359332)
    (15.0, 2.174917257004754)
    (20.0, 2.054772224650667)
    (30.0, 1.8686781672608337)
    (40.0, 1.7286531960615028)
    (50.0, 1.617757605712653)
};
\addplot[cyan, mark=otimes*, mark size=1.5pt, thick] coordinates {
    (1.0, 3.1202299655917076)
    (1.5, 3.1001545016434324)
    (2.0, 3.0800081305031553)
    (3.0, 3.039755643451216)
    (5.0, 2.9610475962341725)
    (7.0, 2.886562123512525)
    (10.0, 2.7841361643621325)
    (15.0, 2.636441629299225)
    (20.0, 2.512080003618733)
    (30.0, 2.3129006240843117)
    (40.0, 2.1583457818162035)
    (50.0, 2.033326858299125)
};
\addplot[brown, mark=star, mark size=1.2pt, thick] coordinates {
    (2.0, 4.1213916353364235)
    (3.0, 4.090999821573109)
    (5.0, 4.0301089523937055)
    (7.0, 3.9699091395436605)
    (10.0, 3.882496009193475)
    (15.0, 3.747043611132273)
    (20.0, 3.6248705483863635)
    (30.0, 3.4150763096597405)
    (40.0, 3.241288999334673)
    (50.0, 3.0941692526453015)
};
\addplot[magenta, mark=oplus*, mark size=3pt, thick] coordinates {
    (1.0, 4.973512854858545)
    (1.5, 4.960930397894087)
    (2.0, 4.948316251680633)
    (3.0, 4.9229970442973165)
    (5.0, 4.87209437873437)
    (7.0, 4.821156379191395)
    (10.0, 4.7455794005679355)
    (15.0, 4.624192095709037)
    (20.0, 4.51028551270734)
    (30.0, 4.305740538150141)
    (40.0, 4.1286581123573605)
    (50.0, 3.973914899356443)
};
\addplot[teal, mark=Mercedes star, mark size=2.7pt, thick] coordinates {
    (1.0, 6.102116789564862)
    (1.5, 6.091865831453676)
    (2.0, 6.081597606104552)
    (3.0, 6.061009356710306)
    (5.0, 6.01963745759735)
    (7.0, 5.978071995044049)
    (10.0, 5.915661122558196)
    (15.0, 5.812724579633928)
    (20.0, 5.712668549072003)
    (30.0, 5.524544034374356)
    (40.0, 5.353472641089845)
    (50.0, 5.198293138710734)
};
\addplot[olive, mark=10-pointed star, mark size=2.4pt, thick] coordinates {
    (1.0, 7.052363357685022)
    (1.5, 7.043495504429208)
    (2.0, 7.03461647340142)
    (3.0, 7.016824741765711)
    (5.0, 6.9811078537997515)
    (7.0, 6.945228687603054)
    (10.0, 6.891208569537071)
    (15.0, 6.801233870153367)
    (20.0, 6.7123029550358755)
    (30.0, 6.540635135296715)
    (40.0, 6.379537740629631)
    (50.0, 6.229581733774812)
};
\addplot[pink, mark=halfcircle*, mark size=2.1pt, thick] coordinates {
    (1.0, 8.644988636429133)
    (1.5, 8.637755988832064)
    (2.0, 8.630517280299989)
    (3.0, 8.616021618290269)
    (5.0, 8.586957007259722)
    (7.0, 8.557795414797523)
    (10.0, 8.513884965660843)
    (15.0, 8.440391534639614)
    (20.0, 8.366851700525354)
    (30.0, 8.221283473599053)
    (40.0, 8.079783854236458)
    (50.0, 7.9437485650738076)
};
\logLogSlopeTriangle[south]{0.22}{0.15}{0.42}{-1/2}{thick, black}
\end{axis}
\end{tikzpicture}    \end{subfigure}\hfill
    \begin{subfigure}[t]{0.48\textwidth}
        \centering
        \begin{tikzpicture}[baseline=(current axis.north)]
\begin{axis}[title={(b)},
    tick label style={font=\normalsize},
    grid=both, width=\linewidth, height=5cm,
    xmode=log, ymode=log,
    xtick={1,2,5,10,20,50}, xticklabel=\empty,
    ytick={2,5,10,20,50},
    log ticks with fixed point,
    legend to name=decouplinglegend,
    legend columns=6,
    legend cell align=left,
]
\addplot[blue, mark=*, mark size=3pt, thick] coordinates {
    (3.0, 1.2645482024079886)
    (5.0, 1.436680371148381)
    (7.0, 1.5392694998928798)
    (10.0, 1.6357489299186545)
    (15.0, 1.7282924844228986)
    (20.0, 1.7826928066751586)
    (30.0, 1.844389372623217)
    (40.0, 1.878639046376772)
    (50.0, 1.9004826457709014)
};
\addplot[red, mark=square*, mark size=2.7pt, thick] coordinates {
    (1.5, 1.7876426715540128)
    (3.0, 2.3565351736635836)
    (5.0, 2.8090521254716165)
    (7.0, 3.107884446472889)
    (10.0, 3.414303273019182)
    (15.0, 3.737949038648159)
    (20.0, 3.945642495111879)
    (30.0, 4.201422908587226)
    (40.0, 4.354982315812034)
    (50.0, 4.458156478745238)
};
\addplot[black!60!green, mark=triangle*, mark size=2.4pt, thick] coordinates {
    (2.0, 2.4538425203022185)
    (3.0, 2.9015531362655134)
    (5.0, 3.514753303517521)
    (7.0, 3.9339510805010356)
    (10.0, 4.376842880842602)
    (15.0, 4.860886959454838)
    (20.0, 5.181671579784621)
    (30.0, 5.58964518176149)
    (40.0, 5.842602389848479)
    (50.0, 6.016477243387379)
};
\addplot[orange, mark=diamond*, mark size=2.1pt, thick] coordinates {
    (2.0, 2.9999755489738487)
    (3.0, 3.5811280390186555)
    (5.0, 4.404260685932946)
    (7.0, 4.986659440621014)
    (10.0, 5.6208928723035445)
    (15.0, 6.338645758313274)
    (20.0, 6.830300610929157)
    (30.0, 7.477089301805434)
    (40.0, 7.892214379296353)
    (50.0, 8.184774907729494)
};
\addplot[violet, mark=pentagon*, mark size=1.8pt, thick] coordinates {
    (1.5, 3.2671118440579017)
    (2.0, 3.7402627232065893)
    (3.0, 4.501897196667544)
    (5.0, 5.616479366954341)
    (7.0, 6.433447418195683)
    (10.0, 7.3523039894280116)
    (15.0, 8.432120159913763)
    (20.0, 9.198819311163227)
    (30.0, 10.245873245547312)
    (40.0, 10.944337254961486)
    (50.0, 11.450995611795564)
};
\addplot[cyan, mark=otimes*, mark size=1.5pt, thick] coordinates {
    (1.0, 3.1202299655917076)
    (2.0, 4.357409192666132)
    (3.0, 5.267678746184039)
    (5.0, 6.625263871907663)
    (7.0, 7.642377767999825)
    (10.0, 8.810678564712635)
    (15.0, 10.218765570844655)
    (20.0, 11.243176774862585)
    (30.0, 12.678327551446191)
    (40.0, 13.661428436344863)
    (50.0, 14.389094306458935)
};
\addplot[brown, mark=star, mark size=1.2pt, thick] coordinates {
    (1.5, 5.066886268091055)
    (2.0, 5.829742559263257)
    (3.0, 7.087816479031267)
    (5.0, 9.014724546174786)
    (7.0, 10.507381579465065)
    (10.0, 12.282538784792244)
    (15.0, 14.518522514538724)
    (20.0, 16.21813254868313)
    (30.0, 18.71368746896329)
    (40.0, 20.50912845781807)
    (50.0, 21.889140587882505)
};
\addplot[magenta, mark=oplus*, mark size=3pt, thick] coordinates {
    (1.0, 4.973512854858545)
    (1.5, 6.076457003742436)
    (2.0, 6.998988188579216)
    (3.0, 8.528542244262223)
    (5.0, 10.89693154252029)
    (7.0, 12.758897658206811)
    (10.0, 15.011027822992734)
    (15.0, 17.914740894797987)
    (20.0, 20.176805567453215)
    (30.0, 23.591009609342755)
    (40.0, 26.12033266814041)
    (50.0, 28.10891756969233)
};
\addplot[teal, mark=Mercedes star, mark size=2.7pt, thick] coordinates {
    (1.0, 6.102116789564862)
    (1.5, 7.46145620515291)
    (2.0, 8.601501541757695)
    (3.0, 10.499325842487178)
    (5.0, 13.462423378250602)
    (7.0, 15.819176475914052)
    (10.0, 18.710358575462916)
    (15.0, 22.51693079596354)
    (20.0, 25.552941940804875)
    (30.0, 30.265483315028778)
    (40.0, 33.86553751934716)
    (50.0, 36.765392299492504)
};
\addplot[olive, mark=10-pointed star, mark size=2.4pt, thick] coordinates {
    (1.0, 7.052363357685022)
    (1.5, 8.626895642081807)
    (2.0, 9.949162143284697)
    (3.0, 12.1546628745509)
    (5.0, 15.612047004808343)
    (7.0, 18.37766077251833)
    (10.0, 21.79483737740811)
    (15.0, 26.344810373751404)
    (20.0, 30.022753715483933)
    (30.0, 35.83004586019951)
    (40.0, 40.35409966211025)
    (50.0, 44.05684373203306)
};
\addplot[pink, mark=halfcircle*, mark size=2.1pt, thick] coordinates {
    (1.0, 8.644988636429133)
    (1.5, 10.579382217772968)
    (2.0, 12.205975098353292)
    (3.0, 14.924336763692395)
    (5.0, 19.202495489803)
    (7.0, 22.643675257825954)
    (10.0, 26.92563629168914)
    (15.0, 32.69252643809543)
    (20.0, 37.42128311609802)
    (30.0, 45.034324287059555)
    (40.0, 51.10628707252677)
    (50.0, 56.1766606346151)
};
\legend{$c_3 = 2$, $c_3 = 5$, $c_3 = 7$, $c_3 = 10$, $c_3 = 15$, $c_3 = 20$, $c_3 = 35$, $c_3 = 50$, $c_3 = 75$, $c_3 = 100$, $c_3 = 150$}
\end{axis}
\end{tikzpicture}    \end{subfigure}
    \par\vspace{0.4em}
    \begin{subfigure}[t]{0.48\textwidth}
        \centering
        \begin{tikzpicture}[baseline=(current axis.north)]
\begin{axis}[title={(c)},
    tick label style={font=\normalsize},
    xlabel={$c_2$}, ylabel={$\ell / \bar{h}$},
    grid=both, width=\linewidth, height=5cm,
    xmode=log, ymode=log,
    xtick={1,2,5,10,20,50},
    ytick={5,7,10,15,20},
    log ticks with fixed point,
]
\addplot[red, mark=square*, mark size=2.7pt, thick] coordinates {
    (3.0, 4.021567551166154)
    (5.0, 4.259158154266508)
    (7.0, 4.489825599484887)
    (10.0, 4.819742274512861)
    (15.0, 5.322596377796446)
};
\addplot[black!60!green, mark=triangle*, mark size=2.4pt, thick] coordinates {
    (5.0, 4.839356465883794)
    (7.0, 5.044340436560297)
    (10.0, 5.329609340788599)
    (15.0, 5.784019318486651)
    (20.0, 6.205375337249956)
    (30.0, 6.974911273589687)
    (40.0, 7.680031017768221)
};
\addplot[orange, mark=diamond*, mark size=2.1pt, thick] coordinates {
    (10.0, 6.018927091290619)
    (15.0, 6.4204376348757695)
    (20.0, 6.802772425321135)
    (30.0, 7.510128910113822)
    (40.0, 8.157451105215433)
    (50.0, 8.756676752560075)
};
\addplot[violet, mark=pentagon*, mark size=1.8pt, thick] coordinates {
    (3.0, 6.559309900111987)
    (5.0, 6.6795519971317106)
    (7.0, 6.815293263816177)
    (10.0, 7.024039528374355)
    (15.0, 7.3667949264483195)
    (20.0, 7.699674792162109)
    (30.0, 8.326805641463999)
    (40.0, 8.91413588078179)
    (50.0, 9.464681365549305)
};
\addplot[cyan, mark=otimes*, mark size=1.5pt, thick] coordinates {
    (1.0, 7.494969088546916)
    (1.5, 7.492947445537763)
    (2.0, 7.500110077168362)
    (3.0, 7.530189543000934)
    (5.0, 7.6213060323695485)
    (7.0, 7.731901967803252)
    (10.0, 7.9082553600518715)
    (15.0, 8.20931742905164)
    (20.0, 8.505532128219295)
    (30.0, 9.07472268702926)
    (40.0, 9.614186617401538)
    (50.0, 10.127595012707092)
};
\addplot[brown, mark=star, mark size=1.2pt, thick] coordinates {
    (2.0, 9.913812304503807)
    (3.0, 9.916120780683238)
    (5.0, 9.955209710350582)
    (7.0, 10.017850131915488)
    (10.0, 10.132247170150809)
    (15.0, 10.34691070093259)
    (20.0, 10.57278588240701)
    (30.0, 11.02602754481105)
    (40.0, 11.468299175522379)
    (50.0, 11.896583199219137)
};
\addplot[magenta, mark=oplus*, mark size=3pt, thick] coordinates {
    (1.0, 11.887990146873143)
    (1.5, 11.870160859832966)
    (2.0, 11.85816042226253)
    (3.0, 11.847127027510457)
    (5.0, 11.857772538417358)
    (7.0, 11.89302188269424)
    (10.0, 11.970827752832507)
    (15.0, 12.132270182254759)
    (20.0, 12.312640677118223)
    (30.0, 12.690043236521662)
    (40.0, 13.069421665415643)
    (50.0, 13.443851186840588)
};
\addplot[teal, mark=Mercedes star, mark size=2.7pt, thick] coordinates {
    (1.0, 14.57830700320956)
    (1.5, 14.559125968634836)
    (2.0, 14.543710840863971)
    (3.0, 14.522593843284648)
    (5.0, 14.508187209209197)
    (7.0, 14.517293674650986)
    (10.0, 14.556779090349242)
    (15.0, 14.66051436873485)
    (20.0, 14.788978887445975)
    (30.0, 15.08021325319856)
    (40.0, 15.385748949404103)
    (50.0, 15.696053238315843)
};
\addplot[olive, mark=10-pointed star, mark size=2.4pt, thick] coordinates {
    (1.0, 16.84577120021239)
    (1.5, 16.82721673116812)
    (2.0, 16.811050396498125)
    (3.0, 16.78590320953534)
    (5.0, 16.758585764989295)
    (7.0, 16.752527843545185)
    (10.0, 16.76848532649843)
    (15.0, 16.834323085788625)
    (20.0, 16.927846078660973)
    (30.0, 17.15637312701666)
    (40.0, 17.41087849906305)
    (50.0, 17.6758890472554)
};
\addplot[pink, mark=halfcircle*, mark size=2.1pt, thick] coordinates {
    (1.0, 20.648070888693486)
    (1.5, 20.631534289358992)
    (2.0, 20.616044411987694)
    (3.0, 20.588969947742395)
    (5.0, 20.55017650098233)
    (7.0, 20.528040949065346)
    (10.0, 20.51695877170457)
    (15.0, 20.536599055288256)
    (20.0, 20.585560020047854)
    (30.0, 20.731802900888457)
    (40.0, 20.913458008271736)
    (50.0, 21.113472686629322)
};
\logLogSlopeTriangle{0.9}{0.15}{0.1}{1/2}{thick, black}
\end{axis}
\end{tikzpicture}    \end{subfigure}\hfill
    \begin{subfigure}[t]{0.48\textwidth}
        \centering
        \begin{tikzpicture}[baseline=(current axis.north)]
\begin{axis}[title={(d)},
    tick label style={font=\normalsize},
    xlabel={$c_2$},
    grid=both, width=\linewidth, height=5cm,
    xmode=log, ymode=log,
    xtick={1,2,5,10,20,50},
    ytick={5,10,20,50,100},
    log ticks with fixed point,
]
\addplot[blue, mark=*, mark size=3pt, thick] coordinates {
    (3.0, 4.985881801738666)
    (5.0, 7.163952118332397)
    (7.0, 9.272153610930065)
    (10.0, 12.374183278514387)
    (15.0, 17.473900313888358)
    (20.0, 22.533162914686873)
    (30.0, 32.601186199343914)
    (40.0, 42.63969523731348)
    (50.0, 52.66463239622592)
};
\addplot[red, mark=square*, mark size=2.7pt, thick] coordinates {
    (1.5, 4.704258883599129)
    (3.0, 6.953944398852188)
    (5.0, 9.507753884296982)
    (7.0, 11.854816565434756)
    (10.0, 15.195583297799564)
    (15.0, 20.542794141732053)
    (20.0, 25.76345498831848)
    (30.0, 36.036866962146846)
    (40.0, 46.20338515308426)
    (50.0, 56.31714752767015)
};
\addplot[black!60!green, mark=triangle*, mark size=2.4pt, thick] coordinates {
    (2.0, 6.408652820174983)
    (3.0, 8.016338928420328)
    (5.0, 10.800234342596125)
    (7.0, 13.30837041960205)
    (10.0, 16.82227081775046)
    (15.0, 22.357637083645603)
    (20.0, 27.702341812293707)
    (30.0, 38.135826574366284)
    (40.0, 48.40679278628011)
    (50.0, 58.59514808084691)
};
\addplot[orange, mark=diamond*, mark size=2.1pt, thick] coordinates {
    (2.0, 7.5723872001250125)
    (3.0, 9.399465876773332)
    (5.0, 12.500156180547823)
    (7.0, 15.238814019090272)
    (10.0, 19.009987100289866)
    (15.0, 24.84088021481707)
    (20.0, 30.385185460325403)
    (30.0, 41.08030734158949)
    (40.0, 51.52193259109819)
    (50.0, 61.83372520822843)
};
\addplot[violet, mark=pentagon*, mark size=1.8pt, thick] coordinates {
    (1.5, 7.95548852340104)
    (2.0, 9.208335354362614)
    (3.0, 11.356790527075196)
    (5.0, 14.93171467770191)
    (7.0, 18.021273374732345)
    (10.0, 22.19640504946554)
    (15.0, 28.513454494991375)
    (20.0, 34.40719327013601)
    (30.0, 45.575747825808804)
    (40.0, 56.32968054063711)
    (50.0, 66.8687663221697)
};
\addplot[cyan, mark=otimes*, mark size=1.5pt, thick] coordinates {
    (1.0, 7.494969088546916)
    (2.0, 10.606377455971435)
    (3.0, 13.040490648261605)
    (5.0, 17.037521951214085)
    (7.0, 20.447546783533777)
    (10.0, 24.99734459472096)
    (15.0, 31.77698769352638)
    (20.0, 38.01703571794891)
    (30.0, 49.67761302328728)
    (40.0, 60.76738587542227)
    (50.0, 71.55354335067764)
};
\addplot[brown, mark=star, mark size=1.2pt, thick] coordinates {
    (1.5, 12.14828057207144)
    (2.0, 14.018019688231336)
    (3.0, 17.171559059117257)
    (5.0, 22.25386123752075)
    (7.0, 26.49532512147088)
    (10.0, 32.030060425126784)
    (15.0, 40.06313217165813)
    (20.0, 47.26921932047778)
    (30.0, 60.37343954541979)
    (40.0, 72.50965552124346)
    (50.0, 84.09721506295217)
};
\addplot[magenta, mark=oplus*, mark size=3pt, thick] coordinates {
    (1.0, 11.887990146873143)
    (1.5, 14.537303525981072)
    (2.0, 16.768966105641606)
    (3.0, 20.51823697541074)
    (5.0, 26.51222821856266)
    (7.0, 31.462758280398905)
    (10.0, 37.84970508387815)
    (15.0, 46.98185260904267)
    (20.0, 55.05429218265656)
    (30.0, 69.49221319818867)
    (40.0, 82.64249882101524)
    (50.0, 95.0406100179735)
};
\addplot[teal, mark=Mercedes star, mark size=2.7pt, thick] coordinates {
    (1.0, 14.57830700320956)
    (1.5, 17.83084846762904)
    (2.0, 20.567178293912765)
    (3.0, 25.15262273847213)
    (5.0, 32.439457402008024)
    (7.0, 38.40675263186346)
    (10.0, 46.02898612661419)
    (15.0, 56.774296195939044)
    (20.0, 66.13108782661281)
    (30.0, 82.58044315403936)
    (40.0, 97.29485452760704)
    (50.0, 110.97471269827687)
};
\addplot[olive, mark=10-pointed star, mark size=2.4pt, thick] coordinates {
    (1.0, 16.84577120021239)
    (1.5, 20.608773042853258)
    (2.0, 23.77399496588343)
    (3.0, 29.073149501767265)
    (5.0, 37.47166647206958)
    (7.0, 44.320703227620896)
    (10.0, 53.02343901501503)
    (15.0, 65.1948661772602)
    (20.0, 75.69867684336093)
    (30.0, 93.96274142140753)
    (40.0, 110.10714965329535)
    (50.0, 124.97574560787676)
};
\addplot[pink, mark=halfcircle*, mark size=2.1pt, thick] coordinates {
    (1.0, 20.648070888693486)
    (1.5, 25.26822902620769)
    (2.0, 29.155109985992066)
    (3.0, 35.66051550605722)
    (5.0, 45.950474672145226)
    (7.0, 54.31054329400744)
    (10.0, 64.87823506779264)
    (15.0, 79.53489650686647)
    (20.0, 92.05733561353621)
    (30.0, 113.54782471961559)
    (40.0, 132.2616308921598)
    (50.0, 149.2859359978872)
};
\logLogSlopeTriangle{0.9}{0.15}{0.1}{1}{thick, black}
\end{axis}
\end{tikzpicture}    \end{subfigure}
    \figurelegend{decouplinglegend}
    \caption{\case{burgers-c4}: (a,b)~peak $\nu_\mathrm{peak}/\bar h$ and (c,d)~half-width $\ell/\bar h$ of the converged AV versus $c_2$, one curve per $c_3$, for the (a,c)~standard and (b,d)~derived parametrizations. Triangles: slopes given by~\eqref{eq:greens_std} and~\eqref{eq:greens_reparam}; none is drawn in~(b), where~\eqref{eq:greens_reparam} gives a peak independent of $c_2$.}
    \label{fig:c4_test_burgers_c2}
\end{figure}

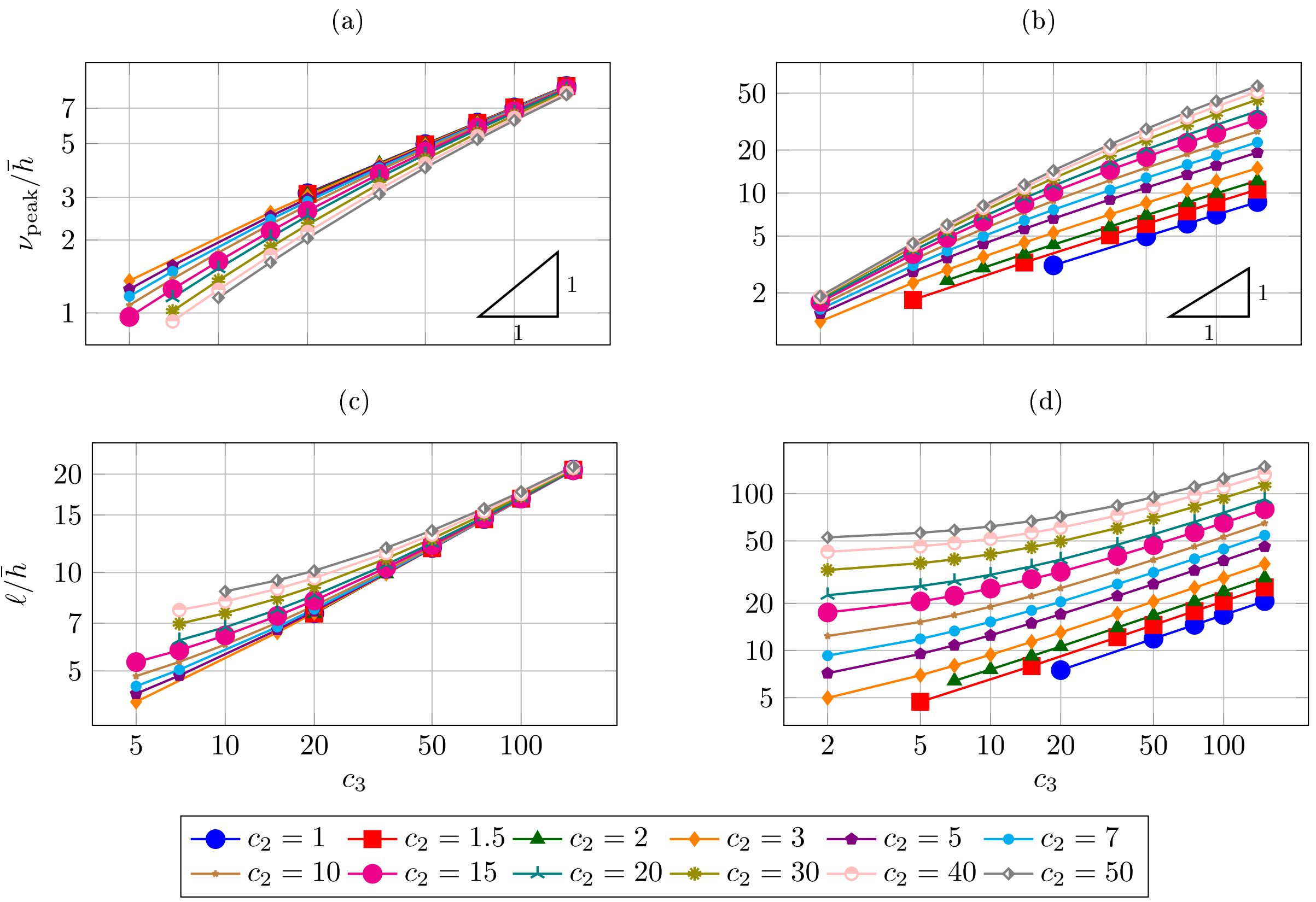
\begin{figure}[hbt!]
    \centering
    \begin{subfigure}[t]{0.48\textwidth}
        \centering
        \begin{tikzpicture}[baseline=(current axis.north)]
\begin{axis}[title={(a)},
    tick label style={font=\normalsize},
    ylabel={$\nu_{\mathrm{peak}} / \bar{h}$},
    grid=both, width=\linewidth, height=5cm,
    xmode=log, ymode=log,
    xtick={1,2,5,10,20,50,100}, xticklabel=\empty,
    ytick={1,2,3,5,7},
    log ticks with fixed point,
]
\addplot[blue, mark=*, mark size=3pt, thick] coordinates {
    (20.0, 3.1202299655917076)
    (50.0, 4.973512854858545)
    (75.0, 6.102116789564862)
    (100.0, 7.052363357685022)
    (150.0, 8.644988636429133)
};
\addplot[red, mark=square*, mark size=2.7pt, thick] coordinates {
    (20.0, 3.1001545016434324)
    (50.0, 4.960930397894087)
    (75.0, 6.091865831453676)
    (100.0, 7.043495504429208)
    (150.0, 8.637755988832064)
};
\addplot[black!60!green, mark=triangle*, mark size=2.4pt, thick] coordinates {
    (20.0, 3.0800081305031553)
    (35.0, 4.1213916353364235)
    (50.0, 4.948316251680633)
    (75.0, 6.081597606104552)
    (100.0, 7.03461647340142)
    (150.0, 8.630517280299989)
};
\addplot[orange, mark=diamond*, mark size=2.1pt, thick] coordinates {
    (5.0, 1.3573742728685891)
    (15.0, 2.5973926521650084)
    (20.0, 3.039755643451216)
    (35.0, 4.090999821573109)
    (50.0, 4.9229970442973165)
    (75.0, 6.061009356710306)
    (100.0, 7.016824741765711)
    (150.0, 8.616021618290269)
};
\addplot[violet, mark=pentagon*, mark size=1.8pt, thick] coordinates {
    (5.0, 1.252993994699859)
    (7.0, 1.5688247022562312)
    (15.0, 2.5096239538931773)
    (20.0, 2.9610475962341725)
    (35.0, 4.0301089523937055)
    (50.0, 4.87209437873437)
    (75.0, 6.01963745759735)
    (100.0, 6.9811078537997515)
    (150.0, 8.586957007259722)
};
\addplot[cyan, mark=otimes*, mark size=1.5pt, thick] coordinates {
    (5.0, 1.171559836810747)
    (7.0, 1.483899032720856)
    (15.0, 2.4293449993205276)
    (20.0, 2.886562123512525)
    (35.0, 3.9699091395436605)
    (50.0, 4.821156379191395)
    (75.0, 5.978071995044049)
    (100.0, 6.945228687603054)
    (150.0, 8.557795414797523)
};
\addplot[brown, mark=star, mark size=1.2pt, thick] coordinates {
    (5.0, 1.0768276688109726)
    (7.0, 1.3812523683890245)
    (10.0, 1.7748517059869944)
    (15.0, 2.32269649359332)
    (20.0, 2.7841361643621325)
    (35.0, 3.882496009193475)
    (50.0, 4.7455794005679355)
    (75.0, 5.915661122558196)
    (100.0, 6.891208569537071)
    (150.0, 8.513884965660843)
};
\addplot[magenta, mark=oplus*, mark size=3pt, thick] coordinates {
    (5.0, 0.9629271938957602)
    (7.0, 1.2525363158640905)
    (10.0, 1.6341396911917803)
    (15.0, 2.174917257004754)
    (20.0, 2.636441629299225)
    (35.0, 3.747043611132273)
    (50.0, 4.624192095709037)
    (75.0, 5.812724579633928)
    (100.0, 6.801233870153367)
    (150.0, 8.440391534639614)
};
\addplot[teal, mark=Mercedes star, mark size=2.7pt, thick] coordinates {
    (7.0, 1.1563409224122223)
    (10.0, 1.524972906373158)
    (15.0, 2.054772224650667)
    (20.0, 2.512080003618733)
    (35.0, 3.6248705483863635)
    (50.0, 4.51028551270734)
    (75.0, 5.712668549072003)
    (100.0, 6.7123029550358755)
    (150.0, 8.366851700525354)
};
\addplot[olive, mark=10-pointed star, mark size=2.4pt, thick] coordinates {
    (7.0, 1.018519813365892)
    (10.0, 1.363100136813772)
    (15.0, 1.8686781672608337)
    (20.0, 2.3129006240843117)
    (35.0, 3.4150763096597405)
    (50.0, 4.305740538150141)
    (75.0, 5.524544034374356)
    (100.0, 6.540635135296715)
    (150.0, 8.221283473599053)
};
\addplot[pink, mark=halfcircle*, mark size=2.1pt, thick] coordinates {
    (7.0, 0.9231417057713385)
    (10.0, 1.2460652333736042)
    (15.0, 1.7286531960615028)
    (20.0, 2.1583457818162035)
    (35.0, 3.241288999334673)
    (50.0, 4.1286581123573605)
    (75.0, 5.353472641089845)
    (100.0, 6.379537740629631)
    (150.0, 8.079783854236458)
};
\addplot[gray, mark=halfsquare right*, mark size=1.8pt, thick] coordinates {
    (10.0, 1.1558645245818309)
    (15.0, 1.617757605712653)
    (20.0, 2.033326858299125)
    (35.0, 3.0941692526453015)
    (50.0, 3.973914899356443)
    (75.0, 5.198293138710734)
    (100.0, 6.229581733774812)
    (150.0, 7.9437485650738076)
};
\logLogSlopeTriangle{0.9}{0.15}{0.1}{1}{thick, black}
\end{axis}
\end{tikzpicture}    \end{subfigure}\hfill
    \begin{subfigure}[t]{0.48\textwidth}
        \centering
        \begin{tikzpicture}[baseline=(current axis.north)]
\begin{axis}[title={(b)},
    tick label style={font=\normalsize},
    grid=both, width=\linewidth, height=5cm,
    xmode=log, ymode=log,
    xtick={1,2,5,10,20,50,100}, xticklabel=\empty,
    ytick={2,5,10,20,50},
    log ticks with fixed point,
    legend to name=nupeakc3legend,
    legend columns=6,
    legend cell align=left,
]
\addplot[blue, mark=*, mark size=3pt, thick] coordinates {
    (20.0, 3.1202299655917076)
    (50.0, 4.973512854858545)
    (75.0, 6.102116789564862)
    (100.0, 7.052363357685022)
    (150.0, 8.644988636429133)
};
\addplot[red, mark=square*, mark size=2.7pt, thick] coordinates {
    (5.0, 1.7876426715540128)
    (15.0, 3.2671118440579017)
    (35.0, 5.066886268091055)
    (50.0, 6.076457003742436)
    (75.0, 7.46145620515291)
    (100.0, 8.626895642081807)
    (150.0, 10.579382217772968)
};
\addplot[black!60!green, mark=triangle*, mark size=2.4pt, thick] coordinates {
    (7.0, 2.4538425203022185)
    (10.0, 2.9999755489738487)
    (15.0, 3.7402627232065893)
    (20.0, 4.357409192666132)
    (35.0, 5.829742559263257)
    (50.0, 6.998988188579216)
    (75.0, 8.601501541757695)
    (100.0, 9.949162143284697)
    (150.0, 12.205975098353292)
};
\addplot[orange, mark=diamond*, mark size=2.1pt, thick] coordinates {
    (2.0, 1.2645482024079886)
    (5.0, 2.3565351736635836)
    (7.0, 2.9015531362655134)
    (10.0, 3.5811280390186555)
    (15.0, 4.501897196667544)
    (20.0, 5.267678746184039)
    (35.0, 7.087816479031267)
    (50.0, 8.528542244262223)
    (75.0, 10.499325842487178)
    (100.0, 12.1546628745509)
    (150.0, 14.924336763692395)
};
\addplot[violet, mark=pentagon*, mark size=1.8pt, thick] coordinates {
    (2.0, 1.436680371148381)
    (5.0, 2.8090521254716165)
    (7.0, 3.514753303517521)
    (10.0, 4.404260685932946)
    (15.0, 5.616479366954341)
    (20.0, 6.625263871907663)
    (35.0, 9.014724546174786)
    (50.0, 10.89693154252029)
    (75.0, 13.462423378250602)
    (100.0, 15.612047004808343)
    (150.0, 19.202495489803)
};
\addplot[cyan, mark=otimes*, mark size=1.5pt, thick] coordinates {
    (2.0, 1.5392694998928798)
    (5.0, 3.107884446472889)
    (7.0, 3.9339510805010356)
    (10.0, 4.986659440621014)
    (15.0, 6.433447418195683)
    (20.0, 7.642377767999825)
    (35.0, 10.507381579465065)
    (50.0, 12.758897658206811)
    (75.0, 15.819176475914052)
    (100.0, 18.37766077251833)
    (150.0, 22.643675257825954)
};
\addplot[brown, mark=star, mark size=1.2pt, thick] coordinates {
    (2.0, 1.6357489299186545)
    (5.0, 3.414303273019182)
    (7.0, 4.376842880842602)
    (10.0, 5.6208928723035445)
    (15.0, 7.3523039894280116)
    (20.0, 8.810678564712635)
    (35.0, 12.282538784792244)
    (50.0, 15.011027822992734)
    (75.0, 18.710358575462916)
    (100.0, 21.79483737740811)
    (150.0, 26.92563629168914)
};
\addplot[magenta, mark=oplus*, mark size=3pt, thick] coordinates {
    (2.0, 1.7282924844228986)
    (5.0, 3.737949038648159)
    (7.0, 4.860886959454838)
    (10.0, 6.338645758313274)
    (15.0, 8.432120159913763)
    (20.0, 10.218765570844655)
    (35.0, 14.518522514538724)
    (50.0, 17.914740894797987)
    (75.0, 22.51693079596354)
    (100.0, 26.344810373751404)
    (150.0, 32.69252643809543)
};
\addplot[teal, mark=Mercedes star, mark size=2.7pt, thick] coordinates {
    (2.0, 1.7826928066751586)
    (5.0, 3.945642495111879)
    (7.0, 5.181671579784621)
    (10.0, 6.830300610929157)
    (15.0, 9.198819311163227)
    (20.0, 11.243176774862585)
    (35.0, 16.21813254868313)
    (50.0, 20.176805567453215)
    (75.0, 25.552941940804875)
    (100.0, 30.022753715483933)
    (150.0, 37.42128311609802)
};
\addplot[olive, mark=10-pointed star, mark size=2.4pt, thick] coordinates {
    (2.0, 1.844389372623217)
    (5.0, 4.201422908587226)
    (7.0, 5.58964518176149)
    (10.0, 7.477089301805434)
    (15.0, 10.245873245547312)
    (20.0, 12.678327551446191)
    (35.0, 18.71368746896329)
    (50.0, 23.591009609342755)
    (75.0, 30.265483315028778)
    (100.0, 35.83004586019951)
    (150.0, 45.034324287059555)
};
\addplot[pink, mark=halfcircle*, mark size=2.1pt, thick] coordinates {
    (2.0, 1.878639046376772)
    (5.0, 4.354982315812034)
    (7.0, 5.842602389848479)
    (10.0, 7.892214379296353)
    (15.0, 10.944337254961486)
    (20.0, 13.661428436344863)
    (35.0, 20.50912845781807)
    (50.0, 26.12033266814041)
    (75.0, 33.86553751934716)
    (100.0, 40.35409966211025)
    (150.0, 51.10628707252677)
};
\addplot[gray, mark=halfsquare right*, mark size=1.8pt, thick] coordinates {
    (2.0, 1.9004826457709014)
    (5.0, 4.458156478745238)
    (7.0, 6.016477243387379)
    (10.0, 8.184774907729494)
    (15.0, 11.450995611795564)
    (20.0, 14.389094306458935)
    (35.0, 21.889140587882505)
    (50.0, 28.10891756969233)
    (75.0, 36.765392299492504)
    (100.0, 44.05684373203306)
    (150.0, 56.1766606346151)
};
\logLogSlopeTriangle{0.9}{0.15}{0.1}{1}{thick, black}
\legend{$c_2 = 1$, $c_2 = 1.5$, $c_2 = 2$, $c_2 = 3$, $c_2 = 5$, $c_2 = 7$, $c_2 = 10$, $c_2 = 15$, $c_2 = 20$, $c_2 = 30$, $c_2 = 40$, $c_2 = 50$}
\end{axis}
\end{tikzpicture}    \end{subfigure}
    \par\vspace{0.4em}
    \begin{subfigure}[t]{0.48\textwidth}
        \centering
        \begin{tikzpicture}[baseline=(current axis.north)]
\begin{axis}[title={(c)},
    tick label style={font=\normalsize},
    xlabel={$c_3$}, ylabel={$\ell / \bar{h}$},
    grid=both, width=\linewidth, height=5cm,
    xmode=log, ymode=log,
    xtick={1,2,5,10,20,50,100},
    ytick={5,7,10,15,20},
    log ticks with fixed point,
]
\addplot[blue, mark=*, mark size=3pt, thick] coordinates {
    (20.0, 7.494969088546916)
    (50.0, 11.887990146873143)
    (75.0, 14.57830700320956)
    (100.0, 16.84577120021239)
    (150.0, 20.648070888693486)
};
\addplot[red, mark=square*, mark size=2.7pt, thick] coordinates {
    (20.0, 7.492947445537763)
    (50.0, 11.870160859832966)
    (75.0, 14.559125968634836)
    (100.0, 16.82721673116812)
    (150.0, 20.631534289358992)
};
\addplot[black!60!green, mark=triangle*, mark size=2.4pt, thick] coordinates {
    (20.0, 7.500110077168362)
    (35.0, 9.913812304503807)
    (50.0, 11.85816042226253)
    (75.0, 14.543710840863971)
    (100.0, 16.811050396498125)
    (150.0, 20.616044411987694)
};
\addplot[orange, mark=diamond*, mark size=2.1pt, thick] coordinates {
    (5.0, 4.021567551166154)
    (15.0, 6.559309900111987)
    (20.0, 7.530189543000934)
    (35.0, 9.916120780683238)
    (50.0, 11.847127027510457)
    (75.0, 14.522593843284648)
    (100.0, 16.78590320953534)
    (150.0, 20.588969947742395)
};
\addplot[violet, mark=pentagon*, mark size=1.8pt, thick] coordinates {
    (5.0, 4.259158154266508)
    (7.0, 4.839356465883794)
    (15.0, 6.6795519971317106)
    (20.0, 7.6213060323695485)
    (35.0, 9.955209710350582)
    (50.0, 11.857772538417358)
    (75.0, 14.508187209209197)
    (100.0, 16.758585764989295)
    (150.0, 20.55017650098233)
};
\addplot[cyan, mark=otimes*, mark size=1.5pt, thick] coordinates {
    (5.0, 4.489825599484887)
    (7.0, 5.044340436560297)
    (15.0, 6.815293263816177)
    (20.0, 7.731901967803252)
    (35.0, 10.017850131915488)
    (50.0, 11.89302188269424)
    (75.0, 14.517293674650986)
    (100.0, 16.752527843545185)
    (150.0, 20.528040949065346)
};
\addplot[brown, mark=star, mark size=1.2pt, thick] coordinates {
    (5.0, 4.819742274512861)
    (7.0, 5.329609340788599)
    (10.0, 6.018927091290619)
    (15.0, 7.024039528374355)
    (20.0, 7.9082553600518715)
    (35.0, 10.132247170150809)
    (50.0, 11.970827752832507)
    (75.0, 14.556779090349242)
    (100.0, 16.76848532649843)
    (150.0, 20.51695877170457)
};
\addplot[magenta, mark=oplus*, mark size=3pt, thick] coordinates {
    (5.0, 5.322596377796446)
    (7.0, 5.784019318486651)
    (10.0, 6.4204376348757695)
    (15.0, 7.3667949264483195)
    (20.0, 8.20931742905164)
    (35.0, 10.34691070093259)
    (50.0, 12.132270182254759)
    (75.0, 14.66051436873485)
    (100.0, 16.834323085788625)
    (150.0, 20.536599055288256)
};
\addplot[teal, mark=Mercedes star, mark size=2.7pt, thick] coordinates {
    (7.0, 6.205375337249956)
    (10.0, 6.802772425321135)
    (15.0, 7.699674792162109)
    (20.0, 8.505532128219295)
    (35.0, 10.57278588240701)
    (50.0, 12.312640677118223)
    (75.0, 14.788978887445975)
    (100.0, 16.927846078660973)
    (150.0, 20.585560020047854)
};
\addplot[olive, mark=10-pointed star, mark size=2.4pt, thick] coordinates {
    (7.0, 6.974911273589687)
    (10.0, 7.510128910113822)
    (15.0, 8.326805641463999)
    (20.0, 9.07472268702926)
    (35.0, 11.02602754481105)
    (50.0, 12.690043236521662)
    (75.0, 15.08021325319856)
    (100.0, 17.15637312701666)
    (150.0, 20.731802900888457)
};
\addplot[pink, mark=halfcircle*, mark size=2.1pt, thick] coordinates {
    (7.0, 7.680031017768221)
    (10.0, 8.157451105215433)
    (15.0, 8.91413588078179)
    (20.0, 9.614186617401538)
    (35.0, 11.468299175522379)
    (50.0, 13.069421665415643)
    (75.0, 15.385748949404103)
    (100.0, 17.41087849906305)
    (150.0, 20.913458008271736)
};
\addplot[gray, mark=halfsquare right*, mark size=1.8pt, thick] coordinates {
    (10.0, 8.756676752560075)
    (15.0, 9.464681365549305)
    (20.0, 10.127595012707092)
    (35.0, 11.896583199219137)
    (50.0, 13.443851186840588)
    (75.0, 15.696053238315843)
    (100.0, 17.6758890472554)
    (150.0, 21.113472686629322)
};
\end{axis}
\end{tikzpicture}    \end{subfigure}\hfill
    \begin{subfigure}[t]{0.48\textwidth}
        \centering
        \begin{tikzpicture}[baseline=(current axis.north)]
\begin{axis}[title={(d)},
    tick label style={font=\normalsize},
    xlabel={$c_3$},
    grid=both, width=\linewidth, height=5cm,
    xmode=log, ymode=log,
    xtick={1,2,5,10,20,50,100},
    ytick={5,10,20,50,100},
    log ticks with fixed point,
]
\addplot[blue, mark=*, mark size=3pt, thick] coordinates {
    (20.0, 7.494969088546916)
    (50.0, 11.887990146873143)
    (75.0, 14.57830700320956)
    (100.0, 16.84577120021239)
    (150.0, 20.648070888693486)
};
\addplot[red, mark=square*, mark size=2.7pt, thick] coordinates {
    (5.0, 4.704258883599129)
    (15.0, 7.95548852340104)
    (35.0, 12.14828057207144)
    (50.0, 14.537303525981072)
    (75.0, 17.83084846762904)
    (100.0, 20.608773042853258)
    (150.0, 25.26822902620769)
};
\addplot[black!60!green, mark=triangle*, mark size=2.4pt, thick] coordinates {
    (7.0, 6.408652820174983)
    (10.0, 7.5723872001250125)
    (15.0, 9.208335354362614)
    (20.0, 10.606377455971435)
    (35.0, 14.018019688231336)
    (50.0, 16.768966105641606)
    (75.0, 20.567178293912765)
    (100.0, 23.77399496588343)
    (150.0, 29.155109985992066)
};
\addplot[orange, mark=diamond*, mark size=2.1pt, thick] coordinates {
    (2.0, 4.985881801738666)
    (5.0, 6.953944398852188)
    (7.0, 8.016338928420328)
    (10.0, 9.399465876773332)
    (15.0, 11.356790527075196)
    (20.0, 13.040490648261605)
    (35.0, 17.171559059117257)
    (50.0, 20.51823697541074)
    (75.0, 25.15262273847213)
    (100.0, 29.073149501767265)
    (150.0, 35.66051550605722)
};
\addplot[violet, mark=pentagon*, mark size=1.8pt, thick] coordinates {
    (2.0, 7.163952118332397)
    (5.0, 9.507753884296982)
    (7.0, 10.800234342596125)
    (10.0, 12.500156180547823)
    (15.0, 14.93171467770191)
    (20.0, 17.037521951214085)
    (35.0, 22.25386123752075)
    (50.0, 26.51222821856266)
    (75.0, 32.439457402008024)
    (100.0, 37.47166647206958)
    (150.0, 45.950474672145226)
};
\addplot[cyan, mark=otimes*, mark size=1.5pt, thick] coordinates {
    (2.0, 9.272153610930065)
    (5.0, 11.854816565434756)
    (7.0, 13.30837041960205)
    (10.0, 15.238814019090272)
    (15.0, 18.021273374732345)
    (20.0, 20.447546783533777)
    (35.0, 26.49532512147088)
    (50.0, 31.462758280398905)
    (75.0, 38.40675263186346)
    (100.0, 44.320703227620896)
    (150.0, 54.31054329400744)
};
\addplot[brown, mark=star, mark size=1.2pt, thick] coordinates {
    (2.0, 12.374183278514387)
    (5.0, 15.195583297799564)
    (7.0, 16.82227081775046)
    (10.0, 19.009987100289866)
    (15.0, 22.19640504946554)
    (20.0, 24.99734459472096)
    (35.0, 32.030060425126784)
    (50.0, 37.84970508387815)
    (75.0, 46.02898612661419)
    (100.0, 53.02343901501503)
    (150.0, 64.87823506779264)
};
\addplot[magenta, mark=oplus*, mark size=3pt, thick] coordinates {
    (2.0, 17.473900313888358)
    (5.0, 20.542794141732053)
    (7.0, 22.357637083645603)
    (10.0, 24.84088021481707)
    (15.0, 28.513454494991375)
    (20.0, 31.77698769352638)
    (35.0, 40.06313217165813)
    (50.0, 46.98185260904267)
    (75.0, 56.774296195939044)
    (100.0, 65.1948661772602)
    (150.0, 79.53489650686647)
};
\addplot[teal, mark=Mercedes star, mark size=2.7pt, thick] coordinates {
    (2.0, 22.533162914686873)
    (5.0, 25.76345498831848)
    (7.0, 27.702341812293707)
    (10.0, 30.385185460325403)
    (15.0, 34.40719327013601)
    (20.0, 38.01703571794891)
    (35.0, 47.26921932047778)
    (50.0, 55.05429218265656)
    (75.0, 66.13108782661281)
    (100.0, 75.69867684336093)
    (150.0, 92.05733561353621)
};
\addplot[olive, mark=10-pointed star, mark size=2.4pt, thick] coordinates {
    (2.0, 32.601186199343914)
    (5.0, 36.036866962146846)
    (7.0, 38.135826574366284)
    (10.0, 41.08030734158949)
    (15.0, 45.575747825808804)
    (20.0, 49.67761302328728)
    (35.0, 60.37343954541979)
    (50.0, 69.49221319818867)
    (75.0, 82.58044315403936)
    (100.0, 93.96274142140753)
    (150.0, 113.54782471961559)
};
\addplot[pink, mark=halfcircle*, mark size=2.1pt, thick] coordinates {
    (2.0, 42.63969523731348)
    (5.0, 46.20338515308426)
    (7.0, 48.40679278628011)
    (10.0, 51.52193259109819)
    (15.0, 56.32968054063711)
    (20.0, 60.76738587542227)
    (35.0, 72.50965552124346)
    (50.0, 82.64249882101524)
    (75.0, 97.29485452760704)
    (100.0, 110.10714965329535)
    (150.0, 132.2616308921598)
};
\addplot[gray, mark=halfsquare right*, mark size=1.8pt, thick] coordinates {
    (2.0, 52.66463239622592)
    (5.0, 56.31714752767015)
    (7.0, 58.59514808084691)
    (10.0, 61.83372520822843)
    (15.0, 66.8687663221697)
    (20.0, 71.55354335067764)
    (35.0, 84.09721506295217)
    (50.0, 95.0406100179735)
    (75.0, 110.97471269827687)
    (100.0, 124.97574560787676)
    (150.0, 149.2859359978872)
};
\end{axis}
\end{tikzpicture}    \end{subfigure}
    \figurelegend{nupeakc3legend}
    \caption{\case{burgers-c4}: (a,b)~peak $\nu_\mathrm{peak}/\bar h$ and (c,d)~half-width $\ell/\bar h$ versus $c_3$, one curve per $c_2$, for the (a,c)~standard and (b,d)~derived parametrizations. Triangles: slopes given by~\eqref{eq:greens_std} and~\eqref{eq:greens_reparam}; none is drawn in~(c,d), where both equations give a width independent of $c_3$.}
    \label{fig:c4_test_burgers_c3}
\end{figure}

Over the full grid, with peak and half-width normalized by $\bar h$ so that they describe the AV's intrinsic shape, the standard parametrization shows the peak decreasing and the half-width increasing with $c_2$, with the signs the Green's-function $\nu_\mathrm{peak} \sim 1/\sqrt{c_2}$ and $\ell \sim \sqrt{c_2}$ predictions require but with slopes well short of $\pm 1/2$ (Figs.~\ref{fig:c4_test_burgers_c2} and~\ref{fig:c4_test_burgers_c3}, Table~\ref{tab:sensitivity_burgers}). The derived parametrization inverts the peak's trend, $\nu_\mathrm{peak}/\bar h$ growing with $c_2$, while $\ell/\bar h$ grows roughly linearly with $c_2$ as the Green's-function prediction $\ell = c_2$ requires. The derived width also varies more strongly with $c_2$ than the standard width does. The peak retains a $c_3$-modulated sensitivity to $c_2$, so the width--magnitude decoupling of the Green's function is not recovered exactly, but the residual trend has the preferred sign: widening the AV reinforces rather than reduces its peak.

To quantify these trends we extract four log--log sensitivity exponents at every point of the $(c_2, c_3)$ grid,
\begin{equation}
    s^{\nu}_{c_3} = \frac{\partial \log \nu_\mathrm{peak}}{\partial \log c_3}, \qquad
    s^{\nu}_{c_2} = \frac{\partial \log \nu_\mathrm{peak}}{\partial \log c_2}, \qquad
    s^{\ell}_{c_2} = \frac{\partial \log \ell}{\partial \log c_2}, \qquad
    s^{\ell}_{c_3} = \frac{\partial \log \ell}{\partial \log c_3},
    \label{eq:sensitivity_exponents}
\end{equation}
the first two for the peak and the last two for the half-width. Each is computed by centered finite differences along its axis (one-sided at the grid edges), with the other parameter held fixed, and the tables report the minimum and maximum over the grid. The Green's function predicts $s^{\nu}_{c_3} = s^{\ell}_{c_2} = 1$ and $s^{\nu}_{c_2} = s^{\ell}_{c_3} = 0$ for the derived parametrization~\eqref{eq:greens_reparam}, and $s^{\nu}_{c_3} = 1$, $s^{\nu}_{c_2} = -1/2$, $s^{\ell}_{c_2} = 1/2$, $s^{\ell}_{c_3} = 0$ for the standard one~\eqref{eq:greens_std}. Ideally both cross exponents vanish, $s^{\nu}_{c_2} = s^{\ell}_{c_3} = 0$, so that $c_2$ controls the width only and $c_3$ the magnitude only, and where they cannot, $s^{\nu}_{c_2}$ should be positive rather than negative (Section~\ref{subsec:greens}). In the coupled problem, the cross exponent $s^{\nu}_{c_2}$ changes sign between the parametrizations, from negative under the standard one to positive under the derived one (Table~\ref{tab:sensitivity_burgers}). The width responds to $c_2$ under the derived parametrization ($s^{\ell}_{c_2} \in [0.50, 0.95]$) but barely under the standard one ($[0.00, 0.32]$), well short of its Green's-function value $1/2$.

\begin{table}[hbt!]
\centering
\caption{\case{burgers-c4}: the four sensitivity exponents from the Green's function,~\eqref{eq:greens_std} and~\eqref{eq:greens_reparam}, and their min--max range in the coupled state--AV problem over the $(c_2, c_3)$ sweep, for the derived ($c_4 = 1/c_2$) and standard ($c_4 = 1$) parametrizations.}
\label{tab:sensitivity_burgers}
\begin{tabular}{@{} l l cc cc @{}}
\toprule
 & & \multicolumn{2}{c}{$c_4 = 1/c_2$} & \multicolumn{2}{c}{$c_4 = 1$} \\
\cmidrule(lr){3-4}\cmidrule(lr){5-6}
Exponent & Role & Green's function & coupled & Green's function & coupled \\
\midrule
$s^{\nu}_{c_3}$ & $c_3 \to$ magnitude (primary) & 1 & $[0.50, 0.93]$ & 1 & $[0.50, 0.82]$ \\
$s^{\nu}_{c_2}$ & $c_2 \to$ magnitude (cross)   & 0 & $[0.05, 0.50]$ & $-1/2$ & $[-0.34, 0.00]$ \\
$s^{\ell}_{c_2}$ & $c_2 \to$ width (primary)    & 1 & $[0.50, 0.95]$ & $1/2$ & $[0.00, 0.32]$ \\
$s^{\ell}_{c_3}$ & $c_3 \to$ width (cross)      & 0 & $[0.07, 0.50]$ & 0 & $[0.20, 0.50]$ \\
\bottomrule
\end{tabular}
\end{table}

\subsubsection{Mach~5 inviscid flow over a half-cylinder}
\label{subsubsec:c4_test_cyl}

The purpose of this experiment (test case: \case{cyl-euler-c4}) is to test whether the inversion confirmed on Burgers' equation carries over to a multidimensional problem with a curved shock. Consider the Mach~5 inviscid flow over a half-cylinder of Section~\ref{subsubsec:selfsim_cyl} on the coarsest mesh, $n_{\mathrm{r}} = n_{\mathrm{a}} = 14$ and $p = 4$. We sweep a logarithmically spaced $(c_2, c_3)$ grid for both parametrizations and extract $\nu_\mathrm{peak}$ and $\ell$ from the AV profile along the stagnation line, with the definitions of Section~\ref{subsubsec:c4_test_burgers}.

The Burgers' pattern is reproduced (Table~\ref{tab:sensitivity_euler}). The sensitivity exponents show the same two features: the width responds to $c_2$ under the derived parametrization but hardly at all under the standard one, and $s^{\nu}_{c_2}$ changes sign from negative to positive. Both parametrizations show the peak growing with $c_3$, consistent with $\nu_\mathrm{peak} \sim c_3$, and a residual $c_3$-dependence of the half-width. The cross-coupling that survives in the derived parametrization is therefore robust to the increased dimensionality and the curved-shock geometry.

\begin{table}[hbt!]
\centering
\caption{\case{cyl-euler-c4}: the four sensitivity exponents from the Green's function and their min--max range in the coupled state--AV problem over the $(c_2, c_3)$ sweep, in the form of Table~\ref{tab:sensitivity_burgers}.}
\label{tab:sensitivity_euler}
\begin{tabular}{@{} l l cc cc @{}}
\toprule
 & & \multicolumn{2}{c}{$c_4 = 1/c_2$} & \multicolumn{2}{c}{$c_4 = 1$} \\
\cmidrule(lr){3-4}\cmidrule(lr){5-6}
Exponent & Role & Green's function & coupled & Green's function & coupled \\
\midrule
$s^{\nu}_{c_3}$ & $c_3 \to$ magnitude (primary) & 1 & $[0.49, 0.76]$ & 1 & $[0.48, 1.00]$ \\
$s^{\nu}_{c_2}$ & $c_2 \to$ magnitude (cross)   & 0 & $[0.25, 0.48]$ & $-1/2$ & $[-0.44, 0.00]$ \\
$s^{\ell}_{c_2}$ & $c_2 \to$ width (primary)    & 1 & $[0.16, 0.90]$ & $1/2$ & $[-0.01, 0.41]$ \\
$s^{\ell}_{c_3}$ & $c_3 \to$ width (cross)      & 0 & $[0.01, 0.48]$ & 0 & $[0.11, 0.48]$ \\
\bottomrule
\end{tabular}
\end{table}

\subsubsection{Mach~5 viscous flow over a half-cylinder}
\label{subsec:c4_practice}

The purpose of this experiment (test case: \case{cyl-ns-widen}) is to show the consequence of the parametrization for the state on a viscous flow: whether widening the AV to damp oscillations near the captured shock matters. Consider the half-cylinder domain of Section~\ref{subsubsec:selfsim_cyl} with the compressible Navier--Stokes equations of Section~\ref{subsubsec:nseuler} at Mach~5 and $\mathrm{Re} = 1000$ based on the cylinder radius~\cite{dong2025rpadapt}. The viscosity is constant and equal to its freestream value $\mu_\infty$. The cylinder is a no-slip isothermal wall at $T_{\mathrm{w}} = 2.5\,T_\infty$, the outer arc carries the freestream state, and the two straight edges are outflow boundaries that impose the back pressure $P_\mathrm{out} = 2.5\,P_\infty$ where the outflow is subsonic, extrapolate the state where it is supersonic, and set the viscous flux to zero. The AV working variable is pinned to $\hat\nu = 0$ at the wall and takes the homogeneous Neumann condition on the other boundaries. The AV filter uses $a = 10^{-6}\mu_\infty$ and $\nu_0 = 0$, and the AV PDE uses $c_1 = 0.1$.

The experiment uses the $n_{\mathrm{r}} = n_{\mathrm{a}} = 14$, $p = 4$ mesh of Fig.~\ref{fig:halfcyl_mesh}. The experiment has three steps. First, choose a standard baseline $(c_2^{\,s}, c_3^{\,s})$ that leaves visible oscillations. Second, construct the derived baseline that produces the same Green's-function AV, by equating the peak and width of Section~\ref{subsec:greens},
\begin{equation}
    c_2^{\,d} = \sqrt{c_2^{\,s}}, \qquad c_3^{\,d} = c_3^{\,s}/\sqrt{c_2^{\,s}}.
    \label{eq:std_der_match}
\end{equation}
Third, widen both by the same factor $W$ applied to their respective $c_2$, with $c_3$ held fixed, and record the change in the peak AV $\nu_\mathrm{peak}$ and in the peak Mach overshoot above the freestream, $\max(M - M_\infty)$, a proxy for the oscillation amplitude. The standard baseline is $(c_2^{\,s}, c_3^{\,s}) = (2, 0.2)$, for which the AV is modest enough that oscillations remain near the captured shock, and the widening factors are $W \in \{\sqrt{2}, 2, 2\sqrt{2}\}$.

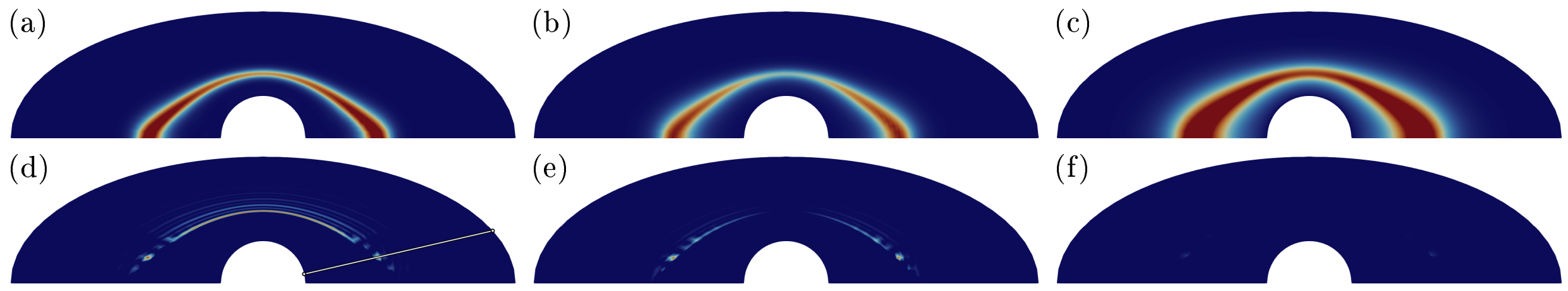
\begin{figure}[hbt!]
    \centering
    \begin{tikzpicture}
        \tikzset{cell/.style={inner sep=2pt, outer sep=0pt, anchor=north west}}
        \node[cell] (A) at (0,0) {\includegraphics[width=0.32\linewidth]{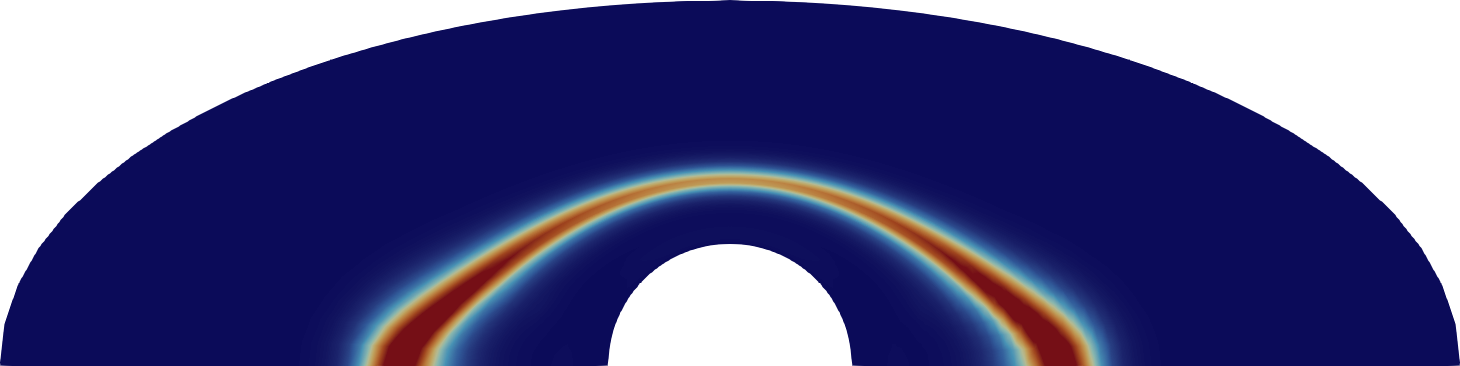}};
        \node[cell] (B) at ([xshift=0.06cm]A.north east) {\includegraphics[width=0.32\linewidth]{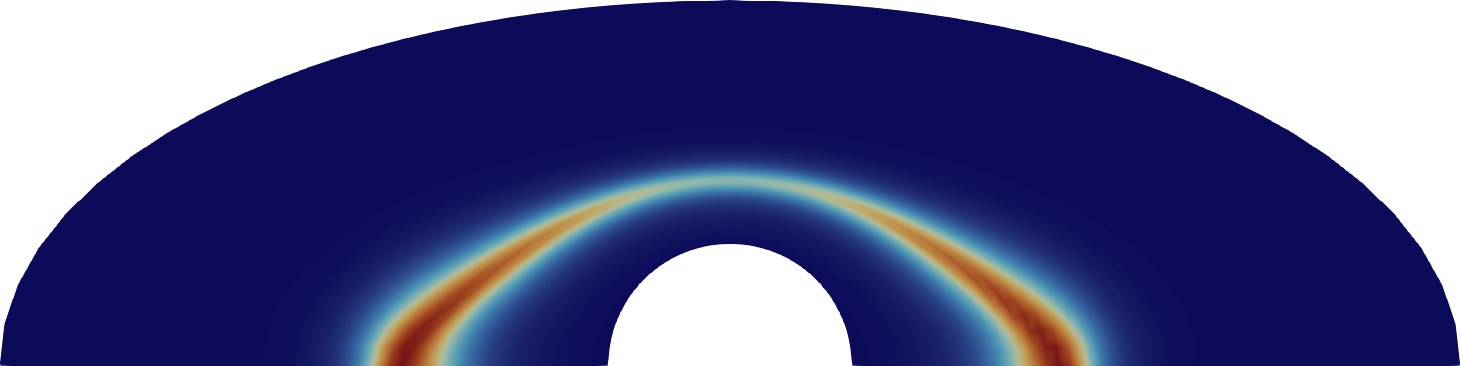}};
        \node[cell] (C) at ([xshift=0.06cm]B.north east) {\includegraphics[width=0.32\linewidth]{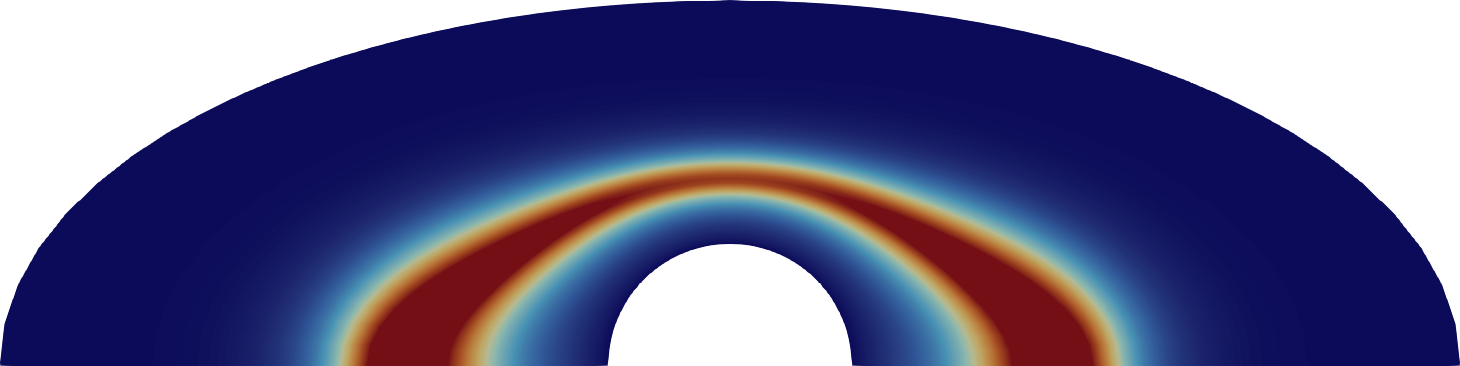}};
        \node[cell] (D) at ([yshift=-0.06cm]A.south west) {\includegraphics[width=0.32\linewidth]{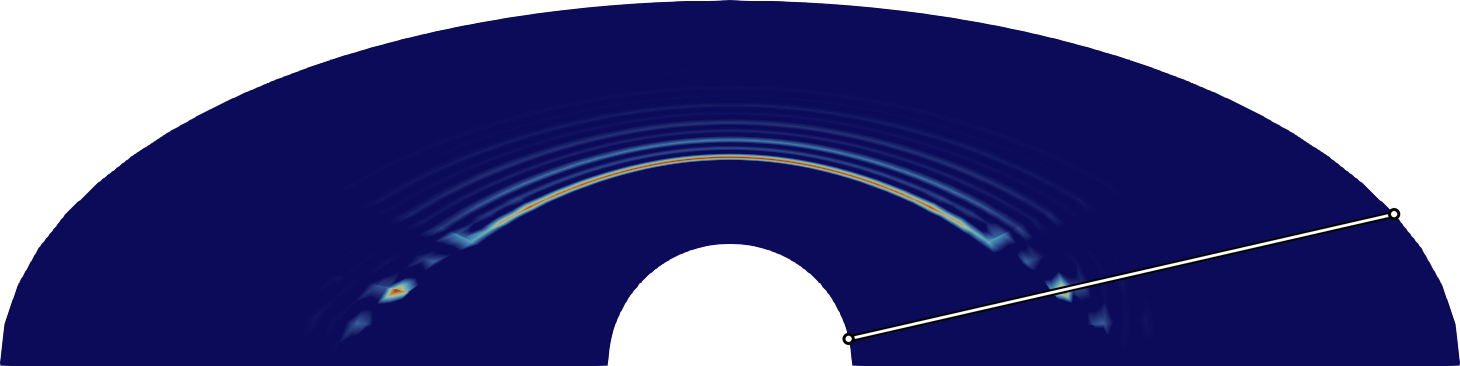}};
        \node[cell] (E) at ([xshift=0.06cm]D.north east) {\includegraphics[width=0.32\linewidth]{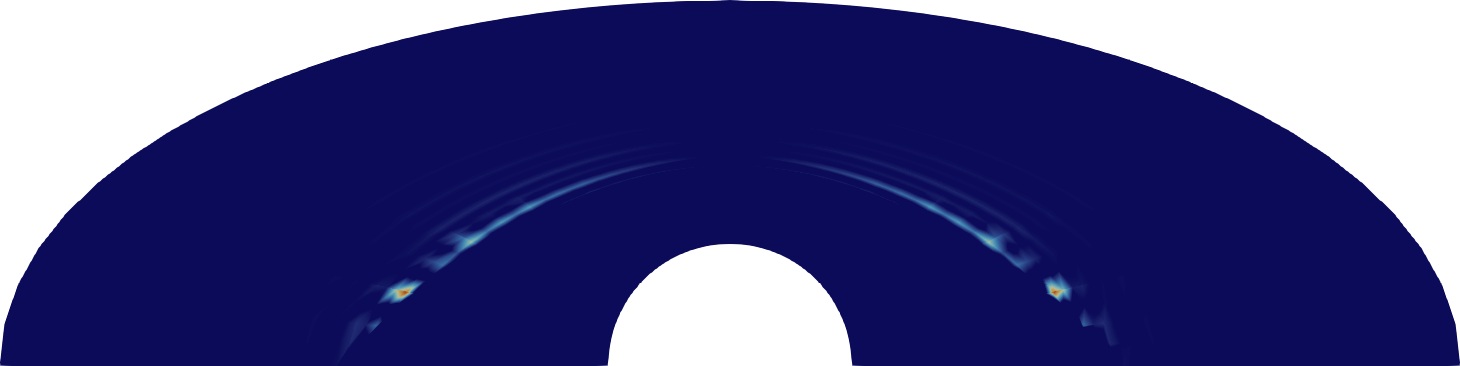}};
        \node[cell] (F) at ([xshift=0.06cm]E.north east) {\includegraphics[width=0.32\linewidth]{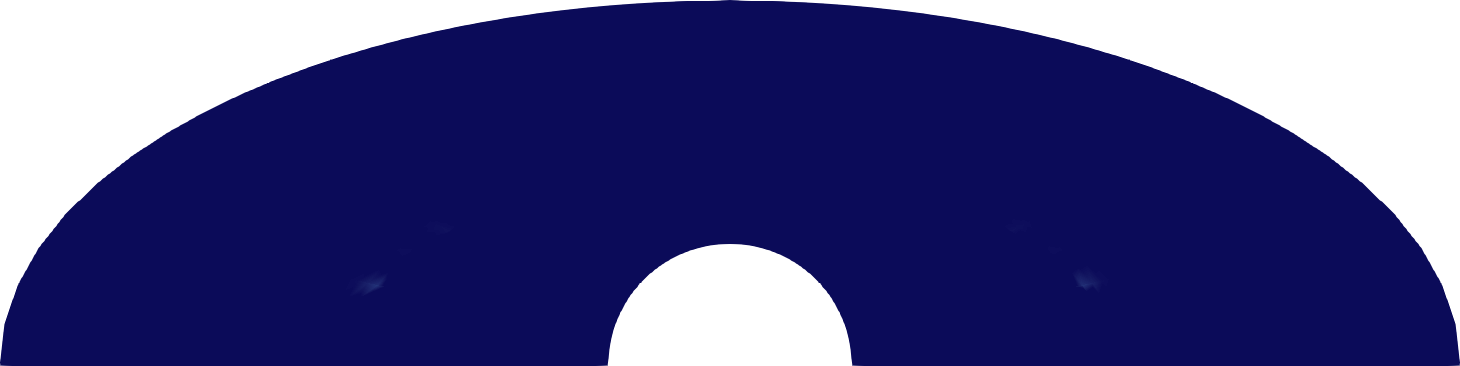}};
        \node[plabel] at (A.north west) {(a)};
        \node[plabel] at (B.north west) {(b)};
        \node[plabel] at (C.north west) {(c)};
        \node[plabel] at (D.north west) {(d)};
        \node[plabel] at (E.north west) {(e)};
        \node[plabel] at (F.north west) {(f)};
    \end{tikzpicture}
    \caption{\case{cyl-ns-widen}, the baseline and its $W = 2\sqrt{2}$ widening: (a--c)~AV $\nu$, on a common color scale from $0$ to $0.22$, and (d--f)~Mach overshoot $\max(M - M_\infty, 0)$, not the Mach number itself, on a color scale from $0$ to $0.085$ to expose the oscillations, for (a,d)~the baseline, (b,e)~the standard parametrization, and (c,f)~the derived parametrization. The white line in~(d) is the off-axis sampling line $\theta = \pi/14$ of Fig.~\ref{fig:367s_offaxis} and Table~\ref{tab:367s_widening}.}
    \label{fig:367s_2d_grid}
\end{figure}

\begin{figure}[hbt!]
    \centering
    \begin{subfigure}[t]{0.48\textwidth}
        \centering
        \begin{tikzpicture}
        \begin{axis}[
            width=\linewidth, height=5.0cm,
            title={(a)},
            xlabel={$r$}, ylabel={$\nu$},
            xmin=1.0, xmax=3.0,
            ymajorgrids=true,
            extra x ticks={1.328, 1.657, 1.985, 2.314, 2.642, 2.971},
            extra x tick style={
                grid=major,
                grid style={black!55, densely dashed, semithick},
                major tick length=0pt,
                tick label style={opacity=0},
            },
            legend to name=offaxis_legend367s,
            legend columns=-1,
            legend style={draw=black, fill=white, /tikz/every even column/.append style={column sep=0.6cm}},
            legend cell align=left,
        ]
        \addplot[thick, black] table[x=r, y=AV] {_tikz/367s_offaxis_baseline.dat};
        \addlegendentry{baseline}
        \addplot[thick, red] table[x=r, y=AV] {_tikz/367s_offaxis_std_wide.dat};
        \addlegendentry{standard parametrization}
        \addplot[thick, blue] table[x=r, y=AV] {_tikz/367s_offaxis_der_wide.dat};
        \addlegendentry{derived parametrization}
        \end{axis}
        \end{tikzpicture}    \end{subfigure}\hfill
    \begin{subfigure}[t]{0.48\textwidth}
        \centering
        \begin{tikzpicture}
        \begin{axis}[
            width=\linewidth, height=5.0cm,
            title={(b)},
            xlabel={$r$}, ylabel={$M$},
            xmin=1.0, xmax=3.0,
            ymajorgrids=true,
            extra x ticks={1.328, 1.657, 1.985, 2.314, 2.642, 2.971},
            extra x tick style={
                grid=major,
                grid style={black!55, densely dashed, semithick},
                major tick length=0pt,
                tick label style={opacity=0},
            },
        ]
        \addplot[thick, black] table[x=r, y=Mach] {_tikz/367s_offaxis_baseline.dat};
        \addplot[thick, red] table[x=r, y=Mach] {_tikz/367s_offaxis_std_wide.dat};
        \addplot[thick, blue] table[x=r, y=Mach] {_tikz/367s_offaxis_der_wide.dat};
        \end{axis}
        \end{tikzpicture}    \end{subfigure}
    \figurelegend{offaxis_legend367s}
    \caption{\case{cyl-ns-widen}: (a)~AV $\nu$ and (b)~Mach number along the off-axis line $\theta = \pi/14$ of Fig.~\ref{fig:367s_2d_grid}(d) for the baseline (black), the standard parametrization (red), and the derived parametrization (blue). Dashed vertical gridlines mark the radial element boundaries.}
    \label{fig:367s_offaxis}
\end{figure}
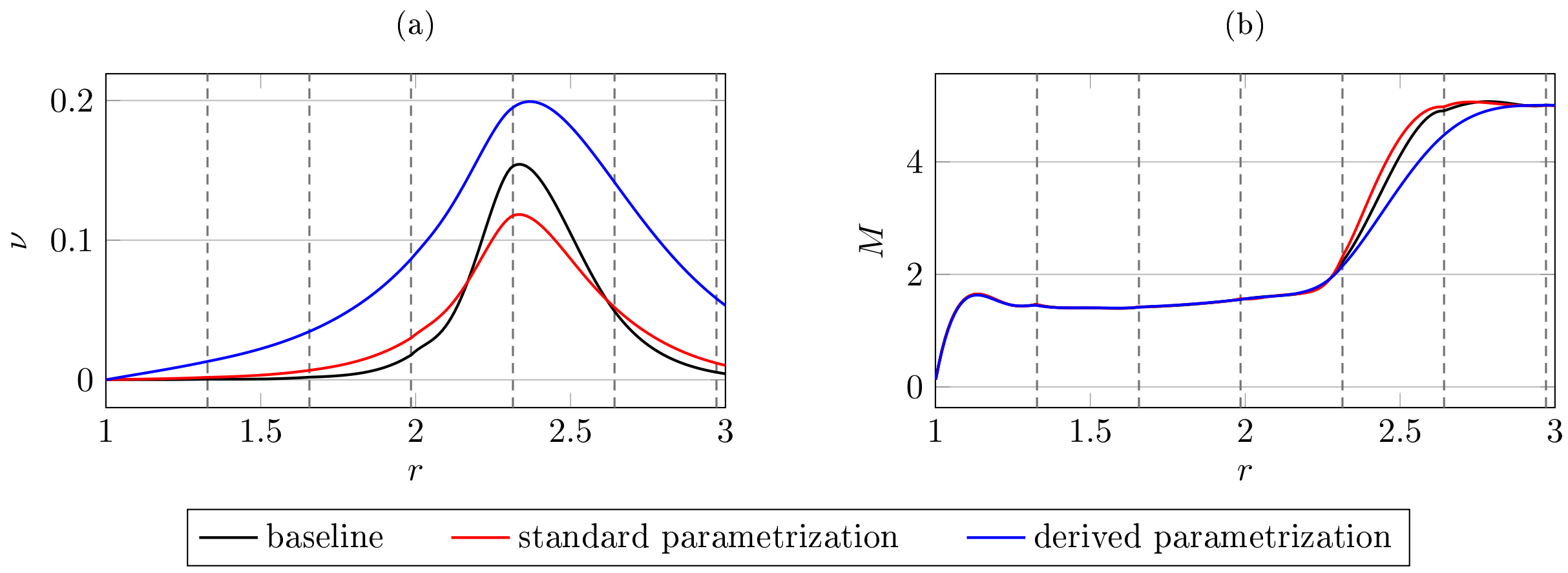

At $W = 2\sqrt{2}$, widening under the standard parametrization thickens the AV band but lowers its peak and leaves the oscillation amplitude comparable to the baseline, while widening under the derived parametrization thickens the band, raises the peak, and largely eliminates the oscillations (Fig.~\ref{fig:367s_2d_grid}). Over the whole domain, the baseline has $\nu_\mathrm{peak} = 0.177$ and overshoot $+0.072$; these change by $-25\%$ and $-10\%$ under the standard parametrization and by $+22\%$ and $-91\%$ under the derived one. Along the off-axis line $\theta = \pi/14$ drawn in Fig.~\ref{fig:367s_2d_grid}(d) through the point of maximum overshoot, where $\theta$ is the polar angle of the mesh with $\theta = \pi/2$ the stagnation streamline, the baseline $\nu_\mathrm{peak}$ is $0.154$ and the trend is monotone in $W$: as $W$ grows the standard peak falls monotonically while the overshoot decreases by at most $10\%$, and the derived peak rises monotonically while the overshoot collapses, by $91\%$ at $W = 2\sqrt{2}$ (Fig.~\ref{fig:367s_offaxis}, Table~\ref{tab:367s_widening}). Together with the Burgers' and inviscid cylinder cases, this confirms on a viscous flow that the derived parametrization widens the AV without lowering its peak or leaving oscillations near the captured shock.

\begin{table}[hbt!]
    \centering
    \caption{\case{cyl-ns-widen}: change in $\nu_\mathrm{peak}$ and Mach overshoot along the off-axis line $\theta = \pi/14$ as functions of the widening factor $W$, for both parametrizations. Baseline values: $\nu_\mathrm{peak} = 0.154$, overshoot $+0.072$.}
    \label{tab:367s_widening}
    \begin{tabular}{l cc cc}
        \toprule
         & \multicolumn{2}{c}{Standard ($c_4 = 1$)} & \multicolumn{2}{c}{Derived ($c_4 = 1/c_2$)} \\
        \cmidrule(lr){2-3}\cmidrule(lr){4-5}
        $W$ & $\Delta\nu_\mathrm{peak}$ & $\Delta$overshoot & $\Delta\nu_\mathrm{peak}$ & $\Delta$overshoot \\
        \midrule
        $\sqrt{2}$  & $-7\%$  & $-4\%$  & $+10\%$ & $-46\%$ \\
        $2$         & $-15\%$ & $-9\%$  & $+20\%$ & $-73\%$ \\
        $2\sqrt{2}$ & $-23\%$ & $-10\%$ & $+29\%$ & $-91\%$ \\
        \bottomrule
    \end{tabular}
\end{table}

\section{Application to hypersonic viscous flows}
\label{sec:validation}
\label{sec:applications}

The two cases in this section exercise the complete method, with the scaling of Section~\ref{sec:avpde_design} and the parametrization of Section~\ref{sec:reparam}, on hypersonic viscous flows with shock interactions and thin boundary layers. For each we ask whether the captured shock remains self-similar under refinement, so that it occupies a fixed number of elements; for the sphere we also compare the wall quantities with reference data. Both cases are solved in entropy variables.

\subsection{Mach 11.68 viscous flow over a double wedge}
\label{subsec:val_dwedge}

\begin{figure}[hbt!]
    \centering
    \includegraphics[width=0.32\linewidth]{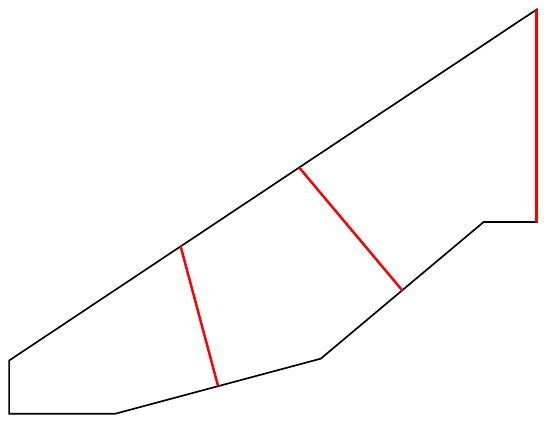}
    \caption{\case{dwedge}: domain and the three sampling lines (red), wall-normal through the midpoints of the two wedge faces and vertical along the outlet. Boundaries: no-slip isothermal walls on the wedges and the trailing surface, symmetry plane ahead of the first wedge, freestream on the inflow and upper boundaries, outflow on the right.}
    \label{fig:dw_meshes}
\end{figure}

Consider Mach~11.68 viscous flow over a double wedge~\cite{razavi2026output}, a $15^\circ$ wedge followed by a $40^\circ$ wedge, governed by the compressible Navier--Stokes equations of Section~\ref{subsubsec:nseuler} at $\mathrm{Re} = 10^4$ based on the length $L$ of each wedge face (Fig.~\ref{fig:dw_meshes}). The viscosity is constant, $\mu = \mu_\infty$. The wedge surfaces and the trailing surface are no-slip isothermal walls at $T_{\mathrm{w}} = 3.75\,T_\infty$, a symmetry plane precedes the first wedge, the inflow and upper boundaries carry the freestream state, and the right boundary is an outflow boundary as in Section~\ref{subsec:c4_practice}, with back pressure $P_\mathrm{out} = 4\,P_\infty$ imposed where the outflow is subsonic. The symmetry plane has length $0.5L$, the trailing surface $0.25L$, and the inflow boundary height $0.25L$, and the outlet rises $L$ above the second wedge. The AV working variable is pinned to $\hat\nu = 0$ on the walls and takes the homogeneous Neumann condition on the other boundaries. The AV filter uses $a = 10^{-4}\mu_\infty$ and $\nu_0 = 10^{-2}\mu_\infty$, and the AV PDE uses $c_1 = 0.01$.

The experiment (test case: \case{dwedge}) tests whether the converged AV remains self-similar under refinement in a hypersonic flow with interacting shocks and a shock--boundary-layer interaction, so that the captured shock occupies a fixed number of elements, and how the coefficients $(c_2, c_3)$ shape it. We use $p = 4$ on a coarse mesh of $1488$ elements and a uniformly refined mesh of $5952$ elements, with the baseline $(c_2, c_3) = (1.7, 1.7)$ unless stated otherwise. The coarse mesh is structured, with 62 elements along the wall and 24 normal to it, and the refined mesh splits each element into four. The wall-normal spacing is clustered toward the wall by the map $s \mapsto 1 + \tanh(s - 1)/\tanh 1$ of a uniform $s \in [0, 1]$. Because the diffusion tensor $\eta$ follows the element shape, the mesh map is made $C^1$ to avoid jumps in the AV at the corner lines: the streamwise spacing varies smoothly across the corners, and the mesh lines round each corner over a half-width that grows from zero at the wall, which stays sharp, to 0.12 at the upper boundary. Profiles are extracted along three lines: wall-normal lines through the midpoints of the first and second wedge faces (\emph{wedge~1}, \emph{wedge~2}) and a vertical line along the outlet (\emph{right outlet}).

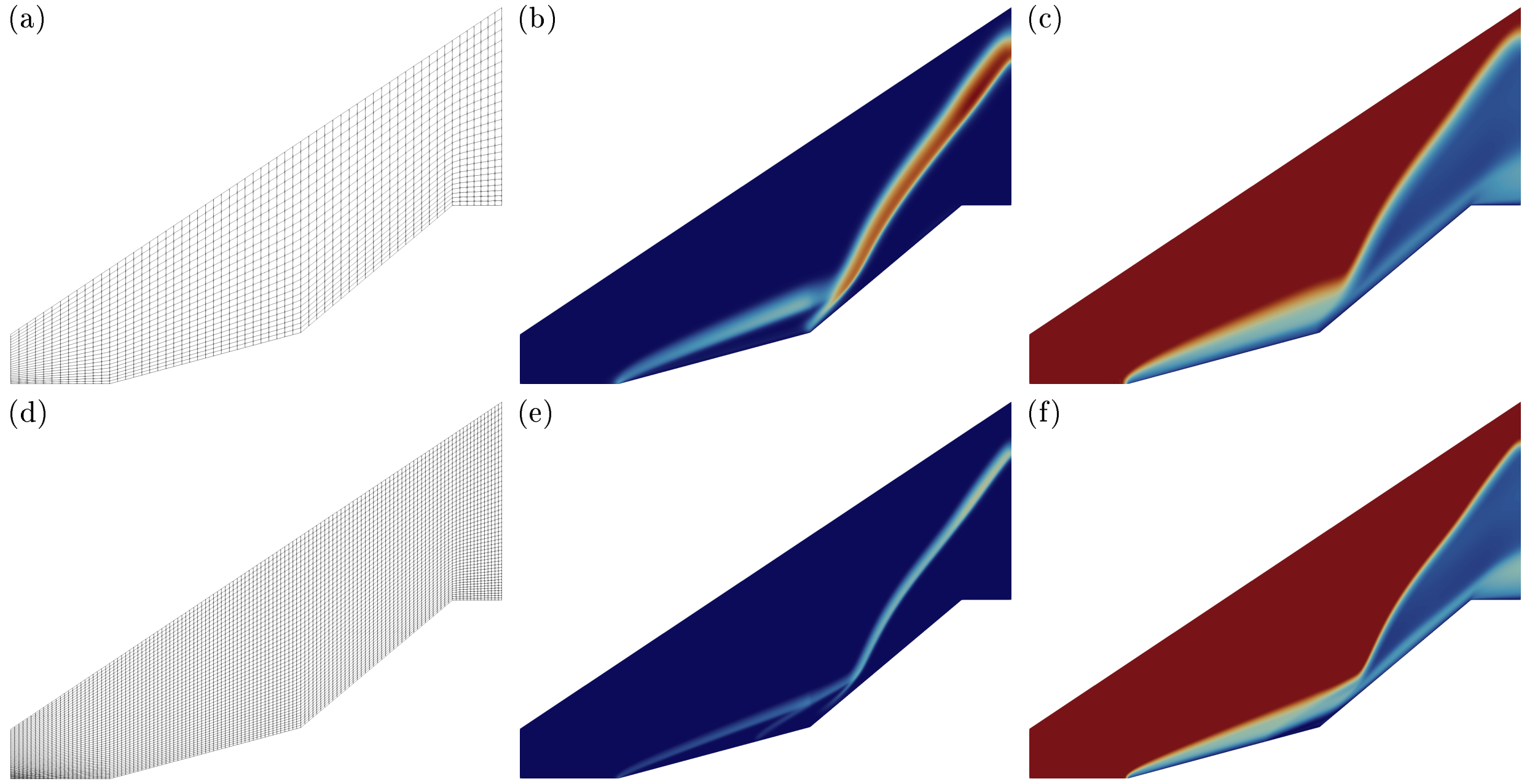
\begin{figure}[hbt!]
    \centering
    \begin{tikzpicture}
        \tikzset{cell/.style={inner sep=2pt, outer sep=0pt, anchor=north west}}
        \node[cell] (A) at (0,0) {\includegraphics[width=0.32\linewidth]{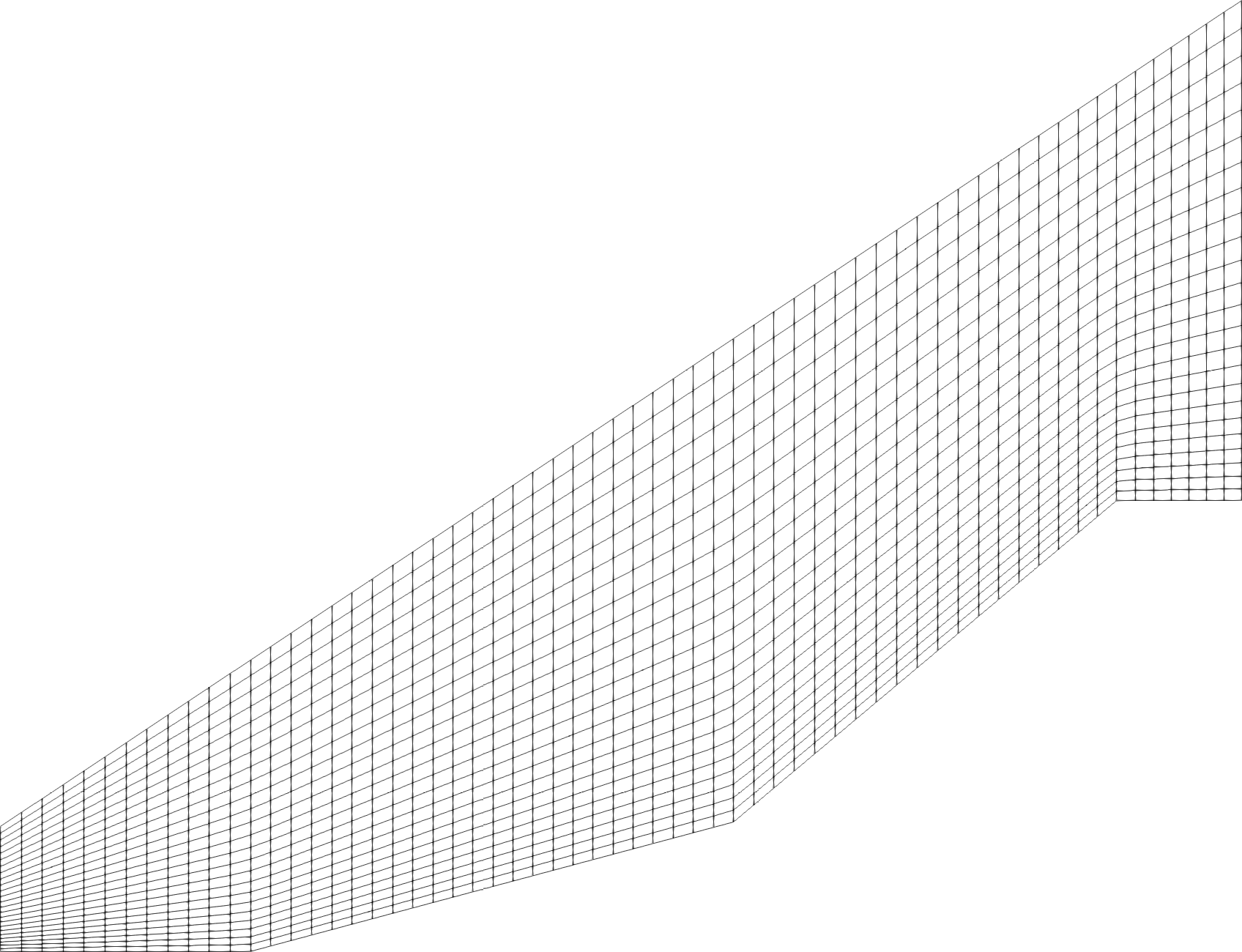}};
        \node[cell] (B) at ([xshift=0.06cm]A.north east) {\includegraphics[width=0.32\linewidth]{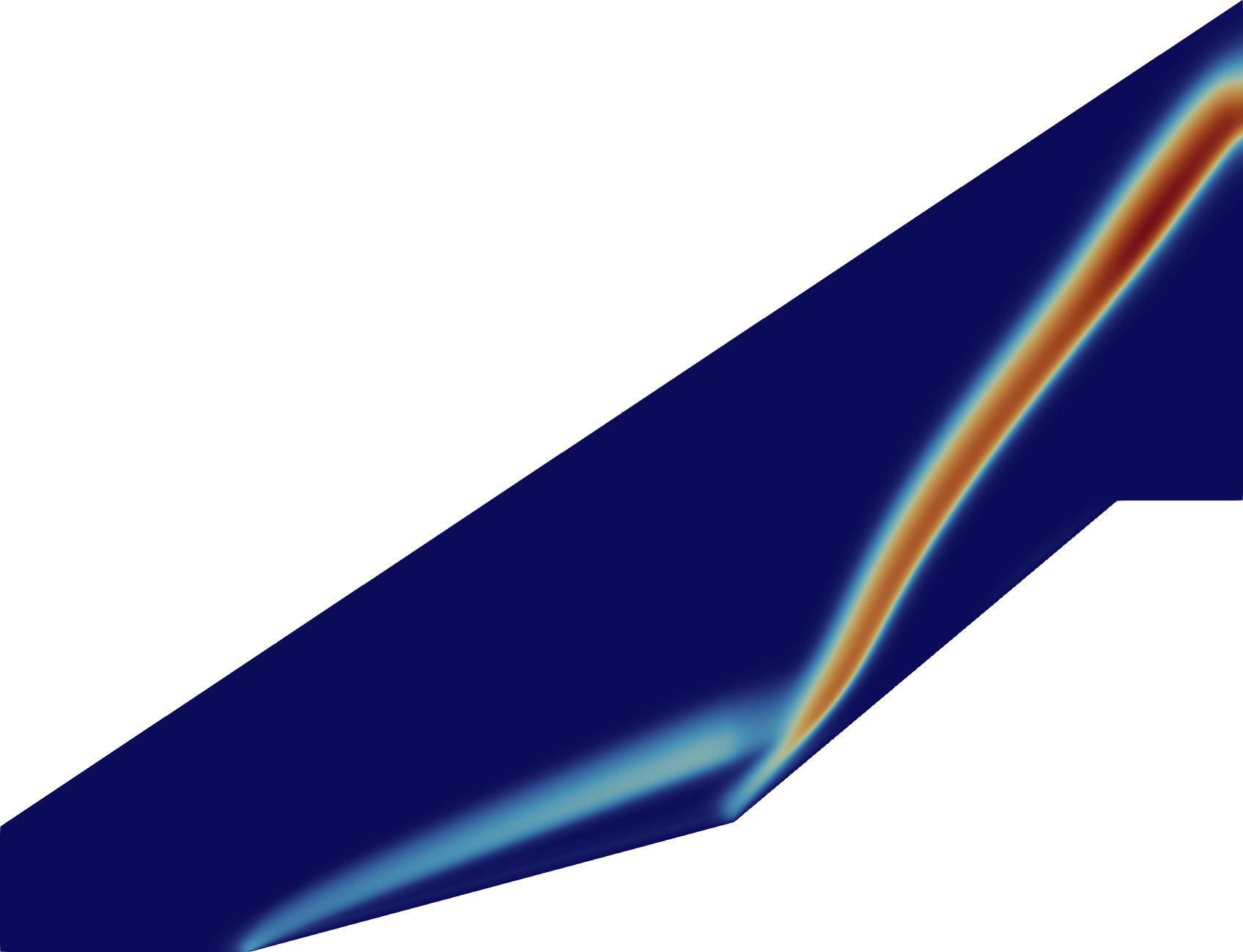}};
        \node[cell] (C) at ([xshift=0.06cm]B.north east) {\includegraphics[width=0.32\linewidth]{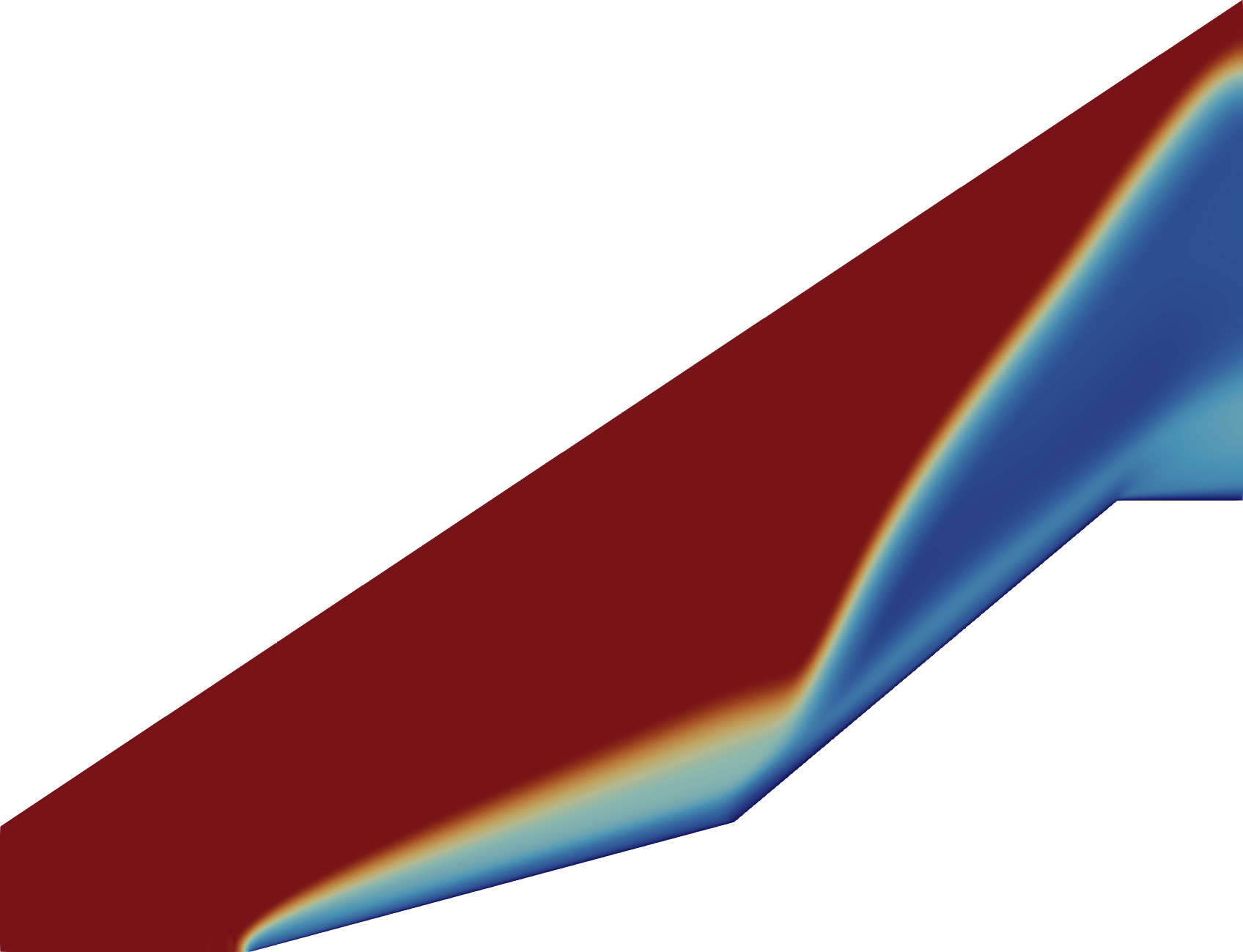}};
        \node[cell] (D) at ([yshift=-0.06cm]A.south west) {\includegraphics[width=0.32\linewidth]{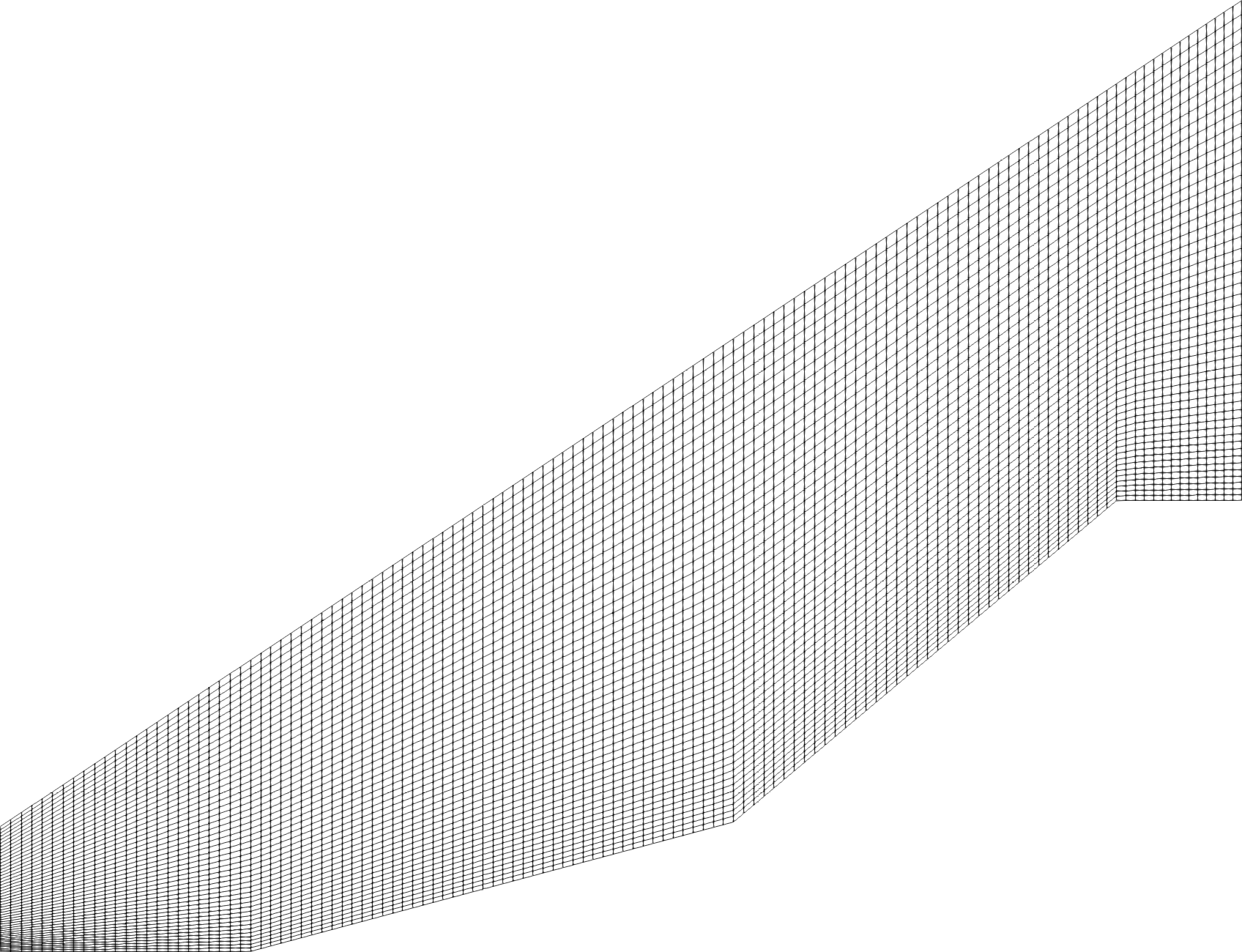}};
        \node[cell] (E) at ([xshift=0.06cm]D.north east) {\includegraphics[width=0.32\linewidth]{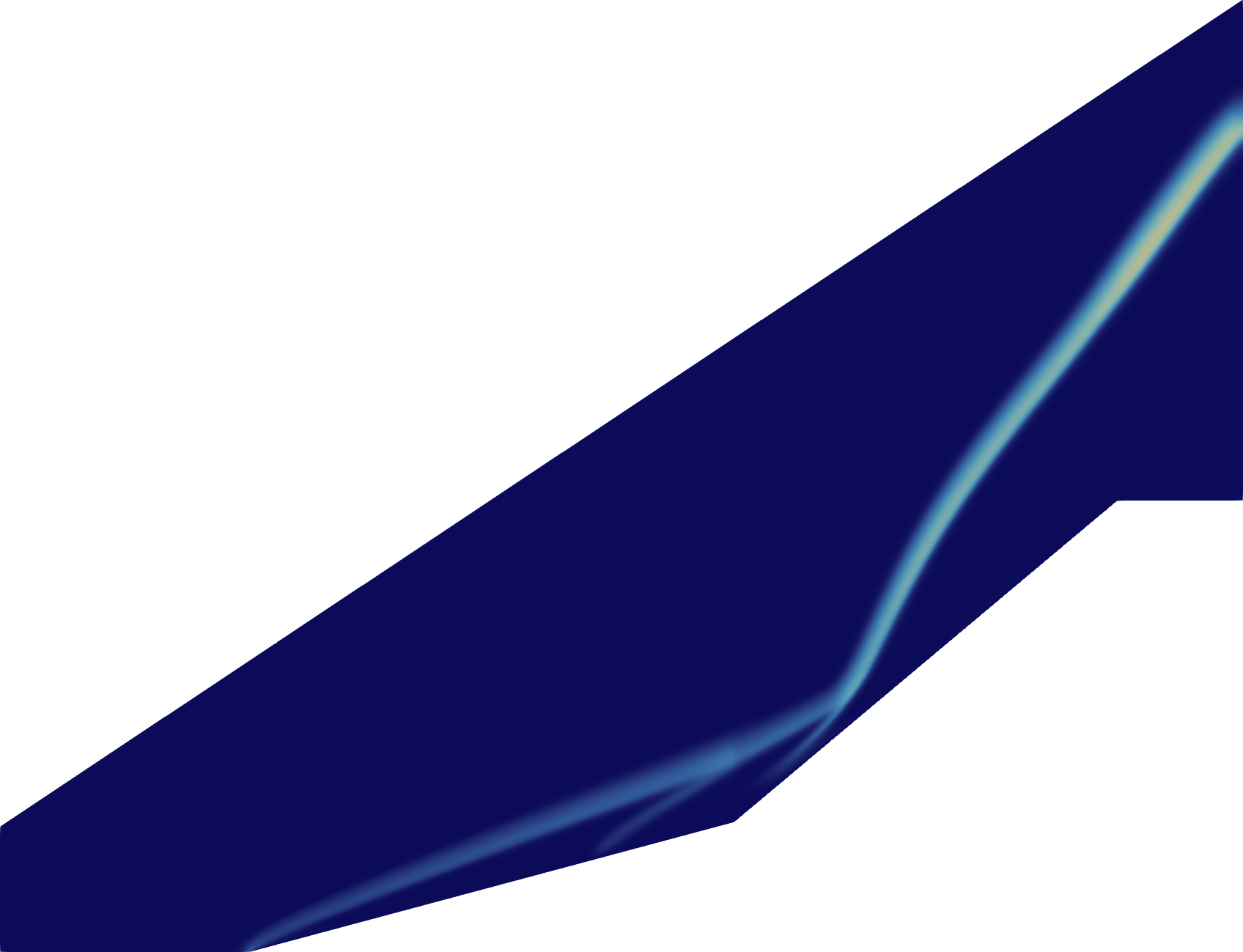}};
        \node[cell] (F) at ([xshift=0.06cm]E.north east) {\includegraphics[width=0.32\linewidth]{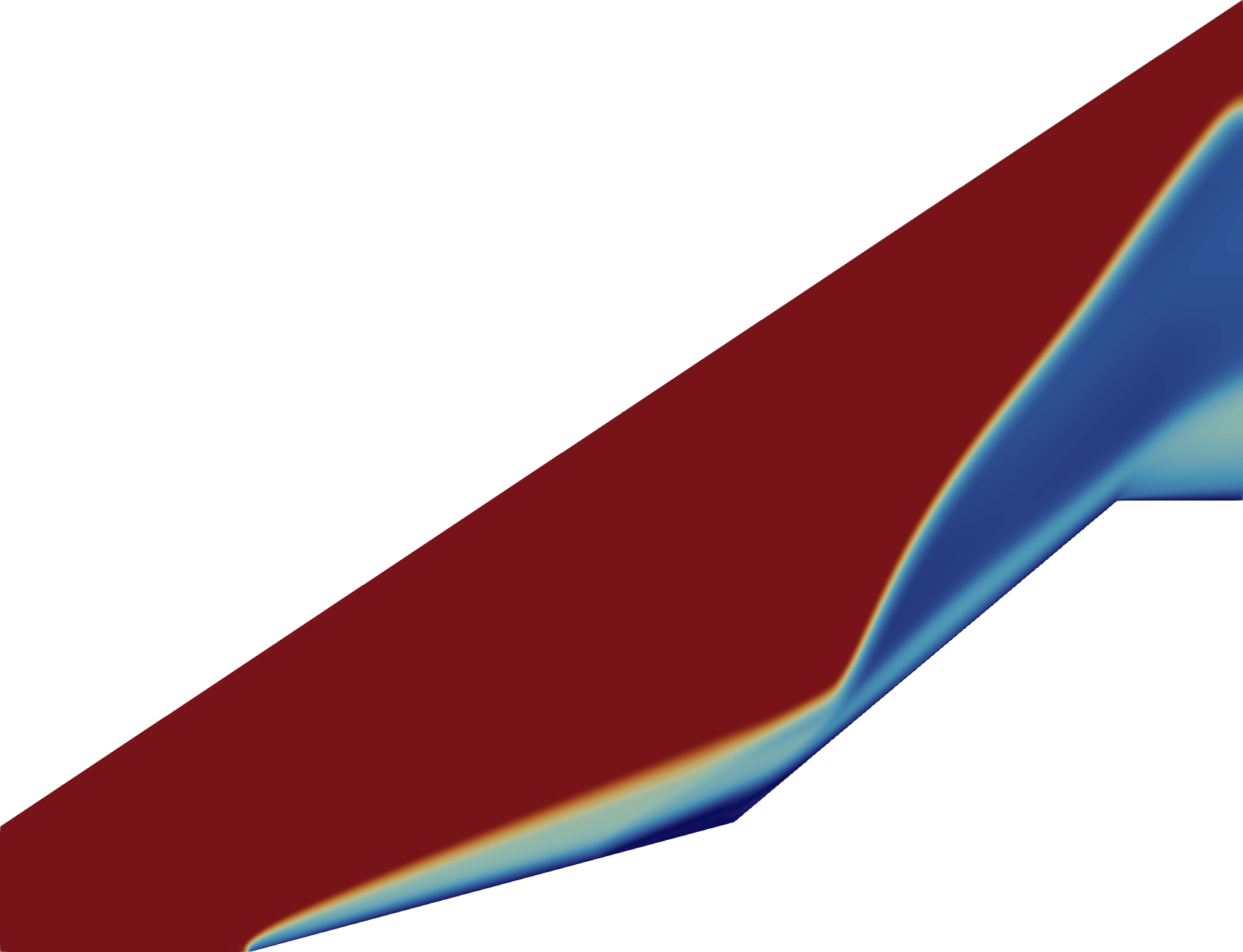}};
        \node[plabel] at (A.north west) {(a)};
        \node[plabel] at (B.north west) {(b)};
        \node[plabel] at (C.north west) {(c)};
        \node[plabel] at (D.north west) {(d)};
        \node[plabel] at (E.north west) {(e)};
        \node[plabel] at (F.north west) {(f)};
    \end{tikzpicture}
    \caption{\case{dwedge}: (a,d)~mesh, (b,e)~AV $\nu$, and (c,f)~Mach number at $(c_2, c_3) = (1.7, 1.7)$ on the (a--c)~coarse and (d--f)~refined meshes.}
    \label{fig:dw_2d_fields}
\end{figure}

\begin{figure}[hbt!]
    \centering
    \begin{tikzpicture}
\pgfplotsset{
  compat=1.18,
  table/col sep=tab,
  every axis plot/.append style={thick},
  shockline/.style={extra x ticks={0}, extra x tick labels={},
    extra x tick style={grid=major, grid style={densely dashed, black!45}}},
  machinf/.style={extra y ticks={11.68}, extra y tick labels={},
    extra y tick style={grid=major, grid style={densely dashed, black!45}}},
}
\pgfplotsset{set layers}   %
\begin{groupplot}[
  group style={
    group size=3 by 2,
    horizontal sep=0.05\textwidth,
    vertical sep=2.1cm,
    yticklabels at=edge left,
  },
  width=0.29\textwidth, height=4.3cm,
  grid=both,
  tick label style={font=\normalsize},
  machinf,
]

\nextgroupplot[title={(a) wedge 1}, ylabel={$M$}, xlabel={$s$ (from wall)},
               xmin=0, xmax=0.3, ymin=0, ymax=12,
               legend to name=dwmachlegend, legend columns=2, legend cell align=left]
\addplot[red]  table[x=s, y=Mach] {_tikz/353s_cmp_wedge1_coarse.dat};
\addlegendentry{coarse ($1488$)}
\addplot[blue] table[x=s, y=Mach] {_tikz/353s_cmp_wedge1_refined.dat};
\addlegendentry{refined ($5952$)}

\nextgroupplot[title={(b) wedge 2}, xlabel={$s$ (from wall)}, xmin=0, xmax=0.4, ymin=0, ymax=12]
\addplot[red]  table[x=s, y=Mach] {_tikz/353s_cmp_wedge2_coarse.dat};
\addplot[blue] table[x=s, y=Mach] {_tikz/353s_cmp_wedge2_refined.dat};

\nextgroupplot[title={(c) right outlet}, xlabel={$s$ (from wall)}, xmin=0, ymin=0, ymax=12]
\addplot[red]  table[x=s, y=Mach] {_tikz/353s_cmp_right_coarse.dat};
\addplot[blue] table[x=s, y=Mach] {_tikz/353s_cmp_right_refined.dat};

\nextgroupplot[title={(d)}, ylabel={$M$}, xlabel={$(s-s_\mathrm{shock})/\bar h$},
               xmax=40, ymin=0, ymax=12, shockline]
\begin{pgfonlayer}{axis background}\fill[black!22] (axis cs:-3.40,0) rectangle (axis cs:10.20,12);\end{pgfonlayer}  %
\addplot[red]  table[x expr={(\thisrow{s} - 0.0937426696081741)/0.007595291689772296},
                     y expr={\thisrow{Mach}}] {_tikz/353s_cmp_wedge1_coarse.dat};
\addplot[blue] table[x expr={(\thisrow{s} - 0.07336382838900578)/0.003695928836310913},
                     y expr={\thisrow{Mach}}] {_tikz/353s_cmp_wedge1_refined.dat};

\nextgroupplot[title={(e)}, xlabel={$(s-s_\mathrm{shock})/\bar h$}, xmax=40, ymin=0, ymax=12, shockline]
\begin{pgfonlayer}{axis background}\fill[black!22] (axis cs:0.83,0) rectangle (axis cs:10.25,12);\end{pgfonlayer}  %
\addplot[red]  table[x expr={(\thisrow{s} - 0.1429796784240934)/0.009106862150894481},
                     y expr={\thisrow{Mach}}] {_tikz/353s_cmp_wedge2_coarse.dat};
\addplot[blue] table[x expr={(\thisrow{s} - 0.11438374273927467)/0.004315252461003034},
                     y expr={\thisrow{Mach}}] {_tikz/353s_cmp_wedge2_refined.dat};

\nextgroupplot[title={(f)}, xlabel={$(s-s_\mathrm{shock})/\bar h$}, ymin=0, ymax=12, shockline]
\begin{pgfonlayer}{axis background}\fill[black!22] (axis cs:0.00,0) rectangle (axis cs:14.74,12);\end{pgfonlayer}  %
\addplot[red]  table[x expr={(\thisrow{s} - 0.763956222974278)/0.010823072602801525},
                     y expr={\thisrow{Mach}}] {_tikz/353s_cmp_right_coarse.dat};
\addplot[blue] table[x expr={(\thisrow{s} - 0.736030930228482)/0.005411536301400804},
                     y expr={\thisrow{Mach}}] {_tikz/353s_cmp_right_refined.dat};

\end{groupplot}
\end{tikzpicture}
    \figurelegend{dwmachlegend}
    \caption{\case{dwedge}: Mach number along the (a,d)~wedge~1, (b,e)~wedge~2, and (c,f)~right-outlet lines on both meshes. (a--c)~against distance from the wall; (d--f)~against $(s - s_\mathrm{shock})/\bar h$, with the captured-shock interval of Table~\ref{tab:dw_shockwidth} shaded for the coarse mesh. Dashed line: freestream $M_\infty = 11.68$.}
    \label{fig:dw_collapse_mach}
\end{figure}
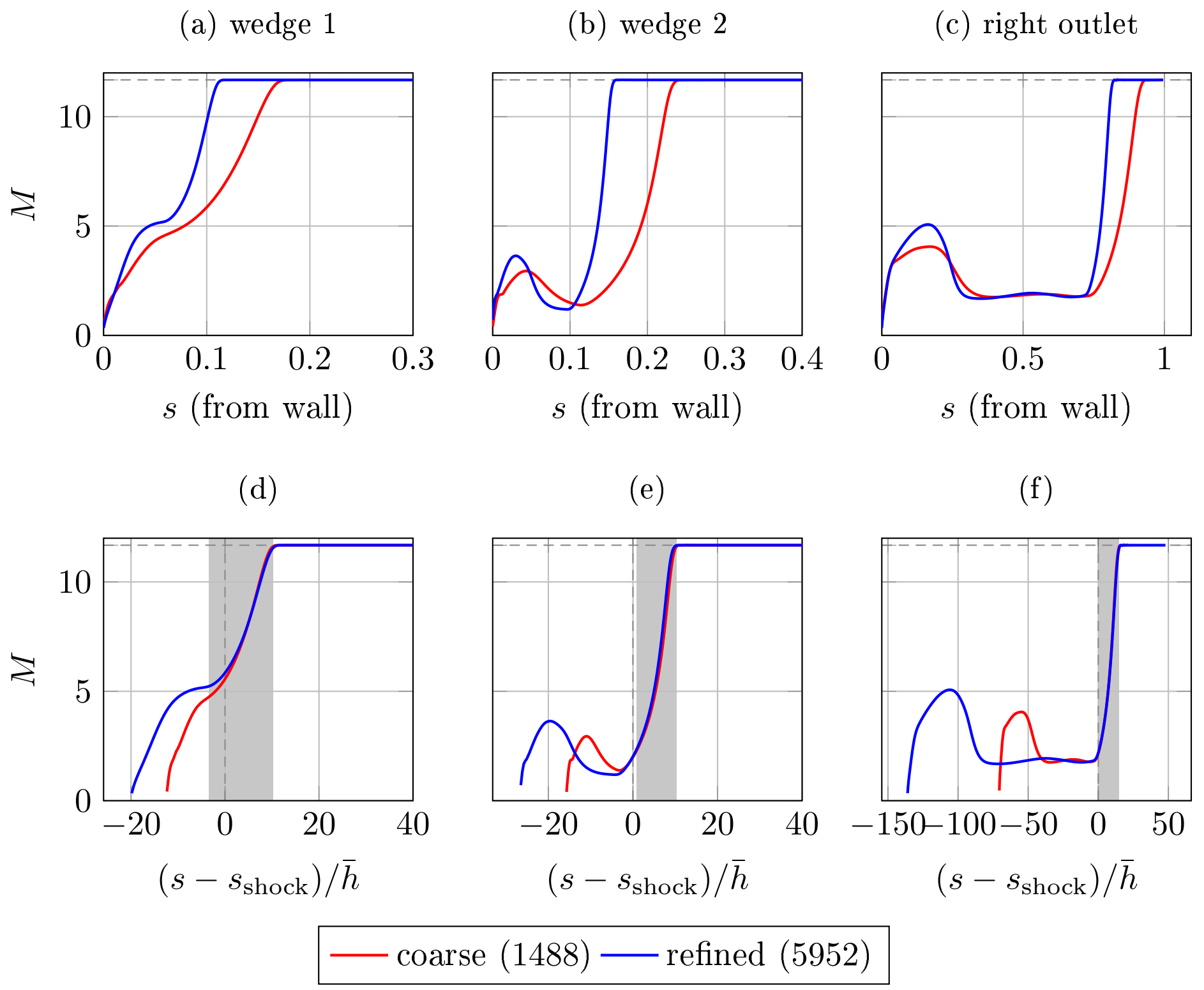

At the baseline coefficients the AV concentrates in thin bands along the captured shocks, and the refined mesh resolves them with markedly less AV (Fig.~\ref{fig:dw_2d_fields}). To test the scaling, we extract $\nu$ along the three sampling lines on both meshes and apply the shock-centered rescaling of Section~\ref{subsubsec:selfsim_cyl}: the captured-shock position $s_\mathrm{shock}$ is the arclength at which $\nu$ peaks, and profiles are shifted to $s - s_\mathrm{shock}$ and rescaled by the local resolution length $\bar h = h/p$ at the shock, with $\nu$ further normalized by $\bar h$. The raw AV peaks roughly halve between the meshes as $\bar h$ halves, and the rescaled coarse and refined profiles collapse onto a common curve on every line (Fig.~\ref{fig:dw_collapse}); the Mach transition thins under refinement in physical coordinates and collapses onto a common transition in the rescaled coordinate (Fig.~\ref{fig:dw_collapse_mach}).

\begin{figure}[hbt!]
    \centering
    \begin{tikzpicture}
\pgfplotsset{
  compat=1.18,
  table/col sep=tab,
  every axis plot/.append style={thick},
  shockline/.style={extra x ticks={0}, extra x tick labels={},
    extra x tick style={grid=major, grid style={densely dashed, black!45}}},
}
\begin{groupplot}[
  group style={
    group size=3 by 2,
    horizontal sep=0.05\textwidth,
    vertical sep=2.1cm,
    yticklabels at=edge left,
  },
  width=0.29\textwidth, height=4.3cm,
  grid=both,
  tick label style={font=\normalsize},
  scaled y ticks=false,
  yticklabel style={/pgf/number format/.cd, fixed, fixed zerofill=false, precision=3},
]

\nextgroupplot[title={(a) wedge 1}, ylabel={$\nu$}, xlabel={$s$ (from wall)},
               xmin=0, xmax=0.3, ymin=0, ymax=0.22,
               legend to name=dwcollapselegend, legend columns=2, legend cell align=left]
\addplot[red]  table[x=s, y=nu] {_tikz/353s_cmp_wedge1_coarse.dat};
\addlegendentry{coarse ($1488$)}
\addplot[blue] table[x=s, y=nu] {_tikz/353s_cmp_wedge1_refined.dat};
\addlegendentry{refined ($5952$)}

\nextgroupplot[title={(b) wedge 2}, xlabel={$s$ (from wall)}, xmin=0, xmax=0.4, ymin=0, ymax=0.22]
\addplot[red]  table[x=s, y=nu] {_tikz/353s_cmp_wedge2_coarse.dat};
\addplot[blue] table[x=s, y=nu] {_tikz/353s_cmp_wedge2_refined.dat};

\nextgroupplot[title={(c) right outlet}, xlabel={$s$ (from wall)}, xmin=0, ymin=0, ymax=0.22]
\addplot[red]  table[x=s, y=nu] {_tikz/353s_cmp_right_coarse.dat};
\addplot[blue] table[x=s, y=nu] {_tikz/353s_cmp_right_refined.dat};

\nextgroupplot[title={(d)}, ylabel={$\nu/\bar h$}, xlabel={$(s-s_\mathrm{shock})/\bar h$},
               xmax=40, ymin=0, ymax=20, shockline]
\addplot[red]  table[x expr={(\thisrow{s} - 0.0937426696081741)/0.007595291689772296},
                     y expr={\thisrow{nu}/0.007595291689772296}] {_tikz/353s_cmp_wedge1_coarse.dat};
\addplot[blue] table[x expr={(\thisrow{s} - 0.07336382838900578)/0.003695928836310913},
                     y expr={\thisrow{nu}/0.003695928836310913}] {_tikz/353s_cmp_wedge1_refined.dat};

\nextgroupplot[title={(e)}, xlabel={$(s-s_\mathrm{shock})/\bar h$}, xmax=40, ymin=0, ymax=20, shockline]
\addplot[red]  table[x expr={(\thisrow{s} - 0.1429796784240934)/0.009106862150894481},
                     y expr={\thisrow{nu}/0.009106862150894481}] {_tikz/353s_cmp_wedge2_coarse.dat};
\addplot[blue] table[x expr={(\thisrow{s} - 0.11438374273927467)/0.004315252461003034},
                     y expr={\thisrow{nu}/0.004315252461003034}] {_tikz/353s_cmp_wedge2_refined.dat};

\nextgroupplot[title={(f)}, xlabel={$(s-s_\mathrm{shock})/\bar h$}, ymin=0, ymax=20, shockline]
\addplot[red]  table[x expr={(\thisrow{s} - 0.763956222974278)/0.010823072602801525},
                     y expr={\thisrow{nu}/0.010823072602801525}] {_tikz/353s_cmp_right_coarse.dat};
\addplot[blue] table[x expr={(\thisrow{s} - 0.736030930228482)/0.005411536301400804},
                     y expr={\thisrow{nu}/0.005411536301400804}] {_tikz/353s_cmp_right_refined.dat};

\end{groupplot}
\end{tikzpicture}
    \figurelegend{dwcollapselegend}
    \caption{\case{dwedge}: AV along the (a,d)~wedge~1, (b,e)~wedge~2, and (c,f)~right-outlet lines on both meshes. (a--c)~$\nu$ against distance from the wall; (d--f)~$\nu/\bar h$ against the shifted, rescaled coordinate $(s - s_\mathrm{shock})/\bar h$.}
    \label{fig:dw_collapse}
\end{figure}
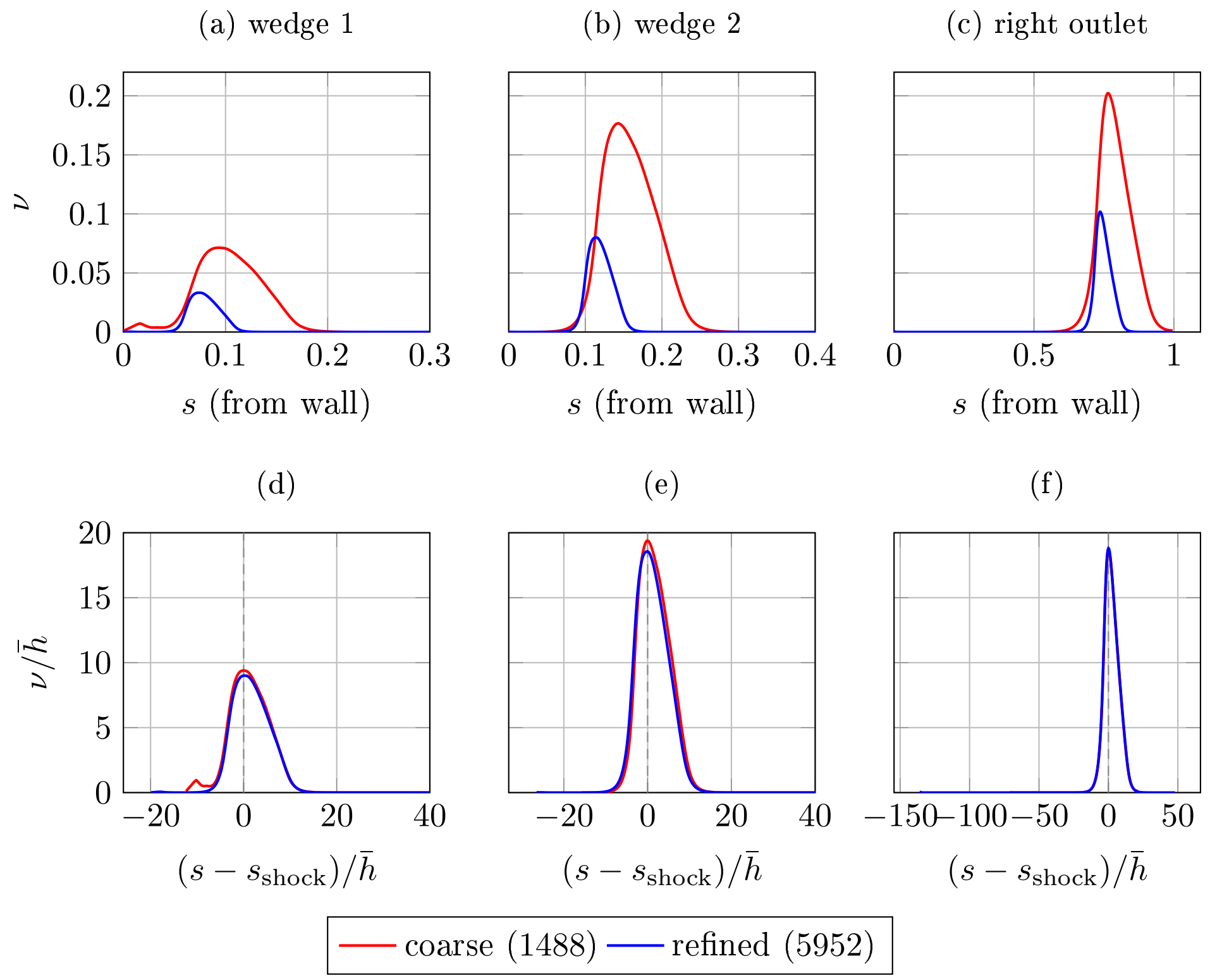

The rescaled coordinate also measures the captured shock in elements. Since $\xi = (s - s_\mathrm{shock})/\bar h$ counts arclength from the shock in units of $\bar h = h/p$, and an element spans $p$ such units, a transition that occupies a width $\Delta\xi$ in the rescaled profile occupies $N_\mathrm{elem} = \Delta\xi / p$ elements. If the collapse holds, $\Delta\xi$ is the same on both meshes and so is $N_\mathrm{elem}$: the captured shock occupies a fixed number of elements as the mesh is refined. The thickness $\delta$ of Section~\ref{subsubsec:dmin} cannot be used for $\Delta\xi$ here. It needs a reference level on each side of the shock, and on these wall-normal lines the Mach number climbs continuously from the wall through the boundary layer and the shock with no post-shock plateau. Instead, $\Delta\xi$ is the width of the interval around the steepest point of the transition over which $\mathrm{d}M/\mathrm{d}\xi$ exceeds one fifth of its maximum. The captured shock occupies $2$--$4$ elements on every line, and $N_\mathrm{elem}$ changes by less than $6\%$ between the meshes although $h$ halves (Table~\ref{tab:dw_shockwidth}).

\begin{table}[hbt!]
    \centering
    \caption{\case{dwedge}: captured-shock thickness from the rescaled Mach profiles of Fig.~\ref{fig:dw_collapse_mach} and the corresponding element count $N_{\mathrm{elem}} = \Delta\xi/p$.}
    \label{tab:dw_shockwidth}
    \begin{tabular}{l cc cc}
        \toprule
         & \multicolumn{2}{c}{coarse ($1488$ elem.)} & \multicolumn{2}{c}{refined ($5952$ elem.)} \\
        \cmidrule(lr){2-3} \cmidrule(lr){4-5}
        sampling line & $\Delta\xi$ & $N_{\mathrm{elem}}$ & $\Delta\xi$ & $N_{\mathrm{elem}}$ \\
        \midrule
        wedge~1      & $13.6$ & $3.40$ & $12.9$ & $3.22$ \\
        wedge~2      & $9.4$  & $2.36$ & $9.4$  & $2.35$ \\
        right outlet & $14.7$ & $3.69$ & $15.1$ & $3.78$ \\
        \bottomrule
    \end{tabular}
\end{table}

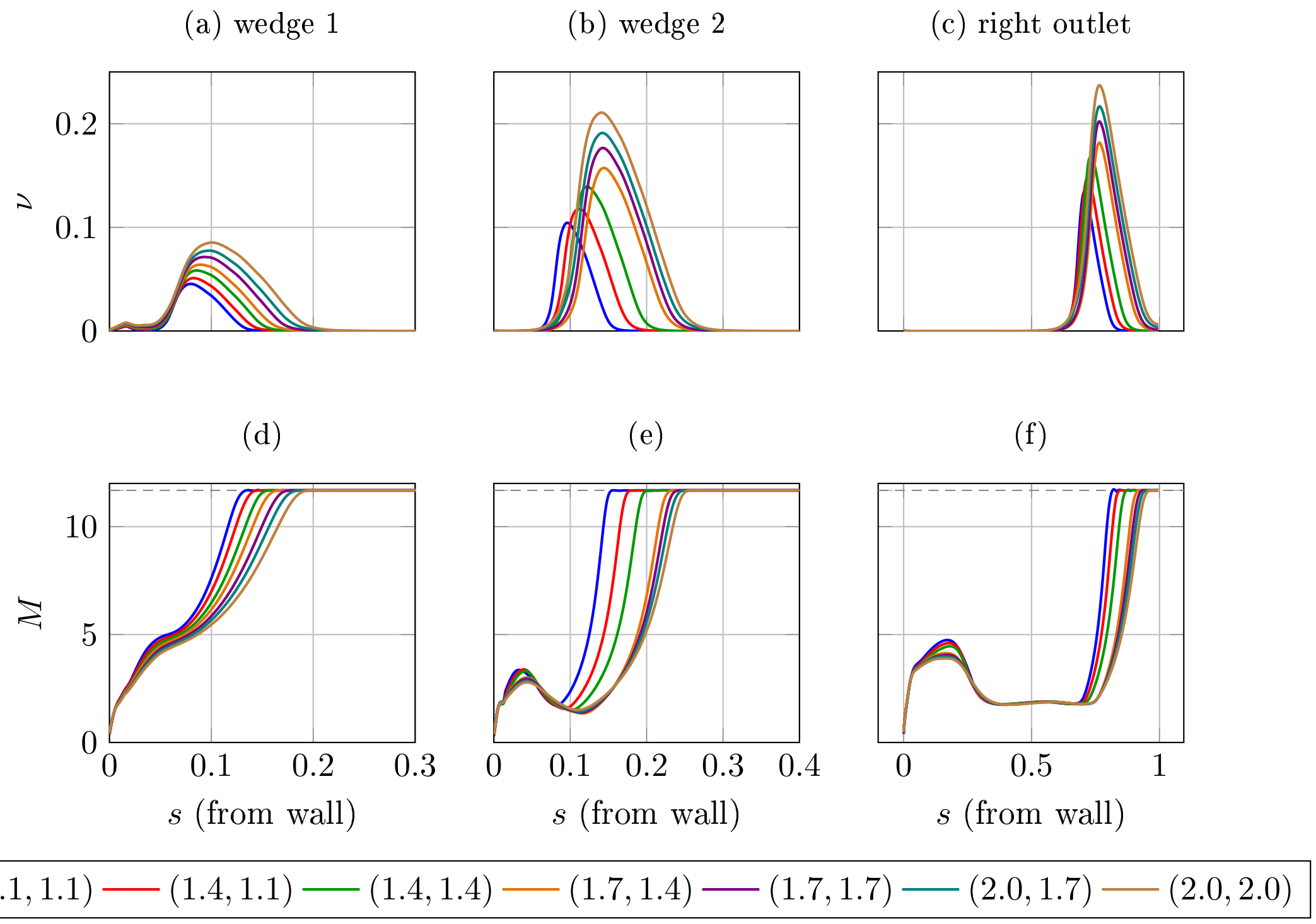
\begin{figure}[hbt!]
    \centering
    \begin{tikzpicture}
\pgfplotsset{
  compat=1.18,
  table/col sep=tab,
  every axis plot/.append style={thick},
  machinf/.style={extra y ticks={11.68}, extra y tick labels={},
    extra y tick style={grid=major, grid style={densely dashed, black!45}}},
}
\begin{groupplot}[
  group style={
    group size=3 by 2,
    horizontal sep=0.05\textwidth,
    vertical sep=1.6cm,
    yticklabels at=edge left,
  },
  width=0.29\textwidth, height=4.3cm,
  grid=both,
  tick label style={font=\normalsize},
  scaled y ticks=false,
  yticklabel style={/pgf/number format/.cd, fixed, fixed zerofill=false, precision=3},
]

\nextgroupplot[title={(a) wedge 1}, ylabel={$\nu$}, xmin=0, xmax=0.3, ymin=0, ymax=0.25, xticklabels={},
               legend to name=dwsweeplegend, legend columns=-1, legend cell align=left]
\addplot[blue]            table[x=s, y=nu] {_tikz/352s_prof_wedge1_c2_1p1_c3_1p1.dat};
\addlegendentry{$(1.1,1.1)$}
\addplot[red]             table[x=s, y=nu] {_tikz/352s_prof_wedge1_c2_1p4_c3_1p1.dat};
\addlegendentry{$(1.4,1.1)$}
\addplot[green!60!black]  table[x=s, y=nu] {_tikz/352s_prof_wedge1_c2_1p4_c3_1p4.dat};
\addlegendentry{$(1.4,1.4)$}
\addplot[orange!90!black] table[x=s, y=nu] {_tikz/352s_prof_wedge1_c2_1p7_c3_1p4.dat};
\addlegendentry{$(1.7,1.4)$}
\addplot[violet]          table[x=s, y=nu] {_tikz/352s_prof_wedge1_c2_1p7_c3_1p7.dat};
\addlegendentry{$(1.7,1.7)$}
\addplot[teal]            table[x=s, y=nu] {_tikz/352s_prof_wedge1_c2_2p0_c3_1p7.dat};
\addlegendentry{$(2.0,1.7)$}
\addplot[brown]           table[x=s, y=nu] {_tikz/352s_prof_wedge1_c2_2p0_c3_2p0.dat};
\addlegendentry{$(2.0,2.0)$}

\nextgroupplot[title={(b) wedge 2}, xmin=0, xmax=0.4, ymin=0, ymax=0.25, xticklabels={}]
\addplot[blue]            table[x=s, y=nu] {_tikz/352s_prof_wedge2_c2_1p1_c3_1p1.dat};
\addplot[red]             table[x=s, y=nu] {_tikz/352s_prof_wedge2_c2_1p4_c3_1p1.dat};
\addplot[green!60!black]  table[x=s, y=nu] {_tikz/352s_prof_wedge2_c2_1p4_c3_1p4.dat};
\addplot[orange!90!black] table[x=s, y=nu] {_tikz/352s_prof_wedge2_c2_1p7_c3_1p4.dat};
\addplot[violet]          table[x=s, y=nu] {_tikz/352s_prof_wedge2_c2_1p7_c3_1p7.dat};
\addplot[teal]            table[x=s, y=nu] {_tikz/352s_prof_wedge2_c2_2p0_c3_1p7.dat};
\addplot[brown]           table[x=s, y=nu] {_tikz/352s_prof_wedge2_c2_2p0_c3_2p0.dat};

\nextgroupplot[title={(c) right outlet}, ymin=0, ymax=0.25, xticklabels={}]
\addplot[blue]            table[x=s, y=nu] {_tikz/352s_prof_right_c2_1p1_c3_1p1.dat};
\addplot[red]             table[x=s, y=nu] {_tikz/352s_prof_right_c2_1p4_c3_1p1.dat};
\addplot[green!60!black]  table[x=s, y=nu] {_tikz/352s_prof_right_c2_1p4_c3_1p4.dat};
\addplot[orange!90!black] table[x=s, y=nu] {_tikz/352s_prof_right_c2_1p7_c3_1p4.dat};
\addplot[violet]          table[x=s, y=nu] {_tikz/352s_prof_right_c2_1p7_c3_1p7.dat};
\addplot[teal]            table[x=s, y=nu] {_tikz/352s_prof_right_c2_2p0_c3_1p7.dat};
\addplot[brown]           table[x=s, y=nu] {_tikz/352s_prof_right_c2_2p0_c3_2p0.dat};

\nextgroupplot[title={(d)}, ylabel={$M$}, xlabel={$s$ (from wall)}, xmin=0, xmax=0.3, ymin=0, ymax=12, machinf]
\addplot[blue]            table[x=s, y=Mach] {_tikz/352s_prof_wedge1_c2_1p1_c3_1p1.dat};
\addplot[red]             table[x=s, y=Mach] {_tikz/352s_prof_wedge1_c2_1p4_c3_1p1.dat};
\addplot[green!60!black]  table[x=s, y=Mach] {_tikz/352s_prof_wedge1_c2_1p4_c3_1p4.dat};
\addplot[orange!90!black] table[x=s, y=Mach] {_tikz/352s_prof_wedge1_c2_1p7_c3_1p4.dat};
\addplot[violet]          table[x=s, y=Mach] {_tikz/352s_prof_wedge1_c2_1p7_c3_1p7.dat};
\addplot[teal]            table[x=s, y=Mach] {_tikz/352s_prof_wedge1_c2_2p0_c3_1p7.dat};
\addplot[brown]           table[x=s, y=Mach] {_tikz/352s_prof_wedge1_c2_2p0_c3_2p0.dat};

\nextgroupplot[title={(e)}, xlabel={$s$ (from wall)}, xmin=0, xmax=0.4, ymin=0, ymax=12, machinf]
\addplot[blue]            table[x=s, y=Mach] {_tikz/352s_prof_wedge2_c2_1p1_c3_1p1.dat};
\addplot[red]             table[x=s, y=Mach] {_tikz/352s_prof_wedge2_c2_1p4_c3_1p1.dat};
\addplot[green!60!black]  table[x=s, y=Mach] {_tikz/352s_prof_wedge2_c2_1p4_c3_1p4.dat};
\addplot[orange!90!black] table[x=s, y=Mach] {_tikz/352s_prof_wedge2_c2_1p7_c3_1p4.dat};
\addplot[violet]          table[x=s, y=Mach] {_tikz/352s_prof_wedge2_c2_1p7_c3_1p7.dat};
\addplot[teal]            table[x=s, y=Mach] {_tikz/352s_prof_wedge2_c2_2p0_c3_1p7.dat};
\addplot[brown]           table[x=s, y=Mach] {_tikz/352s_prof_wedge2_c2_2p0_c3_2p0.dat};

\nextgroupplot[title={(f)}, xlabel={$s$ (from wall)}, ymin=0, ymax=12, machinf]
\addplot[blue]            table[x=s, y=Mach] {_tikz/352s_prof_right_c2_1p1_c3_1p1.dat};
\addplot[red]             table[x=s, y=Mach] {_tikz/352s_prof_right_c2_1p4_c3_1p1.dat};
\addplot[green!60!black]  table[x=s, y=Mach] {_tikz/352s_prof_right_c2_1p4_c3_1p4.dat};
\addplot[orange!90!black] table[x=s, y=Mach] {_tikz/352s_prof_right_c2_1p7_c3_1p4.dat};
\addplot[violet]          table[x=s, y=Mach] {_tikz/352s_prof_right_c2_1p7_c3_1p7.dat};
\addplot[teal]            table[x=s, y=Mach] {_tikz/352s_prof_right_c2_2p0_c3_1p7.dat};
\addplot[brown]           table[x=s, y=Mach] {_tikz/352s_prof_right_c2_2p0_c3_2p0.dat};

\end{groupplot}
\end{tikzpicture}
    \figurelegend{dwsweeplegend}
    \caption{\case{dwedge}: (a--c)~AV and (d--f)~Mach profiles on the coarse mesh along the (a,d)~wedge~1, (b,e)~wedge~2, and (c,f)~right-outlet lines for seven $(c_2, c_3)$ pairs.}
    \label{fig:dw_sweep}
\end{figure}

On the coarse mesh, raising $(c_2, c_3)$ over seven pairs broadens and strengthens the AV band consistently on all three lines, and every Mach profile reaches the freestream value past the shock with at most minor overshoot (Fig.~\ref{fig:dw_sweep}).

\subsection{Mach 5 viscous flow over a sphere}
\label{subsec:val_sphere}

Consider Mach~5 viscous flow over a sphere at $\mathrm{Re} = 9.4375 \times 10^{5}$ based on the freestream velocity and the sphere radius ($1.8875 \times 10^{6}$ on the diameter), the conditions of Blottner~\cite{blottner1989,widhopf1972,fisher2021blottner}. The sphere diameter is $D = 0.127\,\mathrm{m}$ (Blottner's Table~1 misprints $0.027\,\mathrm{m}$). The wall is no-slip and isothermal at $T_{\mathrm{w}} = 1.308\,T_\infty$. The gas is calorically perfect, and the viscosity follows Sutherland's law $\mu/\mu_\infty = (T/T_\infty)^{3/2}(T_\infty + S)/(T + S)$ with $T_\infty = 75.61\,\mathrm{K}$ and $S = 110.4\,\mathrm{K}$. The wall heat flux is converted to $\mathrm{kW\,m^{-2}}$ with the reference flux $\rho_\infty c_\infty^3 = 466.8\,\mathrm{kW\,m^{-2}}$, where $c_\infty = 174.3\,\mathrm{m\,s^{-1}}$ is the freestream sound speed and $\rho_\infty = 0.08813\,\mathrm{kg\,m^{-3}}$ is the measured freestream density of the experiment~\cite{widhopf1972}. The AV filter uses $a = 10^{-6}\mu_\infty$ and $\nu_0 = 0$. Lengths are nondimensionalized by the sphere radius: the wall is the unit sphere centered at the origin, the freestream runs along $+x$, and the stagnation point is $(-1, 0, 0)$. The outer boundary is the sphere centered at $(1.696, 0, 0)$ with radius $3.446$, which stands $0.75$ radii ahead of the nose and $2.0$ from the wall at the equator, and carries the freestream state; the domain is truncated at $x = 0$ by an outflow that extrapolates the state and sets the viscous flux to zero and reduced to an octant by the symmetry planes $y = 0$ and $z = 0$. The AV working variable is pinned to $\hat\nu = 0$ at the wall and on the outer boundary, where the incoming freestream carries no artificial viscosity, and takes the homogeneous Neumann condition of Section~\ref{subsec:av_bnd} on the symmetry and outflow planes.

The experiment (test case: \case{sphere}) tests whether the self-similar scaling of the AV carries over to three dimensions and how closely the wall heat flux, wall pressure, and shock standoff match the reference values~\cite{blottner1989}. We use $p = 2$ on a sequence of three nested meshes of $10\,752$, $21\,504$, and $43\,008$ hexahedra, each obtained from the previous one by bisecting every radial element, with first wall spacing $\Delta r_1 = 1.0 \times 10^{-5}$ on the coarsest; all three use $(c_2, c_3) = (5.0, 2.0)$ and $c_1 = 10^{-4}$. The coarsest mesh is an equiangular cubed-sphere octant of three blocks, each with $8 \times 8$ elements on the wall and 56 along each radial ray. The radial spacing grows geometrically from $\Delta r_1$ and is about twice as fine in a band at $20\%$ of the wall-to-outer distance, where the bow shock lies.

\begin{figure}[hbt!]
    \centering
    \IfFileExists{_imgs/741s_sphere_mach_perspective.png}%
        {\includegraphics[width=0.68\linewidth]{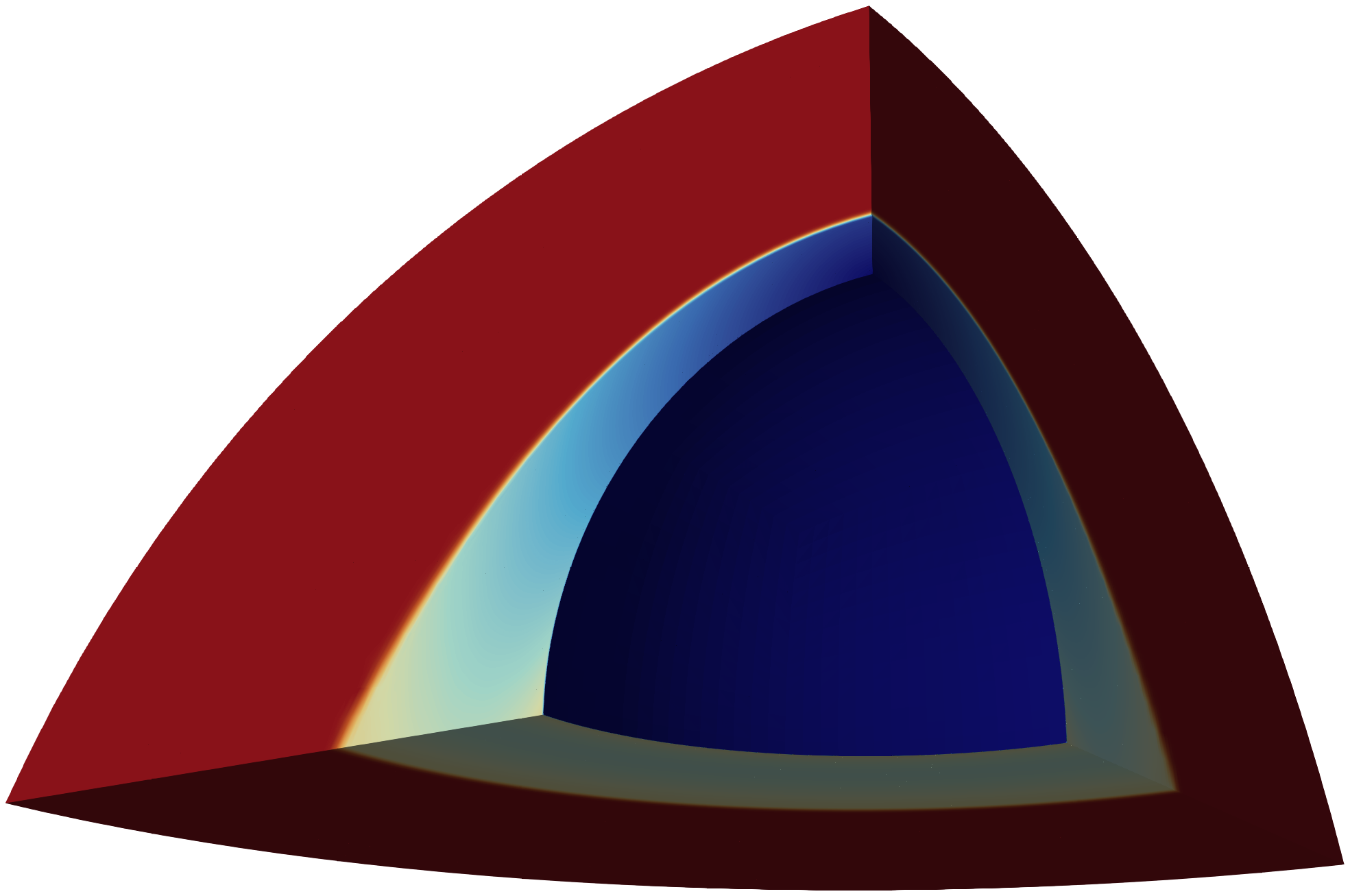}}%
        {\fbox{\parbox[c][5cm][c]{0.68\linewidth}{\centering\scriptsize sphere Mach, oblique view\\[0.3em]\textbf{[Mach on the octant's outer boundary and symmetry planes, finest mesh]}}}}
    \caption{\case{sphere}: Mach number on the finest mesh, viewed obliquely so that the outflow plane and the two symmetry planes of the octant are visible.}
    \label{fig:sphere_mach_3d}
\end{figure}

At each refinement the AV band narrows and weakens while the Mach transition sharpens (Figs.~\ref{fig:sphere_fields}--\ref{fig:sphere_mach_3d}).
\newcommand{\spherepanel}[1]{\includegraphics[width=0.315\linewidth]{#1}}   %
\begin{figure}[hbt!]
    \centering
    \begin{tikzpicture}
        \tikzset{cell/.style={inner sep=2pt, outer sep=0pt, anchor=north west}}
        \node[cell] (A) at (0,0) {\spherepanel{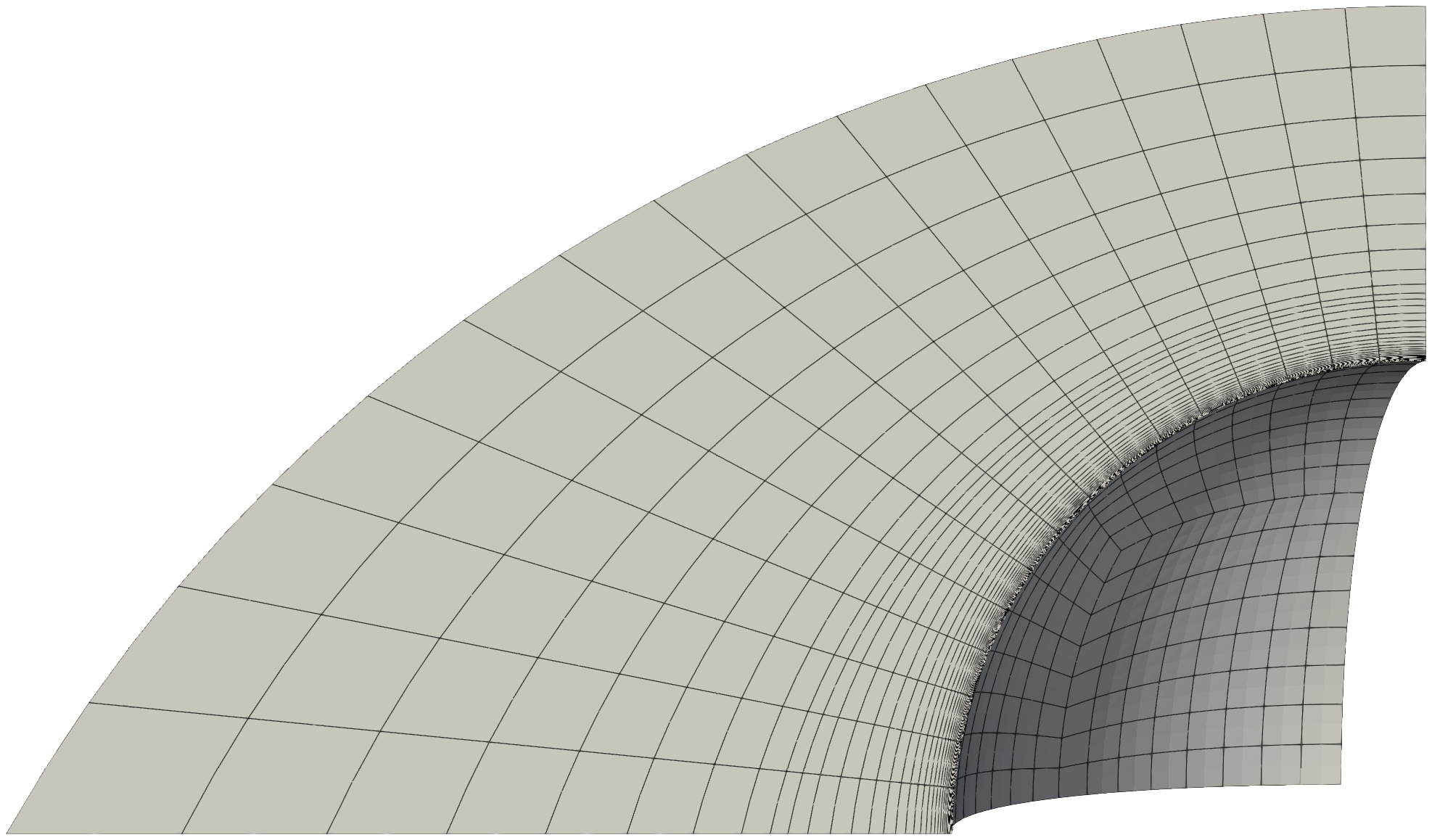}};
        \node[cell] (B) at ([xshift=0.15cm]A.north east) {\spherepanel{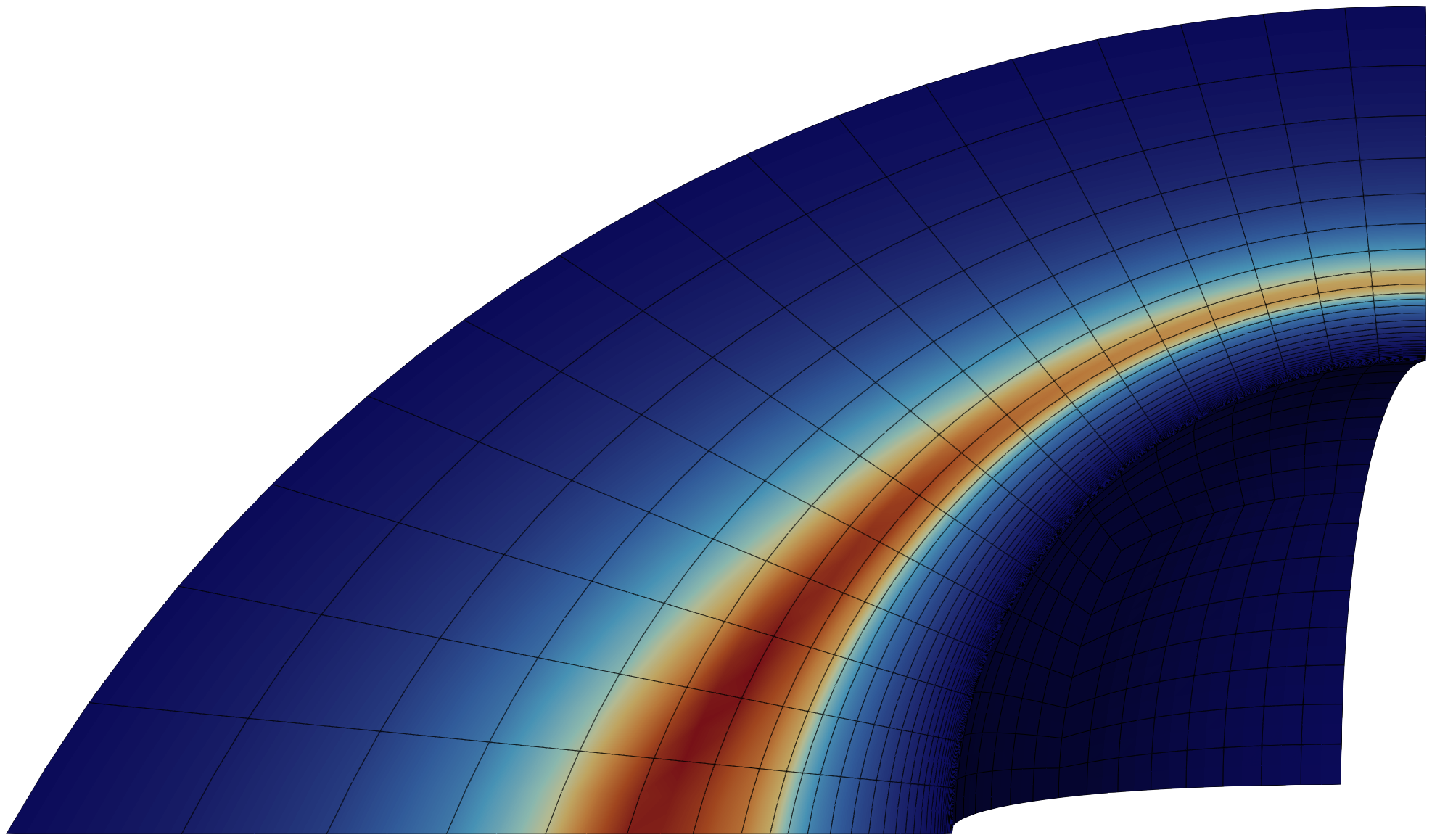}};
        \node[cell] (C) at ([xshift=0.15cm]B.north east) {\spherepanel{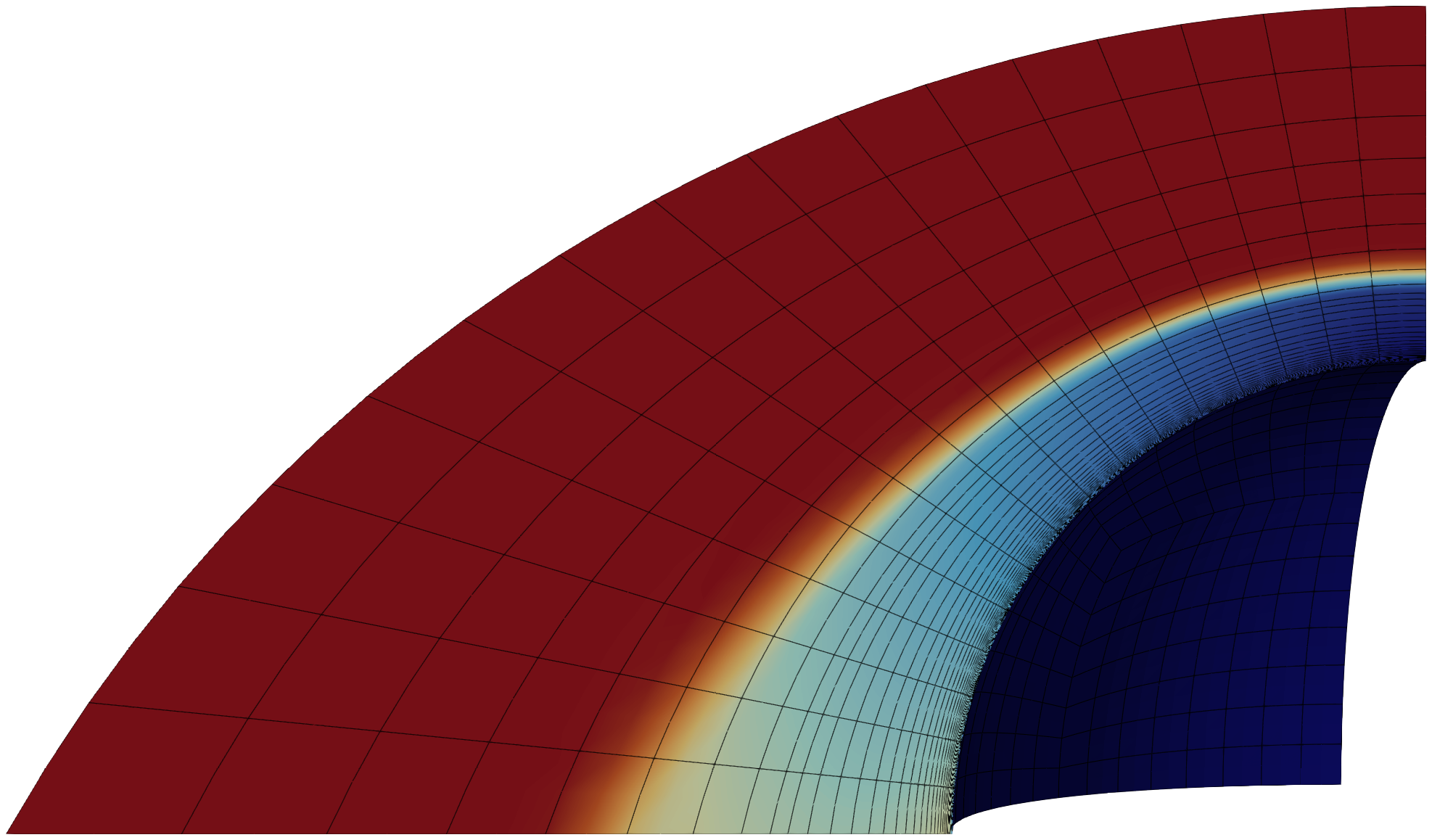}};
        \node[cell] (D) at ([yshift=-0.55cm]A.south west) {\spherepanel{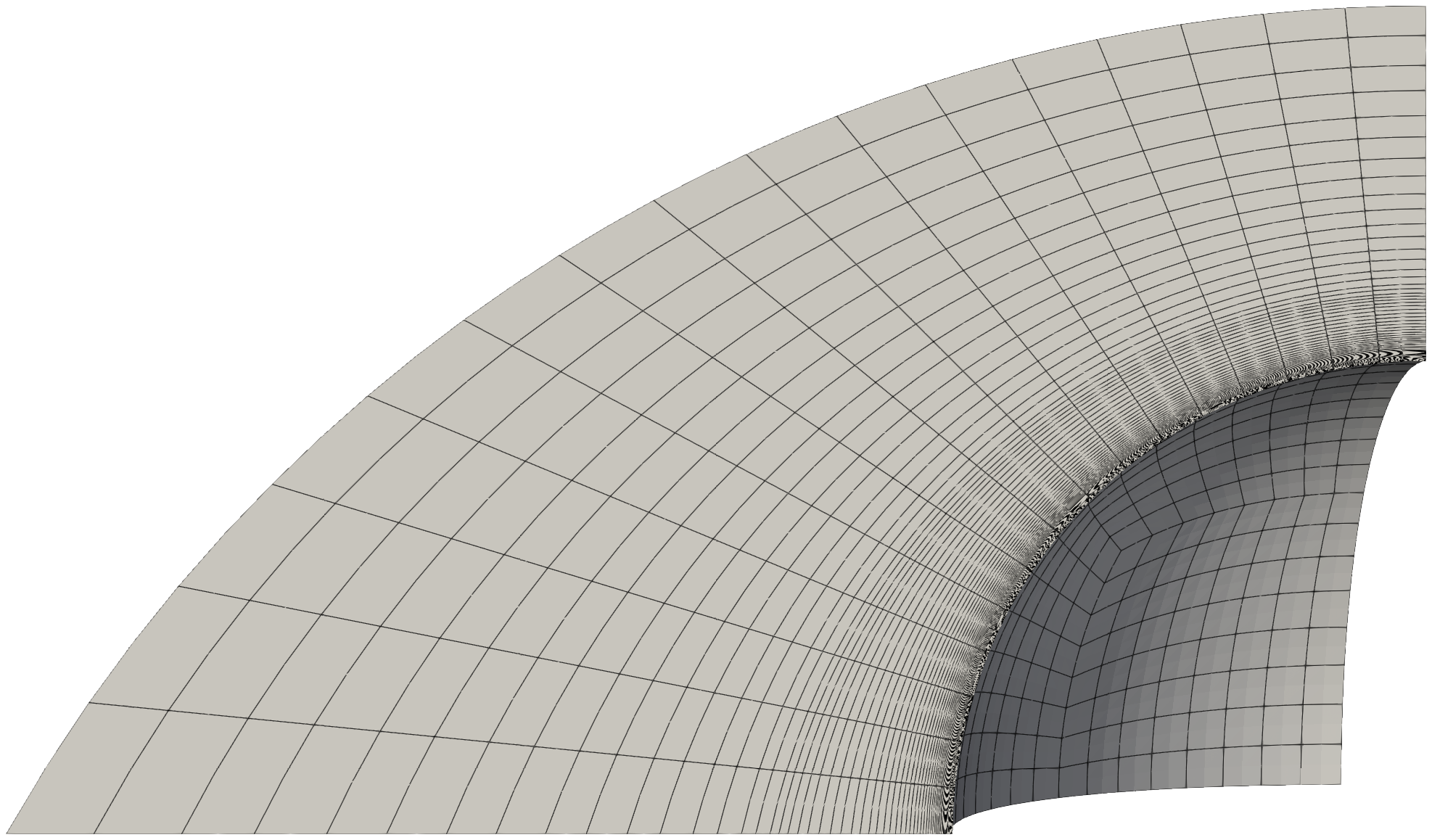}};
        \node[cell] (E) at ([xshift=0.15cm]D.north east) {\spherepanel{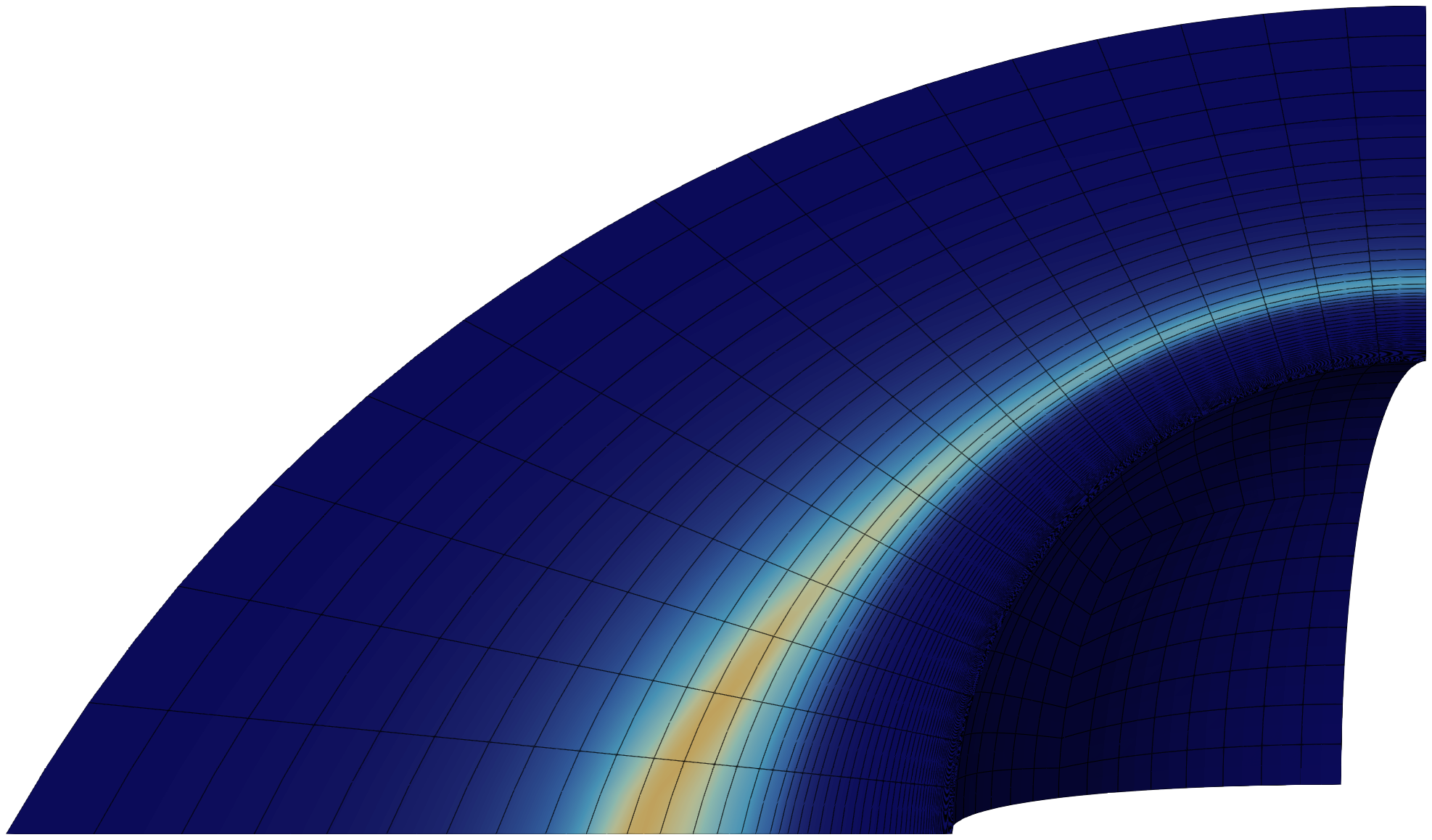}};
        \node[cell] (F) at ([xshift=0.15cm]E.north east) {\spherepanel{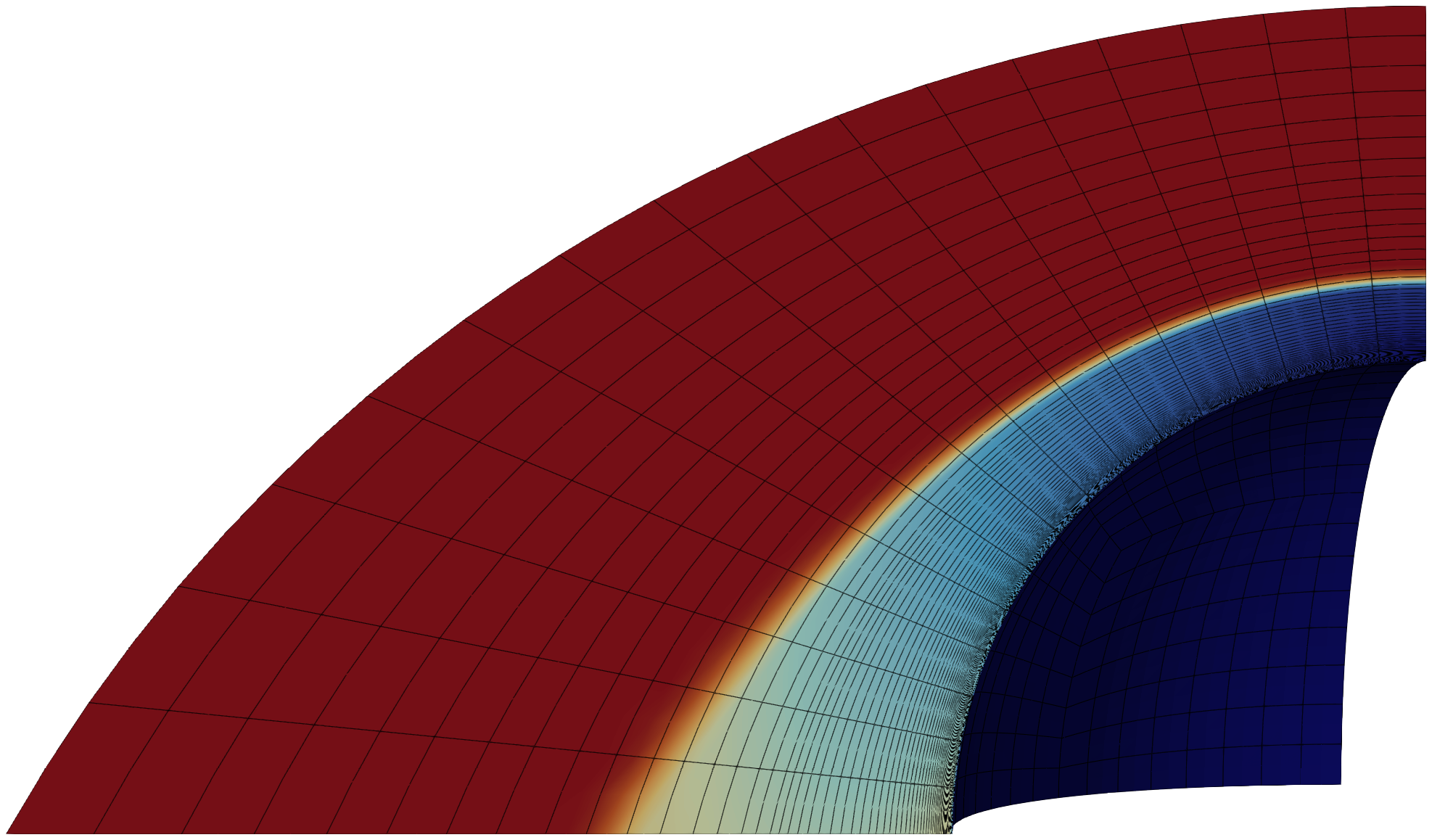}};
        \node[cell] (G) at ([yshift=-0.55cm]D.south west) {\spherepanel{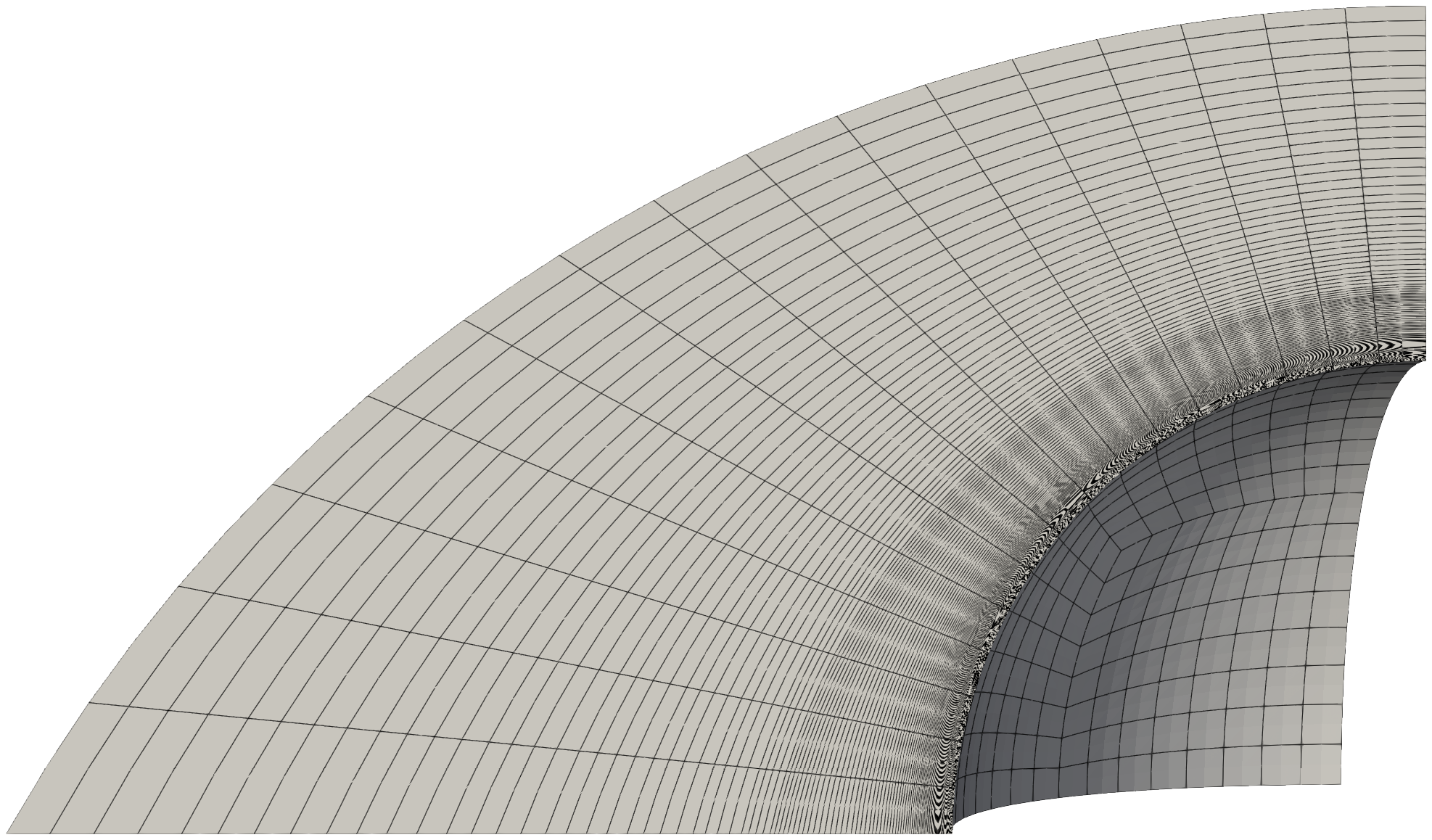}};
        \node[cell] (H) at ([xshift=0.15cm]G.north east) {\spherepanel{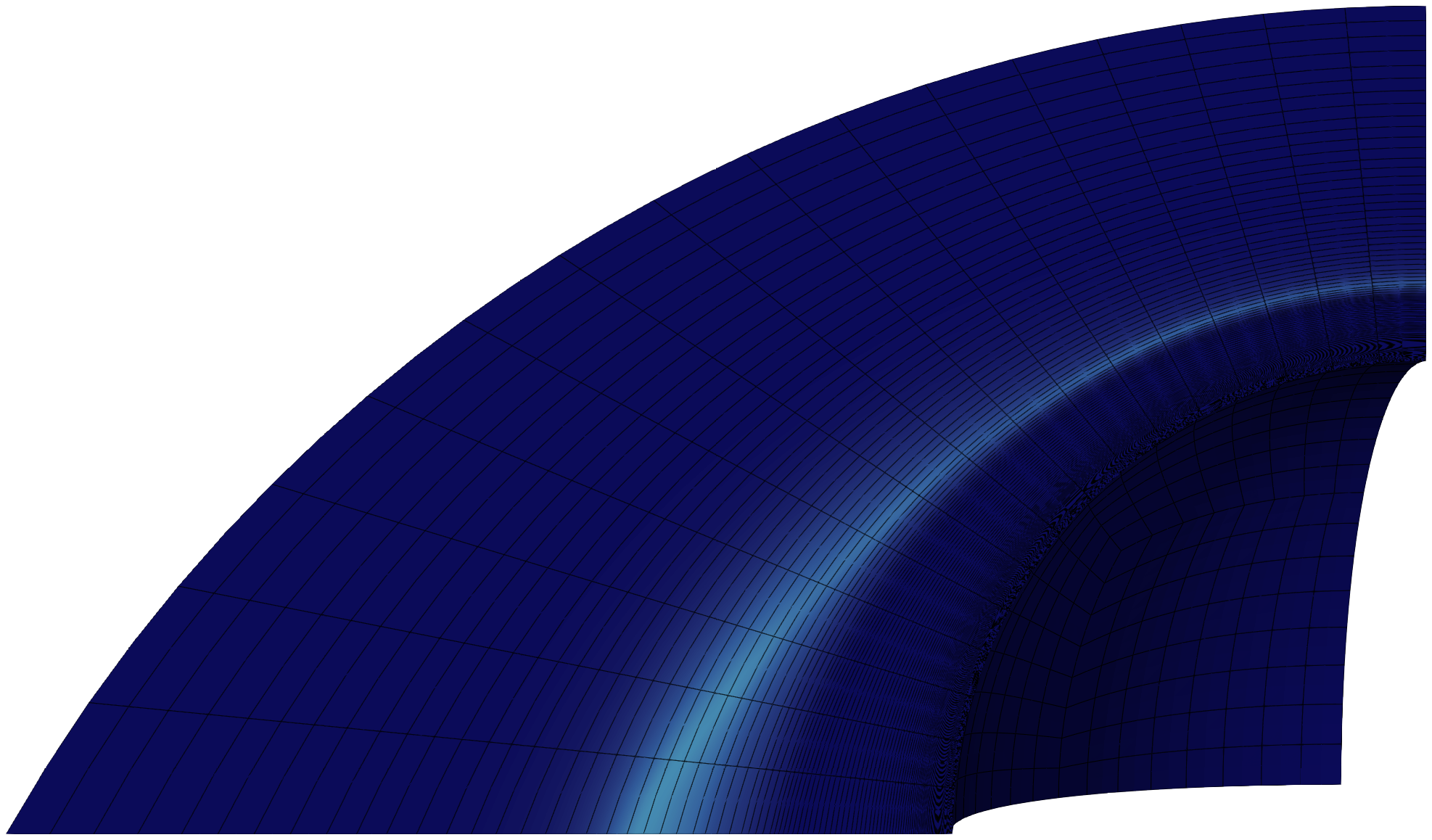}};
        \node[cell] (I) at ([xshift=0.15cm]H.north east) {\spherepanel{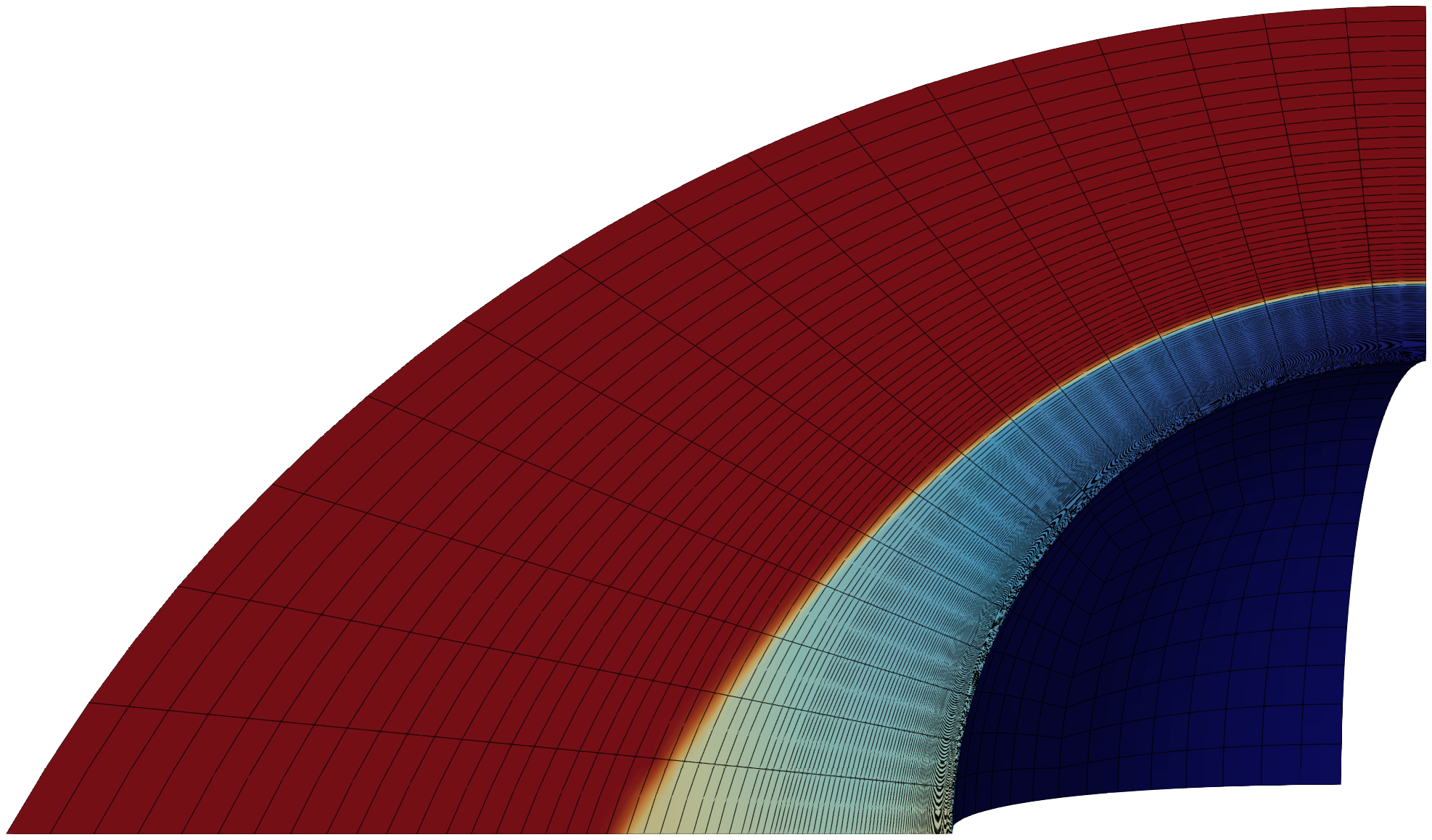}};
        \node[plabel, anchor=south west, yshift=1pt] at (A.north west) {(a)};
        \node[plabel, anchor=south west, yshift=1pt] at (B.north west) {(b)};
        \node[plabel, anchor=south west, yshift=1pt] at (C.north west) {(c)};
        \node[plabel, anchor=south west, yshift=1pt] at (D.north west) {(d)};
        \node[plabel, anchor=south west, yshift=1pt] at (E.north west) {(e)};
        \node[plabel, anchor=south west, yshift=1pt] at (F.north west) {(f)};
        \node[plabel, anchor=south west, yshift=1pt] at (G.north west) {(g)};
        \node[plabel, anchor=south west, yshift=1pt] at (H.north west) {(h)};
        \node[plabel, anchor=south west, yshift=1pt] at (I.north west) {(i)};
    \end{tikzpicture}
    \caption{\case{sphere}: (a,d,g)~mesh, (b,e,h)~AV $\nu$, and (c,f,i)~Mach number on the $z = 0$ symmetry plane for the (a--c)~coarse, (d--f)~refined, and (g--i)~finest meshes, with element edges overlaid. The AV panels share a color range set by the coarse-mesh peak.}
    \label{fig:sphere_fields}
\end{figure}

\subsubsection{Self-similar scaling of the artificial viscosity}

\begin{figure}[hbt!]
    \centering
    \IfFileExists{_tikz/sphere_collapse.tex}%
        {\begin{tikzpicture}
\pgfplotsset{set layers}   %
\pgfplotsset{
  compat=1.18,
  table/col sep=tab,
  every axis plot/.append style={thick},
  shockline/.style={extra x ticks={0}, extra x tick labels={},
    extra x tick style={grid=major, grid style={densely dashed, black!45}}},
  machinf/.style={extra y ticks={5.0}, extra y tick labels={},
    extra y tick style={grid=major, grid style={densely dashed, black!45}}},
}
\begin{groupplot}[
  group style={
    group size=2 by 2,
    horizontal sep=0.085\textwidth,
    vertical sep=1.6cm,
  },
  width=0.34\textwidth, height=4.2cm,
  grid=both,
  tick label style={font=\normalsize},
  scaled y ticks=false,
]

\nextgroupplot[title={(a)}, ylabel={$\nu$}, xlabel={$s$ (from wall)},
               xmin=0, xmax=0.45, ymin=0, ymax=0.06, ytick={0,0.02,0.04,0.06}, yticklabels={0,0.02,0.04,0.06}, legend to name=sphcollapselegend, legend columns=3, legend cell align=left]
\addplot[red] table[x=s, y=AV_deployed] {_tikz/741s_ray_stag_coarse.dat};
\addlegendentry{coarse}
\addplot[black] table[x=s, y=AV_deployed] {_tikz/741s_ray_stag_refined.dat};
\addlegendentry{refined}
\addplot[blue] table[x=s, y=AV_deployed] {_tikz/741s_ray_stag_finest.dat};
\addlegendentry{finest}

\nextgroupplot[title={(b)}, ylabel={$M$}, xlabel={$s$ (from wall)},
               xmin=0, xmax=0.45, ymin=0, ymax=6, machinf, ytick={0,3,5}]
\addplot[red] table[x=s, y=Mach] {_tikz/741s_ray_stag_coarse.dat};
\addplot[black] table[x=s, y=Mach] {_tikz/741s_ray_stag_refined.dat};
\addplot[blue] table[x=s, y=Mach] {_tikz/741s_ray_stag_finest.dat};

\nextgroupplot[title={(c)}, ylabel={$\nu/\bar h$}, xlabel={$(s-s_\mathrm{shock})/\bar h$},
               xmin=-42, xmax=40.0, ymin=0, ymax=3.4, shockline]
\addplot[red] table[x expr={(\thisrow{s} - 0.16369611639000001)/0.015589743634},
                     y expr={\thisrow{AV_deployed}/0.015589743634}] {_tikz/741s_ray_stag_coarse.dat};
\addplot[black] table[x expr={(\thisrow{s} - 0.16275803033)/0.0077948718169999999},
                     y expr={\thisrow{AV_deployed}/0.0077948718169999999}] {_tikz/741s_ray_stag_refined.dat};
\addplot[blue] table[x expr={(\thisrow{s} - 0.16181994427999999)/0.0038974359084999999},
                     y expr={\thisrow{AV_deployed}/0.0038974359084999999}] {_tikz/741s_ray_stag_finest.dat};

\nextgroupplot[title={(d)}, ylabel={$M$}, xlabel={$(s-s_\mathrm{shock})/\bar h$},
               xmin=-42, xmax=40.0, ymin=0, ymax=6, shockline, machinf, ytick={0,3,5}]
\begin{pgfonlayer}{axis background}\fill[black!22] (axis cs:-0.993,0) rectangle (axis cs:4.152,6);\end{pgfonlayer}  %
\addplot[red] table[x expr={(\thisrow{s} - 0.16369611639000001)/0.015589743634},
                     y expr={\thisrow{Mach}}] {_tikz/741s_ray_stag_coarse.dat};
\addplot[black] table[x expr={(\thisrow{s} - 0.16275803033)/0.0077948718169999999},
                     y expr={\thisrow{Mach}}] {_tikz/741s_ray_stag_refined.dat};
\addplot[blue] table[x expr={(\thisrow{s} - 0.16181994427999999)/0.0038974359084999999},
                     y expr={\thisrow{Mach}}] {_tikz/741s_ray_stag_finest.dat};
\end{groupplot}
\end{tikzpicture}}%
        {\fbox{\parbox[c][8cm][c]{0.8\linewidth}{\centering\scriptsize sphere collapse\\[0.3em]\textbf{[(a,c) raw $\nu$ and $M$ vs.\ distance from wall; (b,d) $(s-s_\mathrm{shock})/\bar h$, with $\nu/\bar h$ in (b)]}}}}
    \figurelegend{sphcollapselegend}
    \caption{\case{sphere}: stagnation-ray profiles on the three meshes. (a)~$\nu$ and (b)~Mach number against the distance $s$ from the wall; (c)~$\nu/\bar h$ and (d)~Mach number against $(s - s_\mathrm{shock})/\bar h$. Shaded in~(d): the captured-shock interval on the coarse mesh.}
    \label{fig:sphere_collapse}
\end{figure}
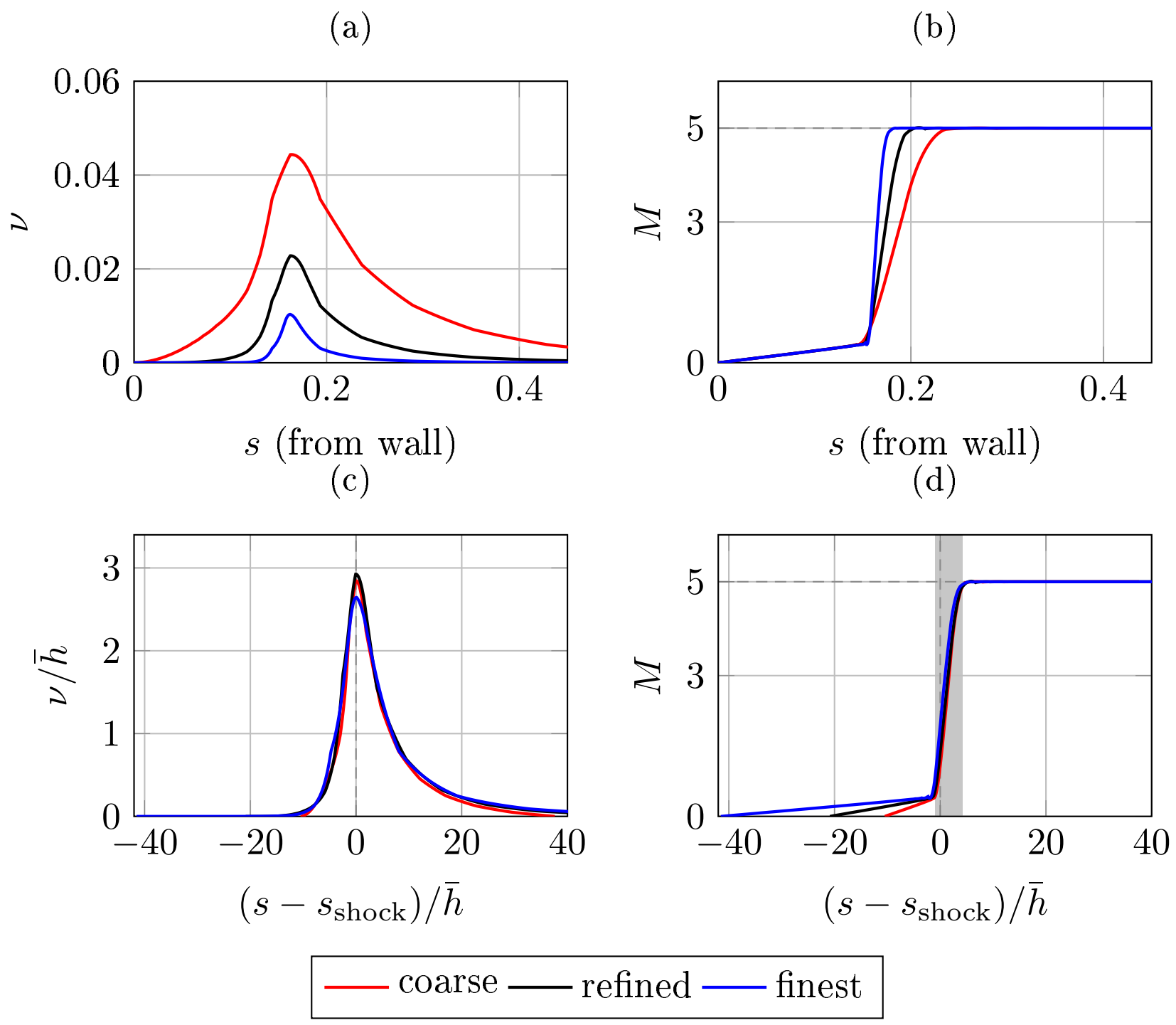

We extract $\nu$ and the Mach number along the stagnation ray and apply the shock-centered rescaling of Section~\ref{subsec:val_dwedge}, with $\bar h$ the smallest principal resolution length of the element. Because the mesh is graded radially, $\bar h$ varies along the ray, and the rescaling needs a single value per mesh; we take the value at the captured shock, $\bar h_0$, read on the coarse mesh. Because the meshes are exactly nested, with every radial element bisected, the resolution length at any fixed location halves exactly from one rung to the next, so on the finer meshes we take $\bar h_k = \bar h_0 / 2^k$ rather than re-reading it at the slightly shifted shock position. The raw AV peak falls from $4.436 \times 10^{-2}$ to $2.280 \times 10^{-2}$ to $1.030 \times 10^{-2}$ across the sequence. The successive ratios, $0.514$ and $0.452$, track the factor $0.500$ by which $\bar h$ halves, and the rescaled AV profiles fall onto a common curve through the shock region. The Mach transition thins in physical coordinates and collapses in the rescaled coordinate (Fig.~\ref{fig:sphere_collapse}). The element count $N_\mathrm{elem} = \Delta\xi/p$, measured as in Section~\ref{subsec:val_dwedge}, is preserved across a fourfold increase in radial resolution, so the captured shock occupies a fixed number of elements in three dimensions as it does in two (Table~\ref{tab:sphere_shockwidth}).

\begin{table}[hbt!]
    \centering
    \caption{\case{sphere}: captured-shock thickness $\Delta\xi$ and element count $N_{\mathrm{elem}} = \Delta\xi/p$ on the stagnation ray of Fig.~\ref{fig:sphere_collapse}(d).}
    \label{tab:sphere_shockwidth}
    \begin{tabular}{l cc cc cc}
        \toprule
         & \multicolumn{2}{c}{coarse ($10\,752$ elem.)} & \multicolumn{2}{c}{refined ($21\,504$ elem.)} & \multicolumn{2}{c}{finest ($43\,008$ elem.)} \\
        \cmidrule(lr){2-3} \cmidrule(lr){4-5} \cmidrule(lr){6-7}
        sampling line & $\Delta\xi$ & $N_{\mathrm{elem}}$ & $\Delta\xi$ & $N_{\mathrm{elem}}$ & $\Delta\xi$ & $N_{\mathrm{elem}}$ \\
        \midrule
        stagnation ray & $5.14$ & $2.57$ & $5.24$ & $2.62$ & $5.30$ & $2.65$ \\
        \bottomrule
    \end{tabular}
\end{table}

\subsubsection{Comparison with reference data}

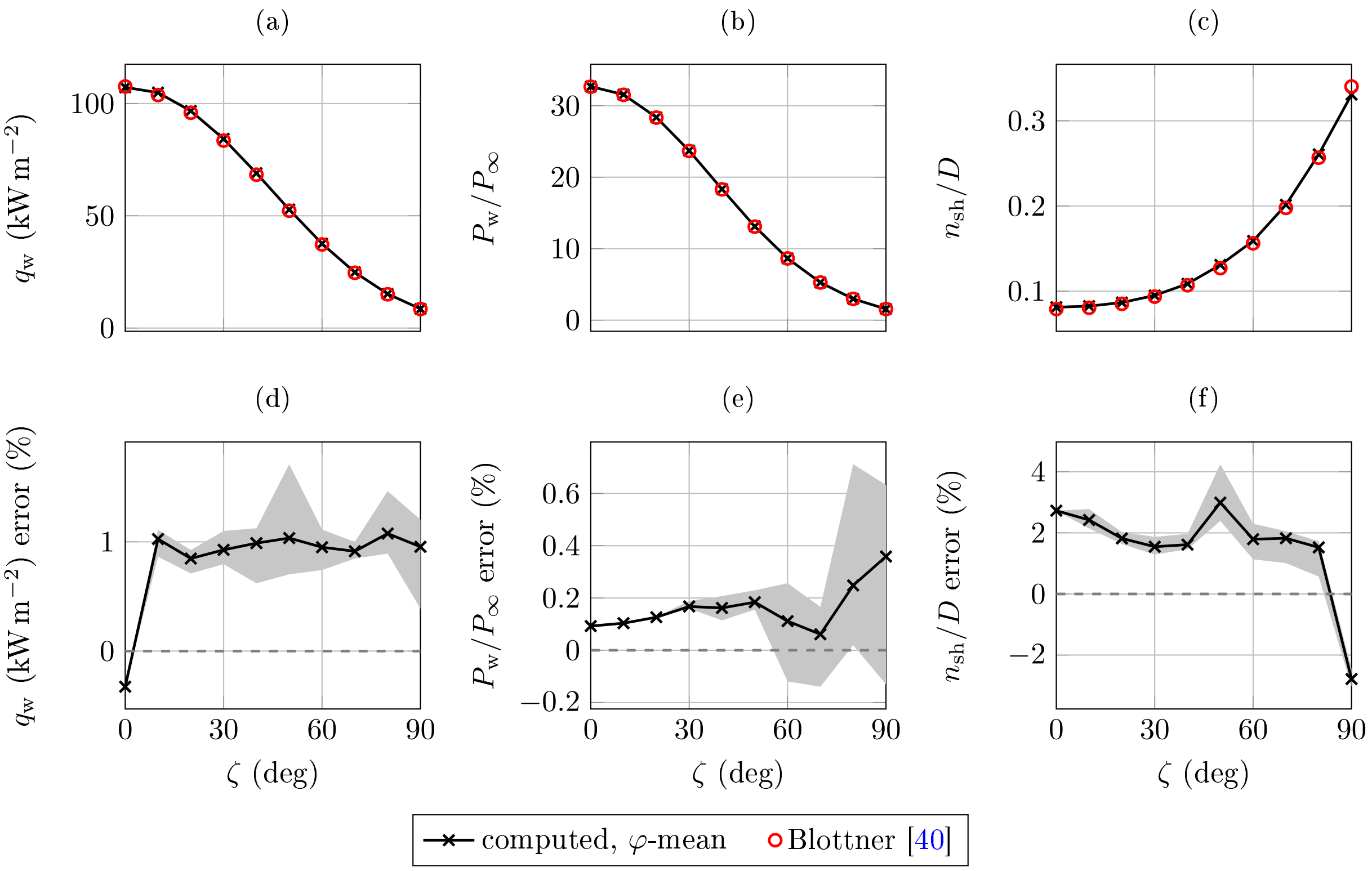
\begin{figure}[hbt!]
    \centering
    \IfFileExists{_tikz/sphere_wall.tex}%
        {\begin{tikzpicture}
\begin{groupplot}[group style={group size=3 by 2, horizontal sep=2.0cm,
    vertical sep=1.3cm}, scale only axis, width=0.21\linewidth, height=0.19\linewidth,
    grid=both, xmin=0, xmax=90, xtick={0,30,60,90},
    ylabel style={text width=}, ylabel absolute, ylabel shift=-30pt,
    every axis plot/.append style={line width=0.9pt}]
    \nextgroupplot[title={(a)}, ylabel={$q_{\mathrm{w}}$ (kW\,m$^{-2}$)}, xticklabels={}, legend to name=sphwalllegend, legend columns=2, legend style={/tikz/every even column/.append style={column sep=12pt}}]
        \addplot[draw=none, name path=qvLo, forget plot] coordinates {(0,107.221) (10,104.675) (20,96.5751) (30,84.1485) (40,68.7166) (50,52.6191) (60,37.5082) (70,24.8235) (80,15.1782) (90,8.50492)};
        \addplot[draw=none, name path=qvHi, forget plot] coordinates {(0,107.221) (10,104.929) (20,96.7814) (30,84.4022) (40,69.0607) (50,53.1459) (60,37.6464) (70,24.861) (80,15.2643) (90,8.57393)};
        \addplot[black!22, forget plot] fill between[of=qvLo and qvHi];

        \addplot[black, mark=x, mark size=2.6pt] coordinates {(0,107.221) (10,104.843) (20,96.7065) (30,84.2576) (40,68.9677) (50,52.7926) (60,37.5859) (70,24.8399) (80,15.2059) (90,8.55289)};
        \addlegendentry{computed, $\varphi$-mean}
        \addplot[only marks, mark=o, red] coordinates {(0,107.572) (10,103.778) (20,95.894) (30,83.484) (40,68.293) (50,52.252) (60,37.232) (70,24.615) (80,15.044) (90,8.472)};
        \addlegendentry{Blottner~\cite{blottner1989}}
    \nextgroupplot[title={(b)}, ylabel={$P_{\mathrm{w}}/P_\infty$}, xticklabels={}]
        \addplot[draw=none, name path=pvLo, forget plot] coordinates {(0,32.6848) (10,31.5565) (20,28.3752) (30,23.7039) (40,18.3263) (50,13.1049) (60,8.62601) (70,5.26055) (80,2.97598) (90,1.55886)};
        \addplot[draw=none, name path=pvHi, forget plot] coordinates {(0,32.6848) (10,31.5571) (20,28.3757) (30,23.7108) (40,18.3433) (50,13.1146) (60,8.65846) (70,5.27665) (80,2.99658) (90,1.57077)};
        \addplot[black!22, forget plot] fill between[of=pvLo and pvHi];

        \addplot[black, mark=x, mark size=2.6pt] coordinates {(0,32.6848) (10,31.5567) (20,28.3754) (30,23.7057) (40,18.335) (50,13.1087) (60,8.64592) (70,5.27112) (80,2.98278) (90,1.5665)};
        \addplot[only marks, mark=o, red] coordinates {(0,32.6546) (10,31.5241) (20,28.3396) (30,23.6661) (40,18.3054) (50,13.0846) (60,8.6363) (70,5.2679) (80,2.9754) (90,1.5609)};
    \nextgroupplot[title={(c)}, ylabel={$n_{\mathrm{sh}}/D$}, xticklabels={}]
        \addplot[draw=none, name path=svLo, forget plot] coordinates {(0,0.0811343) (10,0.0821757) (20,0.0864619) (30,0.0947168) (40,0.108474) (50,0.130027) (60,0.157932) (70,0.199766) (80,0.25814) (90,0.329959)};
        \addplot[draw=none, name path=svHi, forget plot] coordinates {(0,0.0811343) (10,0.082686) (20,0.0867941) (30,0.0952572) (40,0.109013) (50,0.132385) (60,0.159752) (70,0.201828) (80,0.261119) (90,0.331958)};
        \addplot[black!22, forget plot] fill between[of=svLo and svHi];

        \addplot[black, mark=x, mark size=2.6pt] coordinates {(0,0.0811343) (10,0.082398) (20,0.0866169) (30,0.0949567) (40,0.108637) (50,0.130787) (60,0.158959) (70,0.201377) (80,0.260591) (90,0.330816)};
        \addplot[only marks, mark=o, red] coordinates {(0,0.07898) (10,0.08045) (20,0.08507) (30,0.09351) (40,0.10691) (50,0.12699) (60,0.15616) (70,0.19777) (80,0.25668) (90,0.34029)};
    \nextgroupplot[title={(d)}, ylabel={$q_{\mathrm{w}}$ (kW\,m$^{-2}$) error (\%)}, xlabel={$\zeta$ (deg)}]
        \addplot[draw=none, name path=qeLo, forget plot] coordinates {(0,-0.326241) (10,0.864051) (20,0.710292) (30,0.79598) (40,0.620247) (50,0.702477) (60,0.741886) (70,0.847229) (80,0.892156) (90,0.388549)};
        \addplot[draw=none, name path=qeHi, forget plot] coordinates {(0,-0.326241) (10,1.10918) (20,0.925421) (30,1.09983) (40,1.12412) (50,1.71066) (60,1.11296) (70,0.999334) (80,1.46428) (90,1.20308)};
        \addplot[black!22, forget plot] fill between[of=qeLo and qeHi];

        \addplot[black, mark=x, mark size=2.6pt] coordinates {(0,-0.326241) (10,1.02592) (20,0.847239) (30,0.926666) (40,0.987988) (50,1.03457) (60,0.950518) (70,0.913821) (80,1.07613) (90,0.954838)};
        \addplot[dashed, gray, forget plot] coordinates {(0,0) (90,0)};
    \nextgroupplot[title={(e)}, ylabel={$P_{\mathrm{w}}/P_\infty$ error (\%)}, xlabel={$\zeta$ (deg)}]
        \addplot[draw=none, name path=peLo, forget plot] coordinates {(0,0.0925684) (10,0.102832) (20,0.125602) (30,0.159669) (40,0.114425) (50,0.155412) (60,-0.119136) (70,-0.139498) (80,0.0196014) (90,-0.130466)};
        \addplot[draw=none, name path=peHi, forget plot] coordinates {(0,0.0925684) (10,0.104533) (20,0.12726) (30,0.188718) (40,0.207231) (50,0.229349) (60,0.256581) (70,0.166038) (80,0.711821) (90,0.632433)};
        \addplot[black!22, forget plot] fill between[of=peLo and peHi];

        \addplot[black, mark=x, mark size=2.6pt] coordinates {(0,0.0925684) (10,0.103511) (20,0.126258) (30,0.167485) (40,0.161827) (50,0.183925) (60,0.111401) (70,0.0612006) (80,0.24804) (90,0.358733)};
        \addplot[dashed, gray, forget plot] coordinates {(0,0) (90,0)};
    \nextgroupplot[title={(f)}, ylabel={$n_{\mathrm{sh}}/D$ error (\%)}, xlabel={$\zeta$ (deg)}]
        \addplot[draw=none, name path=seLo, forget plot] coordinates {(0,2.72763) (10,2.14504) (20,1.63615) (30,1.29059) (40,1.4627) (50,2.39185) (60,1.13505) (70,1.00941) (80,0.568749) (90,-3.03601)};
        \addplot[draw=none, name path=seHi, forget plot] coordinates {(0,2.72763) (10,2.7794) (20,2.02663) (30,1.86841) (40,1.9671) (50,4.24805) (60,2.30012) (70,2.05206) (80,1.7293) (90,-2.44835)};
        \addplot[black!22, forget plot] fill between[of=seLo and seHi];

        \addplot[black, mark=x, mark size=2.6pt] coordinates {(0,2.72763) (10,2.4214) (20,1.81841) (30,1.54706) (40,1.61493) (50,2.98969) (60,1.79247) (70,1.82386) (80,1.52385) (90,-2.78416)};
        \addplot[dashed, gray, forget plot] coordinates {(0,0) (90,0)};
\end{groupplot}
\end{tikzpicture}}%
        {\fbox{\parbox[c][8cm][c]{0.9\linewidth}{\centering\scriptsize sphere wall comparison\\[0.3em]\textbf{[(a--c) $q_{\mathrm{w}}$, $P_{\mathrm{w}}/P_\infty$, $n_{\mathrm{sh}}/D$ vs.\ $\zeta$ against Blottner; (d--f) the errors]}}}}
    \figurelegend{sphwalllegend}
    \caption{\case{sphere}: (a)~wall heat flux $q_{\mathrm{w}}$, (b)~wall pressure $P_{\mathrm{w}}/P_\infty$, and (c)~shock standoff $n_{\mathrm{sh}}/D$ on the finest mesh against the reference values~\cite{blottner1989} as functions of the polar angle $\zeta$, and (d--f)~the corresponding percentage errors. Shaded: the azimuthal spread.}
    \label{fig:sphere_wall}
\end{figure}

The wall heat flux, wall pressure, and shock standoff on the finest mesh are compared with Blottner's tabulated values~\cite{blottner1989} at ten polar stations $\zeta$ from the nose to the equator; at each station the computed value is the mean over seven azimuths, and errors are percentage differences from the reference, positive where the computation exceeds it. The wall pressure agrees to $0.16\%$ on average and nowhere worse than $0.36\%$, the heat flux to $0.90\%$ and nowhere worse than $1.08\%$, and the standoff, taken as the ridge of the deployed AV since a captured shock has no sharp position, to $2.10\%$ (Fig.~\ref{fig:sphere_wall}, Table~\ref{tab:sphere_wall}). For the heat flux and the standoff, the azimuthal spread (shaded in the figure) is smaller than the discrepancy from the reference at every station, so that discrepancy reflects the discretization error rather than the variation of the solution over the octant.

\begin{table}[hbt!]
    \centering
    \caption{The comparison of Fig.~\ref{fig:sphere_wall}, tabulated: wall heat flux, wall pressure and shock standoff on the finest mesh against the reference values~\cite{blottner1989}, with the error in percent. Computed values are means over the seven azimuthal stations at each $\zeta$.}
    \label{tab:sphere_wall}
    \begin{tabular}{r rrr rrr rrr}
        \toprule
         & \multicolumn{3}{c}{$q_{\mathrm{w}}$ (kW\,m$^{-2}$)} & \multicolumn{3}{c}{$P_{\mathrm{w}}/P_\infty$} & \multicolumn{3}{c}{$n_{\mathrm{sh}}/D$} \\
        \cmidrule(lr){2-4}\cmidrule(lr){5-7}\cmidrule(lr){8-10}
        $\zeta$ & computed & reference & err (\%) & computed & reference & err (\%) & computed & reference & err (\%) \\
        \midrule
        $0$ & $107.221$ & $107.572$ & $-0.326$ & $32.6848$ & $32.6546$ & $0.093$ & $0.08113$ & $0.07898$ & $2.73$ \\
        $10$ & $104.843$ & $103.778$ & $1.026$ & $31.5567$ & $31.5241$ & $0.104$ & $0.08240$ & $0.08045$ & $2.42$ \\
        $20$ & $96.706$ & $95.894$ & $0.847$ & $28.3754$ & $28.3396$ & $0.126$ & $0.08662$ & $0.08507$ & $1.82$ \\
        $30$ & $84.258$ & $83.484$ & $0.927$ & $23.7057$ & $23.6661$ & $0.167$ & $0.09496$ & $0.09351$ & $1.55$ \\
        $40$ & $68.968$ & $68.293$ & $0.988$ & $18.3350$ & $18.3054$ & $0.162$ & $0.10864$ & $0.10691$ & $1.61$ \\
        $50$ & $52.793$ & $52.252$ & $1.035$ & $13.1087$ & $13.0846$ & $0.184$ & $0.13079$ & $0.12699$ & $2.99$ \\
        $60$ & $37.586$ & $37.232$ & $0.951$ & $8.6459$ & $8.6363$ & $0.111$ & $0.15896$ & $0.15616$ & $1.79$ \\
        $70$ & $24.840$ & $24.615$ & $0.914$ & $5.2711$ & $5.2679$ & $0.061$ & $0.20138$ & $0.19777$ & $1.82$ \\
        $80$ & $15.206$ & $15.044$ & $1.076$ & $2.9828$ & $2.9754$ & $0.248$ & $0.26059$ & $0.25668$ & $1.52$ \\
        $90$ & $8.553$ & $8.472$ & $0.955$ & $1.5665$ & $1.5609$ & $0.359$ & $0.33082$ & $0.34029$ & $-2.78$ \\
        \bottomrule
    \end{tabular}
\end{table}

\section{Robustness and cost of the implicit solve}
\label{sec:var_robustness}

The time constant $c_1$ leaves the converged discrete solution unchanged, and the variable formulation changes it only through the variables that the polynomial space approximates. The three formulations of Table~\ref{tab:formulations} are consistent discretizations of the same problem. Both choices, however, govern whether and how fast the implicit solve converges. Section~\ref{subsec:robustness_formulation} measures the robustness of the three formulations of Table~\ref{tab:formulations} over the $(c_2, c_3)$ plane, and Section~\ref{subsec:c1_study} the cost of the solve, in pseudo-transient and GMRES iterations, as a function of $c_1$.

\subsection{Variable formulation robustness}
\label{subsec:robustness_formulation}

The purpose of this experiment (test case: \case{cyl-ns-robust}) is to quantify how the three variable formulations of Table~\ref{tab:formulations} differ in the robustness of pseudo-transient continuation. Consider the Mach~5 viscous flow over a half-cylinder of Section~\ref{subsec:c4_practice} on the $14 \times 14$, $p = 4$ mesh of Fig.~\ref{fig:halfcyl_mesh}. A grid of 1882 $(c_2, c_3)$ combinations is solved with each formulation, all with $c_4 = 1/c_2$, and each combination is classified as \emph{converged} if PTC reaches the tolerance of Section~\ref{subsec:dg}; \emph{CFL$\to$0} if the CFL number collapses during the high-order solve; \emph{maxit} if the iteration cap of 1000 is reached; or \emph{crash} if the Jacobian assembly produces NaN or Inf entries while the state and residual remain finite. Each combination is first solved with $c_1 = 0.1$ and, if it reaches the iteration cap, is repeated with $c_1 = 0.01$.

\begin{figure}[hbt!]
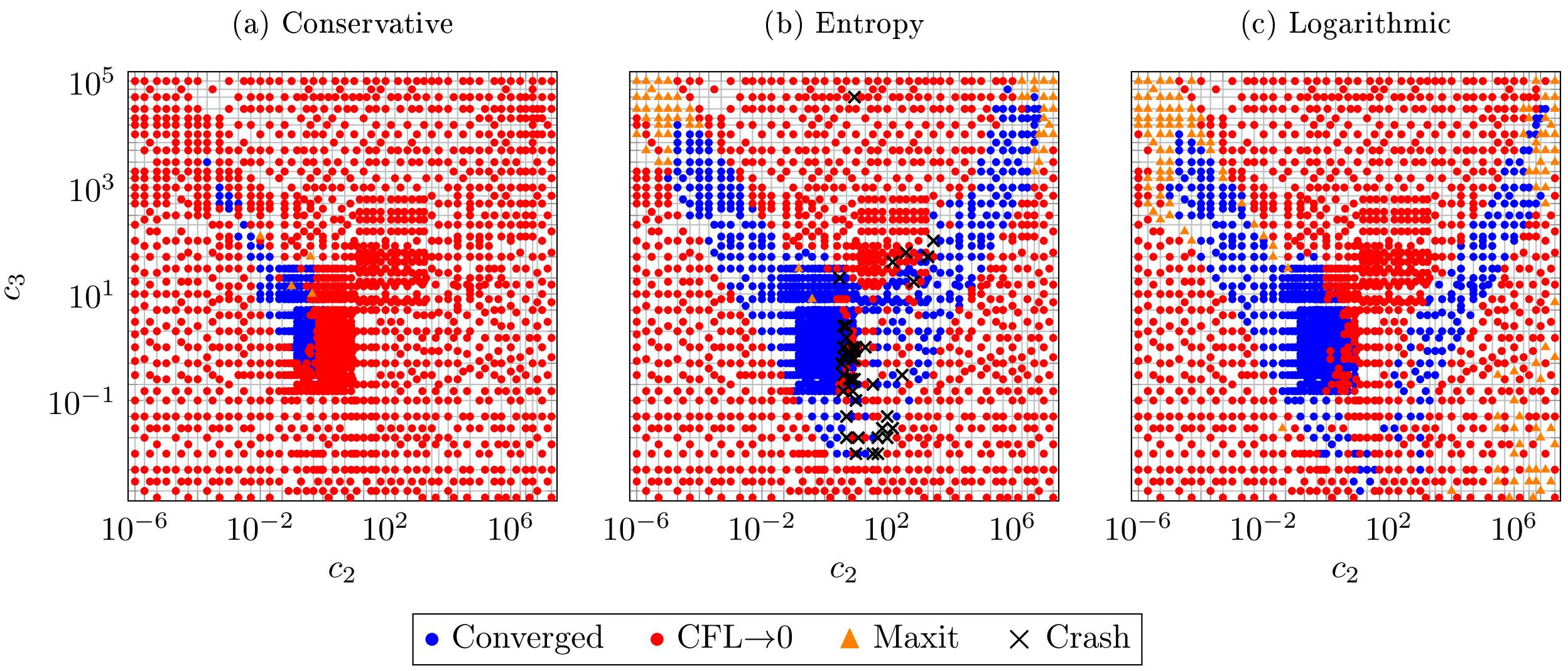

    \centering
    \IfFileExists{_tikz/convmap_conservative.tex}{\input{_tikz/convmap_conservative.tex}}{\fbox{\parbox[c][5cm][c]{0.3\linewidth}{\centering\scriptsize convmap, conservative}}}%
    \hfill
    \IfFileExists{_tikz/convmap_entropy.tex}{\input{_tikz/convmap_entropy.tex}}{\fbox{\parbox[c][5cm][c]{0.3\linewidth}{\centering\scriptsize convmap, entropy}}}%
    \hfill
    \IfFileExists{_tikz/convmap_logarithmic.tex}{\input{_tikz/convmap_logarithmic.tex}}{\fbox{\parbox[c][5cm][c]{0.3\linewidth}{\centering\scriptsize convmap, logarithmic}}}\\[0.5em]
    \fbox{\textcolor{blue}{$\bullet$}~Converged \quad \textcolor{red}{$\bullet$}~CFL$\to$0 \quad \textcolor{orange}{$\blacktriangle$}~Maxit \quad $\boldsymbol{\times}$~Crash}
    \caption{\case{cyl-ns-robust}: convergence outcome of each $(c_2, c_3)$ pair for the (a)~conservative, (b)~entropy, and (c)~logarithmic formulations; both axes logarithmic.}
    \label{fig:convmaps}
\end{figure}

Entropy variables converge on the largest fraction of the grid ($39.3\%$), followed by logarithmic ($34.7\%$) and conservative ($9.2\%$) (Fig.~\ref{fig:convmaps}, Table~\ref{tab:robustness_all}). Conservative variables fail almost entirely by CFL$\to$0, with convergence restricted to a narrow low-$c_2$ strip that collapses for $c_2 \gtrsim 1$. Logarithmic variables converge on a substantially wider window, and their only failures are CFL$\to$0 and maxit. Entropy variables are the most robust but produce 44 crashes at $c_2 \gtrsim 3$. In each crash, oscillations of the high-order approximation drive the last entropy variable, $-\rho/P$, or in one case the pressure, toward zero at a quadrature point. The state remains finite there, but the derivative of the change of variables, which contains inverse powers of these quantities, overflows. This singularity does not occur in conservative or logarithmic variables. Despite these crashes, entropy variables provide the widest usable $(c_2, c_3)$ window on this case, which is why they are the default formulation in Sections~\ref{sec:avpde_design}--\ref{sec:applications}.

\begin{table}[hbt!]
\centering
\caption{\case{cyl-ns-robust}: convergence statistics for the three variable formulations across the $(c_2, c_3)$ parameter sweep.}
\label{tab:robustness_all}
\begin{tabular}{lccccc}
\toprule
Formulation & Converged & CFL$\to$0 & Maxit & Crash & Rate \\
\midrule
Conservative & 173/1882 & 1705 & 4 & 0 & 9.2\% \\
Entropy & 739/1882 & 1045 & 54 & 44 & 39.3\% \\
Logarithmic & 653/1882 & 1104 & 125 & 0 & 34.7\% \\
\bottomrule
\end{tabular}
\end{table}

\begin{remark}
The converged region is V-shaped in the $(c_2, c_3)$ plane, and its two arms correspond to two kinds of AV (Fig.~\ref{fig:convmaps}). At small $c_2$ the AV stays localized at the bow shock and the shock remains sharp. In this localized regime the solve fails on both sides of the arm. When $c_3$ is too small, the AV is too weak to suppress the nonphysical oscillations near the shock. When $c_3$ is too large, the artificial viscous flux dominates the inviscid flux and the governing equations become diffusion dominated, which would require boundary conditions posed for a diffusive rather than an advective problem. At large $c_2$ the AV delocalizes into a smooth, nearly domain-filling field that smears the bow shock, and in the most diffusive cases the Mach number falls below its freestream value ahead of the shock. This broadband regime is nonphysical, but the solver tends to converge in it because the excessive, smoothly distributed dissipation keeps the Jacobian well conditioned. The right arm of the V therefore marks convergence of the solver, not accuracy of the solution.
\end{remark}

\begin{remark}
The robustness of the implicit solve also depends on the boundary treatment, in particular on the outlet condition where the subsonic outlet meets the no-slip wall and the boundary layer leaves the domain. All three variable formulations in this section use the same boundary treatment, so the comparison between them is consistent, but the sensitivity of the convergence window to the boundary treatment itself is not studied here and is left to future work.
\end{remark}

\subsection{Effect of the time constant \texorpdfstring{$c_1$}{c1}}
\label{subsec:c1_study}

\begin{figure}[hbt!]
    \centering
    \IfFileExists{_tikz/325s_byform_no042.tex}{\begin{tikzpicture}
\definecolor{set1}{RGB}{0,114,178}    %
\definecolor{set2}{RGB}{230,159,0}    %
\definecolor{set3}{RGB}{0,158,115}    %
\definecolor{set4}{RGB}{204,121,167}  %
\definecolor{set5}{RGB}{213,94,0}     %
\tikzset{
    cset1/.style={set1, solid,        line width=0.9pt, mark=*,        mark size=1.4pt, mark options={solid, fill=set1}},
    cset2/.style={set2, dashed,       line width=0.9pt, mark=square*,  mark size=1.4pt, mark options={solid, fill=set2}},
    cset3/.style={set3, dotted,       line width=1.0pt, mark=triangle*,mark size=1.7pt, mark options={solid, fill=set3}},
    cset4/.style={set4, dash dot,     line width=0.9pt, mark=diamond*, mark size=1.6pt, mark options={solid, fill=set4}},
    cset5/.style={set5, dash dot dot, line width=0.9pt, mark=pentagon*,mark size=1.7pt, mark options={solid, fill=set5}},
}
\begin{groupplot}[
    group style={group size=3 by 2, horizontal sep=0.9cm, vertical sep=1.6cm,
        x descriptions at=edge bottom, y descriptions at=edge left},
    width=5.9cm,
    height=5.2cm,
    xmode=log,
    xmin=6e-5, xmax=1.6, xtick={1e-4,1e-3,1e-2,1e-1,1},
    ymode=normal,
    xlabel={$c_1$},
    grid=both,
    tick label style={font=\normalsize},
]
\nextgroupplot[title={(a) Conservative}, ylabel={PTC iterations}, ymin=0, ymax=451.44000000000005, ytick={0,150,300,450},
    legend to name=leg325byform,
    legend columns=-1,
    legend style={draw=black, fill=white, /tikz/every even column/.append style={column sep=0.5cm}},
    legend cell align=left,
]
\addplot[cset1] coordinates {
    (0.1, 150)
    (0.31622776601683794, 122)
};
\addlegendentry{$(0.15,\,0.55)$}
\addplot[cset2] coordinates {
    (0.1, 63)
};
\addlegendentry{$(0.2,\,0.75)$}
\addplot[cset3] coordinates {
    (0.1, 52)
};
\addlegendentry{$(0.3,\,1.25)$}
\addplot[cset4] coordinates {
    (0.1, 49)
};
\addlegendentry{$(0.5,\,3)$}
\addplot[cset5] coordinates {
    (0.1, 54)
};
\addlegendentry{$(0.65,\,4)$}
\nextgroupplot[title={(b) Entropy}, ymin=0, ymax=451.44000000000005, ytick={0,150,300,450}]
\addplot[cset1] coordinates {
    (0.0001, 218)
    (0.00031622776601683794, 123)
    (0.001, 71)
    (0.0031622776601683794, 68)
    (0.010000000000000002, 64)
    (0.03162277660168379, 52)
    (0.1, 49)
    (0.31622776601683794, 41)
    (1.0, 37)
};
\addplot[cset2] coordinates {
    (0.0001, 179)
    (0.00031622776601683794, 98)
    (0.001, 52)
    (0.0031622776601683794, 50)
    (0.010000000000000002, 44)
    (0.03162277660168379, 46)
    (0.1, 43)
    (0.31622776601683794, 37)
    (1.0, 31)
};
\addplot[cset3] coordinates {
    (0.0001, 222)
    (0.00031622776601683794, 114)
    (0.001, 102)
    (0.0031622776601683794, 42)
    (0.010000000000000002, 41)
    (0.03162277660168379, 40)
    (0.1, 59)
    (0.31622776601683794, 32)
};
\addplot[cset4] coordinates {
    (0.0001, 357)
    (0.00031622776601683794, 192)
    (0.001, 142)
    (0.0031622776601683794, 48)
    (0.010000000000000002, 52)
    (0.03162277660168379, 29)
    (0.1, 61)
    (0.31622776601683794, 43)
};
\addplot[cset5] coordinates {
    (0.0001, 418)
    (0.00031622776601683794, 202)
    (0.001, 136)
    (0.0031622776601683794, 79)
    (0.010000000000000002, 49)
    (0.03162277660168379, 67)
    (0.1, 38)
    (0.31622776601683794, 42)
};
\nextgroupplot[title={(c) Logarithmic}, ymin=0, ymax=451.44000000000005, ytick={0,150,300,450}]
\addplot[cset1] coordinates {
    (0.0001, 239)
    (0.00031622776601683794, 139)
    (0.001, 90)
    (0.0031622776601683794, 74)
    (0.010000000000000002, 70)
    (0.03162277660168379, 60)
    (0.1, 57)
    (0.31622776601683794, 48)
};
\addplot[cset2] coordinates {
    (0.0001, 205)
    (0.00031622776601683794, 127)
    (0.001, 84)
    (0.0031622776601683794, 61)
    (0.010000000000000002, 73)
    (0.03162277660168379, 67)
    (0.1, 62)
    (0.31622776601683794, 45)
};
\addplot[cset3] coordinates {
    (0.0001, 231)
    (0.00031622776601683794, 125)
    (0.001, 75)
    (0.0031622776601683794, 47)
    (0.010000000000000002, 42)
    (0.03162277660168379, 40)
    (0.1, 37)
    (0.31622776601683794, 40)
};
\addplot[cset4] coordinates {
    (0.0001, 323)
    (0.00031622776601683794, 166)
    (0.001, 92)
    (0.0031622776601683794, 55)
    (0.010000000000000002, 39)
    (0.03162277660168379, 32)
    (0.1, 29)
};
\addplot[cset5] coordinates {
    (0.0001, 359)
    (0.00031622776601683794, 197)
    (0.001, 97)
    (0.0031622776601683794, 64)
    (0.010000000000000002, 41)
    (0.03162277660168379, 32)
    (0.1, 28)
};
\nextgroupplot[title={(d)}, ylabel={Mean GMRES iters / PTC step}, ymin=10, ymax=29.063880000000005]
\addplot[cset1] coordinates {
    (0.1, 16.793)
    (0.31622776601683794, 17.738)
};
\addplot[cset2] coordinates {
    (0.1, 18.032)
};
\addplot[cset3] coordinates {
    (0.1, 19.135)
};
\addplot[cset4] coordinates {
    (0.1, 18.898)
};
\addplot[cset5] coordinates {
    (0.1, 20.019)
};
\nextgroupplot[title={(e)}, ymin=10, ymax=29.063880000000005]
\addplot[cset1] coordinates {
    (0.0001, 15.839)
    (0.00031622776601683794, 15.724)
    (0.001, 14.296)
    (0.0031622776601683794, 13.838)
    (0.010000000000000002, 14.344)
    (0.03162277660168379, 15.288)
    (0.1, 15.51)
    (0.31622776601683794, 17.098)
    (1.0, 18.405)
};
\addplot[cset2] coordinates {
    (0.0001, 17.168)
    (0.00031622776601683794, 16.847)
    (0.001, 17.077)
    (0.0031622776601683794, 16.24)
    (0.010000000000000002, 16.205)
    (0.03162277660168379, 16.957)
    (0.1, 16.581)
    (0.31622776601683794, 17.784)
    (1.0, 20.129)
};
\addplot[cset3] coordinates {
    (0.0001, 18.495)
    (0.00031622776601683794, 19.193)
    (0.001, 18.755)
    (0.0031622776601683794, 18.167)
    (0.010000000000000002, 17.488)
    (0.03162277660168379, 17.85)
    (0.1, 16.237)
    (0.31622776601683794, 19.469)
};
\addplot[cset4] coordinates {
    (0.0001, 19.317)
    (0.00031622776601683794, 20.156)
    (0.001, 22.261)
    (0.0031622776601683794, 23.0)
    (0.010000000000000002, 23.519)
    (0.03162277660168379, 20.31)
    (0.1, 23.295)
    (0.31622776601683794, 22.047)
};
\addplot[cset5] coordinates {
    (0.0001, 21.718)
    (0.00031622776601683794, 23.05)
    (0.001, 24.537)
    (0.0031622776601683794, 26.911)
    (0.010000000000000002, 25.306)
    (0.03162277660168379, 26.537)
    (0.1, 25.5)
    (0.31622776601683794, 26.0)
};
\nextgroupplot[title={(f)}, ymin=10, ymax=29.063880000000005]
\addplot[cset1] coordinates {
    (0.0001, 16.636)
    (0.00031622776601683794, 17.216)
    (0.001, 16.178)
    (0.0031622776601683794, 16.784)
    (0.010000000000000002, 16.4)
    (0.03162277660168379, 16.817)
    (0.1, 17.088)
    (0.31622776601683794, 18.854)
};
\addplot[cset2] coordinates {
    (0.0001, 16.995)
    (0.00031622776601683794, 17.661)
    (0.001, 18.226)
    (0.0031622776601683794, 18.721)
    (0.010000000000000002, 18.438)
    (0.03162277660168379, 16.866)
    (0.1, 19.274)
    (0.31622776601683794, 19.022)
};
\addplot[cset3] coordinates {
    (0.0001, 18.147)
    (0.00031622776601683794, 18.344)
    (0.001, 19.547)
    (0.0031622776601683794, 16.83)
    (0.010000000000000002, 18.357)
    (0.03162277660168379, 17.75)
    (0.1, 18.865)
    (0.31622776601683794, 21.475)
};
\addplot[cset4] coordinates {
    (0.0001, 19.892)
    (0.00031622776601683794, 21.09)
    (0.001, 21.565)
    (0.0031622776601683794, 22.218)
    (0.010000000000000002, 21.436)
    (0.03162277660168379, 20.406)
    (0.1, 19.793)
};
\addplot[cset5] coordinates {
    (0.0001, 20.989)
    (0.00031622776601683794, 22.36)
    (0.001, 24.804)
    (0.0031622776601683794, 24.625)
    (0.010000000000000002, 24.854)
    (0.03162277660168379, 23.531)
    (0.1, 21.179)
};
\end{groupplot}
\coordinate (leg) at ($(group c1r2.south)!0.5!(group c3r2.south) - (0,1.1cm)$);
\node[anchor=north] at (leg) {\ref*{leg325byform}};
\end{tikzpicture}}{\fbox{\parbox[c][5cm][c]{0.6\linewidth}{\centering\scriptsize 325s byform}}}
    \caption{\case{cyl-ns-robust}: (a--c)~PTC iterations and (d--f)~mean GMRES iterations per PTC step against $c_1$ for the (a,d)~conservative, (b,e)~entropy, and (c,f)~logarithmic formulations, one curve per $(c_2, c_3)$ setting. Only converged combinations are shown.}
    \label{fig:c1_byform_no042}
\end{figure}
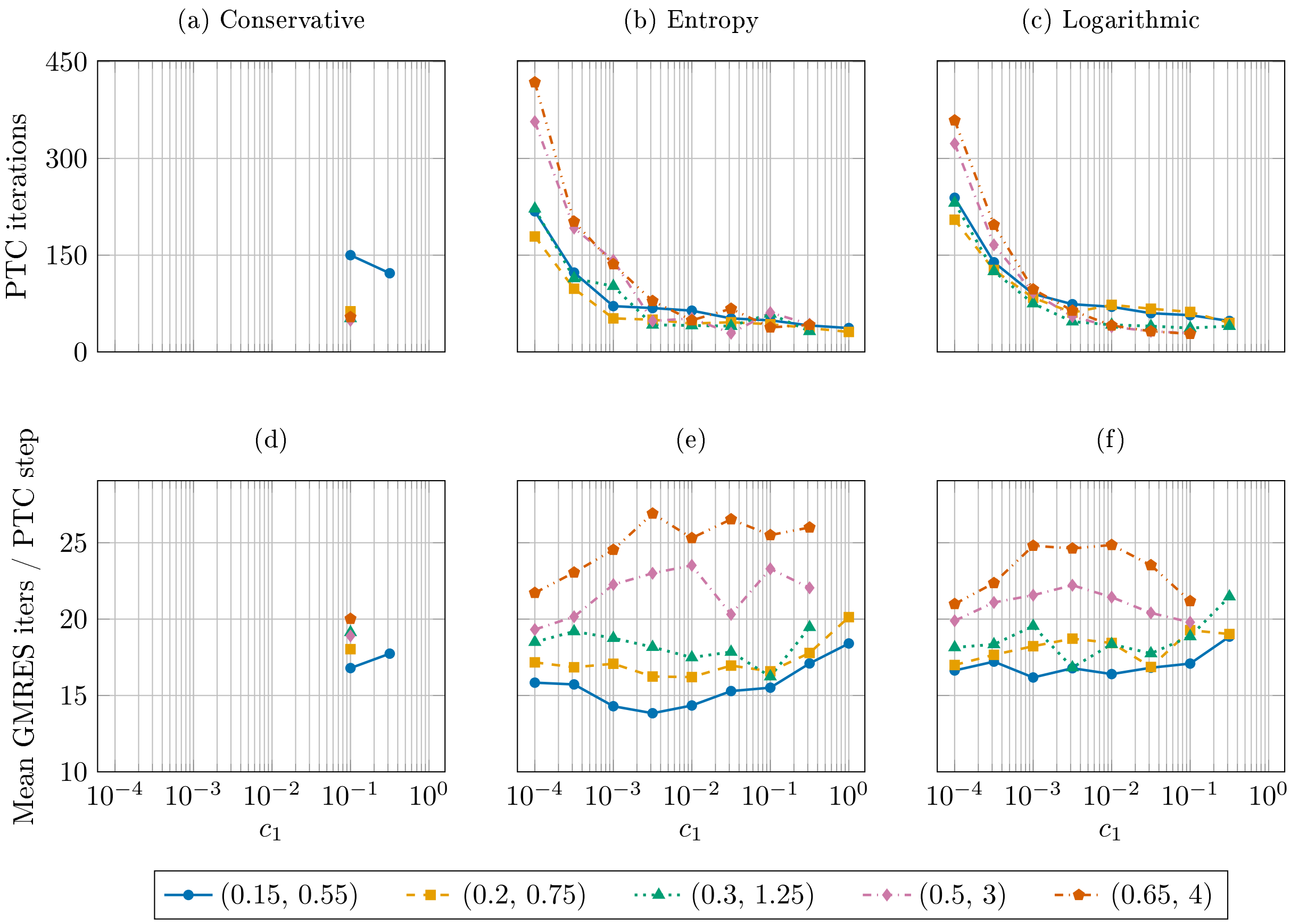

On the same case, a $c_1$ sweep is performed for the three formulations at five $(c_2, c_3)$ settings of increasing AV, $(c_2, c_3) \in \{(0.15, 0.55), (0.2, 0.75), (0.3, 1.25), (0.5, 3.0), (0.65, 4.0)\}$, chosen so that both coefficients increase from one setting to the next. The purpose is to measure how $c_1$ trades robustness against the cost of the solve. The conservative formulation converges only in a narrow $c_1$ window near $0.1$, collapsing by CFL$\to$0 at both smaller and larger $c_1$, and requires more PTC iterations where it does converge, whereas the entropy and logarithmic formulations converge across most of the range (Fig.~\ref{fig:c1_byform_no042}). Within each formulation's window the PTC count drops as $c_1$ increases, with too large a $c_1$ eventually stalling, and the mean GMRES iterations per step grow with the AV magnitude, the larger-$(c_2, c_3)$ settings sitting above the smaller ones in every panel.

\begin{remark}
The time constant acts as a damping coefficient on the AV, whose relaxation time is $c_1/c_4$. A larger $c_1$ slows the AV's response to transients in the state, so the AV does not feed the oscillations of the pseudo-transient iterates back into the flow through a rapidly changing artificial viscous flux; fewer steps are rejected and the CFL grows. Too large a $c_1$ leaves the AV lagging the state, the shock is not stabilized in time, and the iteration stalls. This is the robustness-versus-cost trade-off noted in Section~\ref{subsec:dg}. The growth of the GMRES count with the AV magnitude reflects the preconditioner: the block ILU(0) preconditioner of Section~\ref{subsec:dg} suits the advection-dominated Jacobian of the flow but is less effective once the AV adds a strong diffusion-dominated term.
\end{remark}

\section{Conclusions}
\label{sec:conclude}

This work made three contributions to the use of PDE-based artificial viscosity in high-order DG methods, each validated on a hierarchy of flows of increasing complexity. First, an $L^2$-projection analysis identified the smallest shock thickness the polynomial space can monotonically resolve as scaling with the resolution length $\bar h = h/p$, and the AV PDE was scaled so that its solution is self-similar in $\bar h$. The result is a captured shock that thins with the resolution length rather than persisting at a fixed physical width. Under joint $hp$-refinement the rescaled AV and state profiles collapsed onto a single curve on every test case, from the 1D Burgers' equation to the Mach~5 sphere. On the double wedge and the sphere, where the captured shock was measured as an element count, it spanned $2$--$4$ elements, unchanged to within $6\%$ as the mesh was halved. Second, a Green's-function analysis of the steady AV PDE exposed the width--magnitude coupling of the standard parametrization, in which widening the AV by raising $c_2$ also lowers its peak; the reparametrization $c_4 = 1/c_2$ removes this coupling in the Green's-function response and inverts it in the coupled state--AV problem. The practical consequence, demonstrated on the viscous flow over a half-cylinder, is that widening under the derived parametrization raised the AV peak and cut the captured-shock Mach overshoot by up to $91\%$, whereas the same widening under the standard parametrization lowered the peak and reduced the overshoot by at most $10\%$. Third, a study of how the variable formulation and the parameters $c_1$, $c_2$, and $c_3$ govern the robustness and cost of the implicit solve showed that variable choice dominates robustness. Entropy variables converged on $39.3\%$ of the $(c_2, c_3)$ space tested, against $34.7\%$ for logarithmic and only $9.2\%$ for conservative. The time constant $c_1$ traded robustness against per-solve cost, and the conservative formulation converged only in a narrow $c_1$ window.

Future work will pursue two directions. The first is to integrate PDE-based shock capturing with implicit shock tracking~\cite{huang2022robust, dong2025rpadapt}, so that shock fitting resolves the strong bow shock while shock capturing handles the weaker and more complex secondary shocks. The second is to select the coefficients $c_2$ and $c_3$, which are still chosen by hand, by PDE-constrained optimization, for example minimization of a dual-weighted-residual error estimate for a given quantity of interest.

\section*{Acknowledgments}
This work is supported by AFOSR award numbers FA9550-20-1-0236,
FA9550-22-1-0002, FA9550-22-1-0004, ONR award number
N00014-22-1-2299, and NSF award number CBET-2338843. A. H. Razavi and M. Yano acknowledge financial support by the Natural Sciences and Engineering Research Council of Canada.
The content of this publication does not necessarily reflect the position
or policy of any of these supporters, and no official endorsement should
be inferred.

\section*{Declaration of generative AI and AI-assisted technologies in the manuscript preparation process}
During the preparation of this work, the authors used Claude Fable 5.1
(Anthropic), accessed through Claude Code, in order to revise and reorganize
the text of the manuscript for clarity and conciseness. After using this tool, the
authors reviewed and edited the content as needed and take full responsibility
for the content of the published article.

\bibliographystyle{elsarticle-num}
\bibliography{biblio}

\begin{thebibliography}{10}
\expandafter\ifx\csname url\endcsname\relax
  \def\url#1{\texttt{#1}}\fi
\expandafter\ifx\csname urlprefix\endcsname\relax\def\urlprefix{URL }\fi
\expandafter\ifx\csname href\endcsname\relax
  \def\href#1#2{#2} \def\path#1{#1}\fi

\bibitem{hesthaven2007nodal}
J.~S. Hesthaven, T.~Warburton, Nodal Discontinuous {G}alerkin Methods:
  Algorithms, Analysis, and Applications, Springer, New York, 2008.
\newblock \href {https://doi.org/10.1007/978-0-387-72067-8}
  {\path{doi:10.1007/978-0-387-72067-8}}.

\bibitem{van1979towards}
B.~Van~Leer, Towards the ultimate conservative difference scheme. {V}. a
  second-order sequel to {G}odunov's method, Journal of Computational Physics
  32~(1) (1979) 101--136.

\bibitem{harten1997uniformly}
A.~Harten, B.~Engquist, S.~Osher, S.~R. Chakravarthy, Uniformly high order
  accurate essentially non-oscillatory schemes, {III}, Journal of Computational
  Physics 131~(1) (1997) 3--47.

\bibitem{liu1994weighted}
X.-D. Liu, S.~Osher, T.~Chan, Weighted essentially non-oscillatory schemes,
  Journal of Computational Physics 115~(1) (1994) 200--212.

\bibitem{jiang1996efficient}
G.-S. Jiang, C.-W. Shu, Efficient implementation of weighted {ENO} schemes,
  Journal of Computational Physics 126~(1) (1996) 202--228.

\bibitem{cockburn2000development}
B.~Cockburn, G.~E. Karniadakis, C.-W. Shu, The development of discontinuous
  {G}alerkin methods, in: Discontinuous {G}alerkin Methods: Theory, Computation
  and Applications, Springer, 2000, pp. 3--50.
\newblock \href {https://doi.org/10.1007/978-3-642-59721-3_1}
  {\path{doi:10.1007/978-3-642-59721-3_1}}.

\bibitem{dumbser2014posteriori}
M.~Dumbser, O.~Zanotti, R.~Loub{\`e}re, S.~Diot, A posteriori subcell limiting
  of the discontinuous {G}alerkin finite element method for hyperbolic
  conservation laws, Journal of Computational Physics 278 (2014) 47--75.
\newblock \href {https://doi.org/10.1016/j.jcp.2014.08.009}
  {\path{doi:10.1016/j.jcp.2014.08.009}}.

\bibitem{vilar2019posteriori}
F.~Vilar, A posteriori local subcell correction of high-order discontinuous
  {G}alerkin scheme for conservation laws on two-dimensional unstructured
  grids, Journal of Computational PhysicsPreprint (2019).

\bibitem{persson2006sub}
P.-O. Persson, J.~Peraire, Sub-cell shock capturing for discontinuous
  {G}alerkin methods, in: 44th AIAA Aerospace Sciences Meeting and Exhibit,
  Reno, Nevada, 2006, {AIAA} Paper 2006-112.

\bibitem{barter2010shock}
G.~E. Barter, D.~L. Darmofal, Shock capturing with {PDE}-based artificial
  viscosity for {DGFEM}: {P}art {I}. {F}ormulation, Journal of Computational
  Physics 229~(5) (2010) 1810--1827.

\bibitem{persson2013shock}
P.-O. Persson, Shock capturing for high-order discontinuous {G}alerkin
  simulation of transient flow problems, in: 21st AIAA Computational Fluid
  Dynamics Conference, 2013, {AIAA} Paper 2013-3061.

\bibitem{moro2016dilation}
D.~Moro, N.~C. Nguyen, J.~Peraire, Dilation-based shock capturing for
  high-order methods, International Journal for Numerical Methods in Fluids
  82~(7) (2016) 398--416.

\bibitem{fernandez2018physics}
P.~Fern{\'a}ndez, N.~C. Nguyen, J.~Peraire, A physics-based shock capturing
  method for unsteady laminar and turbulent flows, in: 2018 AIAA Aerospace
  Sciences Meeting, Kissimmee, Florida, 2018, {AIAA} Paper 2018-0062.

\bibitem{ching2019shock}
E.~J. Ching, Y.~Lv, P.~Gnoffo, M.~Barnhardt, M.~Ihme, Shock capturing for
  discontinuous {G}alerkin methods with application to predicting heat transfer
  in hypersonic flows, Journal of Computational Physics 376 (2019) 54--75.

\bibitem{moretti2002thirty}
G.~Moretti, Thirty-six years of shock fitting, Computers \& Fluids 31~(4-7)
  (2002) 719--723.

\bibitem{salas2009shock}
M.~D. Salas, A Shock-Fitting Primer, CRC Press, 2009.

\bibitem{bonfiglioli2016unsteady}
A.~Bonfiglioli, R.~Paciorri, L.~Campoli, Unsteady shock-fitting for
  unstructured grids, International Journal for Numerical Methods in Fluids
  81~(4) (2016) 245--261.

\bibitem{zou2021embedded}
D.~Zou, A.~Bonfiglioli, R.~Paciorri, J.~Liu, An embedded shock-fitting
  technique on unstructured dynamic grids, Computers \& Fluids 218 (2021)
  104847.

\bibitem{assonitis2021numerical}
A.~Assonitis, R.~Paciorri, A.~Bonfiglioli, Numerical simulation of
  shock/boundary-layer interaction using an unstructured shock-fitting
  technique, Computers \& Fluids 228 (2021) 105058.

\bibitem{corrigan2019moving}
A.~Corrigan, A.~D. Kercher, D.~A. Kessler, A moving discontinuous {G}alerkin
  finite element method for flows with interfaces, International Journal for
  Numerical Methods in Fluids 89~(9) (2019) 362--406.

\bibitem{kercher2021least}
A.~D. Kercher, A.~Corrigan, A least-squares formulation of the moving
  discontinuous {G}alerkin finite element method with interface condition
  enforcement, Computers \& Mathematics with Applications 95 (2021) 143--171.

\bibitem{zahr2018shktrk}
M.~J. Zahr, P.-O. Persson,
  \href{http://www.sciencedirect.com/science/article/pii/S002199911830189X}{An
  optimization-based approach for high-order accurate discretization of
  conservation laws with discontinuous solutions}, Journal of Computational
  Physics 365 (2018) 105 -- 134.
\newblock \href {https://doi.org/https://doi.org/10.1016/j.jcp.2018.03.029}
  {\path{doi:https://doi.org/10.1016/j.jcp.2018.03.029}}.
\newline\urlprefix\url{http://www.sciencedirect.com/science/article/pii/S002199911830189X}

\bibitem{zahr2020implicit}
M.~J. Zahr, A.~Shi, P.-O. Persson, Implicit shock tracking using an
  optimization-based high-order discontinuous {G}alerkin method, Journal of
  Computational Physics 410 (2020) 109385.

\bibitem{huang2022robust}
T.~Huang, M.~J. Zahr, A robust, high-order implicit shock tracking method for
  simulation of complex, high-speed flows, Journal of Computational Physics 454
  (2022) 110981.

\bibitem{shi2022implicit}
A.~Shi, P.-O. Persson, M.~J. Zahr, Implicit shock tracking for unsteady flows
  by the method of lines, Journal of Computational Physics 454 (2022) 110906.

\bibitem{naudet2024space}
C.~J. Naudet, M.~J. Zahr, A space-time high-order implicit shock tracking
  method for shock-dominated unsteady flows, Journal of Computational Physics
  501 (2024) 112792.

\bibitem{yano2012optimization}
M.~Yano, An optimization framework for adaptive higher-order discretizations of
  partial differential equations on anisotropic simplex meshes, Ph.D. thesis,
  Massachusetts Institute of Technology (2012).

\bibitem{nguyen2024optimal}
N.~C. Nguyen, R.~L. Van~Heyningen, J.~Vila-P{\'e}rez, J.~Peraire, Optimal
  transport for mesh adaptivity and shock capturing of compressible flows,
  Journal of Computational Physics 508 (2024) 113005.
\newblock \href {https://doi.org/10.1016/j.jcp.2024.113005}
  {\path{doi:10.1016/j.jcp.2024.113005}}.

\bibitem{razavi2026output}
A.~H. Razavi, A.~Ambardekar, M.~J. Zahr, M.~Yano, Output error estimation and
  adaptive mesh refinement for {DGFEM} with {PDE}-based artificial viscosity
  for shock capturing, Journal of Computational Physics (2026) 115514\href
  {https://doi.org/10.1016/j.jcp.2026.115514}
  {\path{doi:10.1016/j.jcp.2026.115514}}.

\bibitem{alberti2024comparative}
L.~Alberti, F.~Bassi, E.~Carnevali, A.~Colombo, A.~Crivellini, A.~Nigro, A
  comparative study of different sets of variables in a discontinuous
  {Galerkin} method with entropy balance enforcement, International Journal of
  Computational Fluid Dynamics 37~(6) (2023) 487--508.
\newblock \href {https://doi.org/10.1080/10618562.2024.2310537}
  {\path{doi:10.1080/10618562.2024.2310537}}.

\bibitem{hughes1986new}
T.~J.~R. Hughes, L.~P. Franca, M.~Mallet, A new finite element formulation for
  computational fluid dynamics: {I}. {S}ymmetric forms of the compressible
  {E}uler and {N}avier--{S}tokes equations and the second law of
  thermodynamics, Computer Methods in Applied Mechanics and Engineering 54~(2)
  (1986) 223--234.
\newblock \href {https://doi.org/10.1016/0045-7825(86)90127-1}
  {\path{doi:10.1016/0045-7825(86)90127-1}}.

\bibitem{hartmann2008optimal}
R.~Hartmann, P.~Houston, An optimal order interior penalty discontinuous
  {G}alerkin discretization of the compressible {N}avier--{S}tokes equations,
  Journal of Computational Physics 227~(22) (2008) 9670--9685.

\bibitem{roe1981approximate}
P.~L. Roe, Approximate {R}iemann solvers, parameter vectors, and difference
  schemes, Journal of Computational Physics 43~(2) (1981) 357--372.
\newblock \href {https://doi.org/10.1016/0021-9991(81)90128-5}
  {\path{doi:10.1016/0021-9991(81)90128-5}}.

\bibitem{harten1983self}
A.~Harten, J.~M. Hyman, Self adjusting grid methods for one-dimensional
  hyperbolic conservation laws, Journal of Computational Physics 50~(2) (1983)
  235--269.
\newblock \href {https://doi.org/10.1016/0021-9991(83)90066-9}
  {\path{doi:10.1016/0021-9991(83)90066-9}}.

\bibitem{kelley1998convergence}
C.~T. Kelley, D.~E. Keyes, Convergence analysis of pseudo-transient
  continuation, SIAM Journal on Numerical Analysis 35~(2) (1998) 508--523.

\bibitem{modisette2011automated}
J.~M. Modisette, An automated reliable method for two-dimensional
  {R}eynolds-averaged {N}avier-{S}tokes simulations, Ph.D. thesis,
  Massachusetts Institute of Technology (2011).

\bibitem{persson_newton-gmres_2008}
P.-O. Persson, J.~Peraire,
  \href{http://epubs.siam.org/doi/10.1137/070692108}{Newton-{GMRES}
  preconditioning for discontinuous {Galerkin} discretizations of the
  {Navier}--{Stokes} equations}, SIAM Journal on Scientific Computing 30~(6)
  (2008) 2709--2733.
\newblock \href {https://doi.org/10.1137/070692108}
  {\path{doi:10.1137/070692108}}.
\newline\urlprefix\url{http://epubs.siam.org/doi/10.1137/070692108}

\bibitem{whitham1974linear}
G.~B. Whitham, Linear and Nonlinear Waves, John Wiley \& Sons, New York, 1974.

\bibitem{dong2025rpadapt}
H.~Dong, M.~Yano, T.~Huang, M.~J. Zahr,
  \href{https://www.sciencedirect.com/science/article/pii/S0021999125006199}{An
  $rp$-adaptive method for accurate resolution of shock-dominated viscous flow
  based on implicit shock tracking}, Journal of Computational Physics 543
  (2025) 114337.
\newblock \href {https://doi.org/https://doi.org/10.1016/j.jcp.2025.114337}
  {\path{doi:https://doi.org/10.1016/j.jcp.2025.114337}}.
\newline\urlprefix\url{https://www.sciencedirect.com/science/article/pii/S0021999125006199}

\bibitem{blottner1989}
F.~G. Blottner, Accurate {Navier}--{Stokes} results for the hypersonic flow
  over a spherical nosetip, in: AIAA 24th Thermophysics Conference, Buffalo,
  New York, 1989, {AIAA} Paper 89-1671.

\bibitem{widhopf1972}
G.~F. Widhopf, Laminar, transitional, and turbulent heat transfer measurements
  on a yawed blunt conical nosetip, Tech. Rep. TR-0073(S3450-60)-1, The
  Aerospace Corporation, El Segundo, California, also SAMSO-TR-72-209; DTIC
  AD-748292 (September 1972).

\bibitem{fisher2021blottner}
T.~Fisher,
  \href{https://turbmodels.larc.nasa.gov/highfidelitycfd_workshop2022.html}{High
  fidelity {CFD} workshop 2021: High speed steady advanced case: {Blottner}
  sphere}, NASA Langley Research Center Turbulence Modeling Resource (2021).
\newline\urlprefix\url{https://turbmodels.larc.nasa.gov/highfidelitycfd_workshop2022.html}

\end{thebibliography}

\end{document}